\documentclass[12pt]{amsart}

\usepackage{amsmath}
\usepackage{amssymb}
\usepackage{graphicx}
\usepackage[a4paper, margin=1in]{geometry}
\usepackage{float} 
\usepackage{bm}
\usepackage{subcaption}
\usepackage{xcolor}
\usepackage{booktabs}
\usepackage{array} 
\usepackage{indentfirst}

\usepackage[colorlinks=true, linkcolor=blue, citecolor=blue, urlcolor=blue]{hyperref}

\usepackage{cleveref}

\makeatletter
\renewcommand{\tocsection}[3]{%
  \vspace{0.8ex} 
  \indentlabel{\@ifnotempty{#2}{\bfseries\ignorespaces#1 #2.\quad}}\bfseries#3}

\renewcommand{\tocsubsection}[3]{%
  \indentlabel{\@ifnotempty{#2}{\mdseries\hspace{2.5em}\ignorespaces#1 #2.\quad}}\mdseries#3}
\makeatother

\renewcommand{\vec}[1]{\mathbf{#1}}
\newcommand{\tens}[1]{\boldsymbol{#1}}

\title[Numerical PDEs Approach to Shape Evolution]{A Numerical PDEs Approach to Evolution Equations in Shape Analysis Based on Regularized Morphoelasticity}

\author{Ziqin Zhou}
\address{Department of Applied Mathematics and Statistics, Johns Hopkins University.}
\email{zhouzq4@gmail.com}
\date{}

\begin{document}
	
	\maketitle

    \begin{abstract}
    \normalsize
    This work studies a variational formulation and numerical solution of a regularized morphoelasticity problem of shape evolution. The foundation of our analysis is based on the governing equations of linear elasticity, extended to account for volumetric growth. In the morphoelastic framework, the total deformation is decomposed into an elastic component and a growth component, represented by a growth tensor $\mathbf{G}$. While the forward one-step problem—computing displacement given a growth tensor—is well-established, a more challenging and relevant question in biological modeling is the inverse problem in a continuous sense. While this problem is fundamentally ill-posed without additional constraints, we will explore parametrized growth models inscribed within an optimal control problem inspired by the Large Deformation Diffeomorphic Metric Mapping (LDDMM) framework. By treating the growth process as a path within a shape space, we can define a physically meaningful metric and seek the most plausible, energy-efficient trajectory between configurations. In the construction, a high-order regularization term is introduced. This elevates the governing equations to a high-order elliptic system, ensuring the existence of a smooth solution. This dissertation focuses on the issue of solving this equation efficiently, as this is a key requirement for the feasibility of the overall approach. This will be achieved with the help of finite element solvers, notably from the FEniCSx library in Python. Also, we implement a Mixed Finite Element Method, which decomposes the problem into a system of coupled second-order equations as a treatment of these high-order systems that have significant computational challenges.

    \vspace{0.5cm} 

    \noindent \textbf{Keywords:} Morphoelasticity, Shape Analysis, Shape Evolution, Variational Method, Large Deformation Diffeomorphic Metric Mapping, Mixed Finite Element Method.
    \end{abstract}
    
    \clearpage

    \tableofcontents 
    \clearpage 

    \section{Introduction}
    Over the past few decades, advancements in medical imaging technologies (e.g., MRI, CT, and ultrasound) have provided unprecedented access to high-resolution 3D anatomical data. A central challenge in modern computational anatomy is to extract meaningful, reproducible geometric information from these raw scans to quantify anatomical variations, track disease progression, and perform biomechanical simulations \cite{feydy2020geometric}. While the recent deep learning revolution has yielded extraordinary results in texture-based medical image segmentation, standard Convolutional Neural Networks (CNNs) operate on fixed Eulerian grids, pixels, or voxels and inherently lack the mathematical structure to rigorously reason about continuous geometry and topology. Even with the rapid rise of Geometric Deep Learning \cite{bronstein2021geometricdeeplearninggrids}, capturing the complex symmetries and large deformations of anatomical structures remains a significant challenge. To accurately capture large, free-form anatomical deformations and guarantee robustness to topological noise, it is strictly preferable to encode biological structures using explicit Lagrangian coordinates, such as continuous curves, surface meshes, or point clouds.

    However, performing statistical or dynamic analysis on such explicit geometric data presents a profound mathematical challenge: shapes do not belong to a standard Euclidean vector space. One cannot simply compute the linear addition or scalar multiplication of two distinct complex organs (a problem often summarized as ``statistics without a `+'" \cite{feydy2020geometric}). Overcoming this fundamental limitation requires endowing sets of anatomical shapes with rigorous non-linear metric structures. 
    
     This necessity motivates the application of the theoretical framework of shape spaces, which offers a powerful paradigm by abstracting entire geometric configurations as individual points on a high-dimensional manifold. Well-known finite-dimensional labeled examples include Kendall’s shape manifolds \cite{Kendall1984SHAPEMP} with a large variety of applications in shape data sets \cite{dryden2016statistical}, while contemporary shape spaces of continuous curves and surfaces form infinite-dimensional Riemannian manifolds \cite{younes2019shapes,dupuis1998variational}. Traditionally, the construction and metrics of these spaces have been rooted in pure geometry \cite{gu2008computational,srivastava2016functional,michor2005vanishinggeodesicdistancespaces,Kendall1984SHAPEMP}. However, accurately modeling how physical bodies evolve over time—whether through natural development, tumor expansion, or surgical intervention—necessitates a rigorous mathematical framework that integrates fundamental physical principles with the complex dynamics of biology. The biological growth model provides a compelling basis for bridging continuum mechanics and shape spaces as suggested in \cite{charon2023shape}. The fundamental physical principles of this approach are deeply rooted in the extensive literature on the mathematics and mechanics of biological growth through the framework of morphoelasticity \cite{goriely2017mathematics}. The work presents several models and a variety of examples, establishing a comprehensive three-dimensional theory of volumetric growth.

    Building upon this rigorous biomechanical foundation, Charon and Younes \cite{charon2023shape} theoretically translated the morphoelastic growth process into a formal Riemannian shape space framework. Crucially, they demonstrated that the continuous structural remodeling and elastic energy dissipation characterizing biological growth can be formulated as an action functional over time. Because the underlying elastic strain energy density inherently forms a symmetric and strictly positive-definite quadratic form, this physical energy functional naturally acts as a valid Riemannian metric tensor on the infinite-dimensional shape manifold. Consequently, the minimum total mechanical work required to drive an initial shape to a target configuration mathematically defines a rigorous geodesic distance. This profound connection allows us to reframe complex biological development as an optimal, energy-efficient continuous trajectory within a shape space.

     Our work follows the biological view introduced above, explores a numerical approach to solve an evolution equation derived from a simplified model, and establishes the foundational mathematical and computational framework for the forward problem, which is necessary to rigorously address the inverse problem, grounded in the theory of morphoelasticity.

  The manuscript is organized as follows: Section 1 is the introduction. Sections 2 and 3 review the standard governing equations of linear elasticity and perform the corresponding numerical simulations. Sections 4 and 5 introduce the morphoelastic model accounting for volumetric growth and do the simulations. Sections 6 through 8 do the preliminary analysis of the inverse problem, detail the transition to the dynamic shape space framework, and propose the solution as an optimization problem. Sections 9 to 10 detail the mathematical derivation of the regularized variational forms and the Mixed Finite Element implementation in FEniCSx, then present simulation results across varying geometries. Finally, Section 11 discusses the current computation challenges and future perspectives.

  The primary contribution of this work lies in the development and computational implementation of a unified framework that couples the mechanics of biological growth with the geometry of shape spaces. While the physical principles of morphoelasticity \cite{goriely2017mathematics}, the regularization operators \cite{beg2005computing}, and the mathematical formulation of growth patterns as an optimal control problem or geodesics in shape space \cite{charon2023shape} are well-established, this study introduces a novel numerical strategy.

	\section{Governing Equations of Linear Elasticity}
	\quad The behavior of a deformable solid body under load is described by a set of governing partial differential equations. These equations represent fundamental physical principles of elasticity. For the linear elasticity model, we have the following five assumptions.
	
	\begin{itemize}
		\item \textbf{Continuum Mechanics:} The displacement field is continuously differentiable.
		\item \textbf{Homogeneous Material:} The material properties are uniform and do not change with position.
		\item \textbf{Isotropic Material:} The material properties are the same in all directions.
		\item \textbf{Linear Elasticity:} The material follows the Generalized Hooke's Law, where stress is linearly proportional to strain.
		\item \textbf{Small Deformations:} The strains and rotations are small enough to neglect geometric nonlinearities.
	\end{itemize}
	
	\subsection{Model Form}
	The strong form of the linear elasticity problem consists of three main relations defined within a domain $\Omega$, under the assumptions above:
	
	\begin{enumerate}
		\item \textbf{Equilibrium Conditions (Force Balance):} This relates the internal stresses to the external body forces.
		\begin{equation}
			-\nabla \cdot \tens{\sigma} = \vec{f} \quad \text{in } \Omega
			\label{eq:equilibrium}
		\end{equation}
		In index notation, this is expressed as:
		\begin{equation}
			\sigma_{ij,i} + f_j = 0
		\end{equation}
		where the comma notation $\sigma_{ij,i}$ represents the partial derivative $\dfrac{\partial \sigma_{ij}}{\partial x_i}$ (summation over the repeated index $i$ is implied) \cite{gould2018introduction}.
		Here $\tens{\sigma}$ is the stress tensor and $\vec{f}$ denotes the body force acting on the object.
		\item \textbf{Kinematic Relations:} This defines the small strain tensor $\tens{\varepsilon}$ in terms of the displacement field $\vec{u}$.
		\begin{equation}
			\tens{\varepsilon}(\vec{u}) = \frac{1}{2}\left(\nabla \vec{u} + (\nabla \vec{u})^T\right) \quad (\text{or in index notation: } \varepsilon_{ij} = \frac{1}{2}(u_{i,j} + u_{j,i}))
		\end{equation}
		
		\item \textbf{Constitutive Law:} This relates stress to strain for an isotropic, linear elastic material.
		\begin{equation}
			\tens{\sigma}(\tens{\varepsilon}(\vec{u})) = \lambda \text{tr}(\tens{\varepsilon})I + 2\mu \tens{\varepsilon} \quad (\text{or in index notation: } \sigma_{ij} = \lambda \delta_{ij} \varepsilon_{kk} + 2\mu \varepsilon_{ij})
		\end{equation}
	\end{enumerate}
	Here, $\lambda$ and $\mu$ are the Lamé constants \cite{gould2018introduction}, which are related to the more common Young's Modulus ($E$) and Poisson's ratio ($\nu$) \cite{gould2018introduction} by:
	\begin{equation}
		E = \frac{\mu(3\lambda + 2\mu)}{\mu + \lambda}, \qquad \nu = \frac{\lambda}{2(\mu+\lambda)}
	\end{equation}
	Where the third equation can also be written as:
	\begin{equation}
		\varepsilon_{ij} = \frac{1}{E} \left[ (1 + \nu)\sigma_{ij} - \nu \delta_{ij} \sigma_{kk} \right]
	\end{equation}
		A fourth condition, the compatibility constraint, known as Saint-Venant compatibility equations \cite{gould2018introduction} or more generally the Bianchi formulas \cite{lee2019introduction}, ensures that the strain field corresponds to a single-valued continuous displacement field. Mathematically, the requirement for a single-valued displacement field addresses the path-independence of the integration process. Physically, this corresponds to the constrain of no overlaps of the material. However, for numerical purposes of FeniCSx where the displacement field is the primary unknown, this condition is implicitly satisfied and is not considered directly.
	\subsection{Navier-Cauchy (Lamé–Navier) Equation and Additional Settings}
	
	The final governing equation for linear elasticity, expressed in terms of displacements, is known as the Navier-Cauchy or Lamé–Navier \cite{sadd2020elasticity} equation, which can be derived by combining the equilibrium, constitutive, and kinematic equations (\cref{appendix:Lamé–Navier}).
    
    \begin{equation}
			(\lambda + \mu) \nabla(\nabla \cdot \mathbf{u}) + \mu \nabla^2 \mathbf{u} + \mathbf{f} = \mathbf{0}
	\end{equation}  

   Within our setting (\cref{appendix:Lamé–Navier}), the governing system departs from a pure Neumann formulation, thereby securing the uniqueness of the solution by constraining the potential kernels associated with rigid motions.

More specifically, to complete the boundary value problem and ensure a unique solution, appropriate boundary conditions must be specified on the whole domain boundary $\partial\Omega$. We partition the boundary into two disjoint sets, $\partial\Omega = \Gamma_D \cup \Gamma_N$ (with $\Gamma_D \cap \Gamma_N = \emptyset$), avoiding a pure Neumann formulation means $\Gamma_D \neq \emptyset$.

\begin{itemize}
	\item \textbf{Dirichlet BC, $\Gamma_D$:} Prescribes that the displacement field in the region is zero ($\vec{u} = \vec{0}$).
	
	\item \textbf{Neumann BC, $\Gamma_N$:} Prescribes the natural boundary condition, where the surface traction is specified on \(\Gamma_N\). In the traction-free case, \(T=0\).

\end{itemize}

	\subsection{Derivation of the Variational Formulation}
    To solve these equations numerically using the Finite Element Method and to implement simulations using FEniCSx \cite{baratta_2025_18101307}, we first derive the variational (or weak) formulation. While this can be achieved through direct integration by parts and Stokes' theorem (see \cref{appendix: Variational_Linear} for this derivation), this thesis adopts the \textbf{Principle of Minimum Potential Energy} as the unifying framework.
    
    We prefer the energy method for two distinct reasons. First, it reveals that the variational problem is equivalent to minimizing a convex energy functional, which explicitly clarifies that the associated bilinear form $a(\cdot, \cdot)$ in the next section is \textbf{symmetric and positive semidefinite}. Establishing these properties is crucial, as they not only ensure the well-posedness of the static elastic problem but also provide the necessary mathematical foundation for defining the Riemannian metric in the shape space analysis presented later. Second, this framework allows us to treat linear elasticity consistently as a special case of the morphoelastic model (where the growth tensor $\tens{G}$ vanishes), which is introduced later \cref{sec:growth_model}.
    
	% \subsection{Summary of the Variational Problem}
	% \quad Denote $[H^1(\Omega)]^d$ is the vector-valued Sobolev space (where $d=3$ is the spatial dimension). The problem derived above can be stated concisely: find the displacement field $\vec{u} \in \mathcal{V}$ such that for all test functions $\vec{v} \in \hat{\mathcal{V}}$: 
    
	% \begin{equation}
	% 	a(\vec{u}, \vec{v}) = \mathcal{L}(\vec{v})
	% \end{equation}

 %    $\mathcal{V}$ and $\hat{\mathcal{V}}$ are defined as:
	% \begin{align}
	% 	\mathcal{V} &= \left\{ \vec{u} \in [H^1(\Omega)]^d \mid \vec{u} = \vec{u}_D \text{ on } \partial\Omega \right\} \\
	% 	\hat{\mathcal{V}} &= \left\{ \vec{v} \in [H^1(\Omega)]^d \mid \vec{v} = \vec{0} \text{ on } \partial\Omega \right\}
	% \end{align}
    
	% This formulation ensures a unique solution as introduced earlier. The bilinear form $a(\vec{u}, \vec{v})$ and the linear form $\mathcal{L}(\vec{v})$ are defined as:
	% \begin{align}
	% 	a(\vec{u}, \vec{v}) &= \int_{\Omega} \tens{\sigma(\varepsilon(\vec{u})}) : \nabla \vec{v} \, dV \\
	% 	\mathcal{L}(\vec{v}) &= \int_{\Omega} \vec{f} \cdot \vec{v} \, dV + \int_{\Gamma_N} \vec{T} \cdot \vec{v} \, dS
	% \end{align}
	% where the stress tensor $\tens{\sigma(\varepsilon(\vec{u})})$ is given by the constitutive law:
	% \begin{equation}
	% 	\tens{\sigma(\varepsilon(\vec{u})}) = \lambda(\nabla \cdot \vec{u})\mathbf{I} + \mu(\nabla \vec{u} + (\nabla \vec{u})^T).
	% \end{equation}

 %    Here, the notation $A:B = \mathrm{tr}(A^TB)$ denotes the usual inner product between matrices of same size.

    \subsection{Summary of the Variational Problem}

We now summarize the weak formulation of the linear elasticity problem. 
Let $\Omega \subset \mathbb{R}^d$ be a bounded physical domain, where in our numerical examples $d=3$. 
The displacement field is a vector-valued function
\[
    \vec{u}:\Omega \to \mathbb{R}^d,
\]
where $\vec{u}(x)$ describes the displacement of the material point originally located at $x\in\Omega$. 
The test function is also vector-valued:
\[
    \vec{v}:\Omega \to \mathbb{R}^d.
\]
The body force is given by
\[
    \vec{f}:\Omega \to \mathbb{R}^d,
\]
and the boundary traction on the Neumann part of the boundary is
\[
    \vec{T}:\Gamma_N \to \mathbb{R}^d.
\]

The strain tensor associated with a displacement field is
\[
    \tens{\varepsilon}(\vec{u})
    = \frac{1}{2}\left(\nabla \vec{u}+(\nabla \vec{u})^T\right),
\]
so that
\[
    \tens{\varepsilon}(\vec{u})(x)\in \mathbb{R}^{d\times d}_{\mathrm{sym}}
    \qquad \text{for each } x\in\Omega.
\]
The stress tensor is defined by the constitutive law
\[
    \tens{\sigma}(\tens{\varepsilon}(\vec{u}))
    = \lambda \operatorname{tr}(\tens{\varepsilon}(\vec{u}))\mathbf{I}
    +2\mu \tens{\varepsilon}(\vec{u}),
\]
and therefore
\[
    \tens{\sigma}(\tens{\varepsilon}(\vec{u})):\Omega
    \to \mathbb{R}^{d\times d}_{\mathrm{sym}}.
\]

The natural function space for the displacement field is the vector-valued Sobolev space
\[
    [H^1(\Omega)]^d
    =
    \left\{
    \vec{u}:\Omega\to\mathbb{R}^d
    \mid
    u_i\in H^1(\Omega),\ i=1,\dots,d
    \right\}.
\]
The admissible space and the test space are defined as
\begin{align}
    \mathcal{V}
    &=
    \left\{
    \vec{u}\in [H^1(\Omega)]^d
    \mid
    \vec{u}=\vec{u}_D \text{ on } \Gamma_D
    \right\},\\
    \hat{\mathcal{V}}
    &=
    \left\{
    \vec{v}\in [H^1(\Omega)]^d
    \mid
    \vec{v}=\vec{0} \text{ on } \Gamma_D
    \right\}.
\end{align}
Here $\vec{u}_D$ denotes the prescribed displacement on the Dirichlet boundary.
In the numerical examples considered in this work, the fixed boundary condition
corresponds to the homogeneous case $\vec{u}_D=\vec{0}$ on $\Gamma_D$.
Thus, in these examples, the admissible displacement fields satisfy
$\vec{u}=\vec{0}$ on $\Gamma_D$.
The test functions vanish on $\Gamma_D$ because variations must preserve the
prescribed Dirichlet boundary condition. However, they are free on $\Gamma_N$,
where the Neumann boundary condition appears naturally in the weak formulation.

The weak formulation is: find $\vec{u}\in\mathcal{V}$ such that for all $\vec{v}\in\hat{\mathcal{V}}$,
\begin{equation}
    a(\vec{u},\vec{v})=\mathcal{L}(\vec{v}).
\end{equation}
The bilinear form $a(\vec{u},\vec{v})$ and the linear functional $\mathcal{L}(\vec{v})$ are defined by
\begin{align}
    a(\vec{u},\vec{v})
    &=
    \int_{\Omega}
    \tens{\sigma}(\tens{\varepsilon}(\vec{u})):\nabla \vec{v}\,dV,\\
    \mathcal{L}(\vec{v})
    &=
    \int_{\Omega}\vec{f}\cdot \vec{v}\,dV
    +
    \int_{\Gamma_N}\vec{T}\cdot \vec{v}\,dS.
\end{align}
Equivalently, since $\tens{\sigma}(\tens{\varepsilon}(\vec{u}))$ is symmetric, one may also write
\[
    \tens{\sigma}(\tens{\varepsilon}(\vec{u})):\nabla \vec{v}
    =
    \tens{\sigma}(\tens{\varepsilon}(\vec{u})):\tens{\varepsilon}(\vec{v}).
\]
Thus the more standard elasticity form is
\begin{equation}
    a(\vec{u},\vec{v})
    =
    \int_{\Omega}
    \tens{\sigma}(\tens{\varepsilon}(\vec{u})):
    \tens{\varepsilon}(\vec{v})\,dV.
\end{equation}

Here, for matrices $A,B\in\mathbb{R}^{d\times d}$, the notation
\[
    A:B=\operatorname{tr}(A^TB)=\sum_{i,j=1}^d A_{ij}B_{ij}
\]
denotes the Frobenius inner product.

\clearpage

	\section{Simulation Examples based on the Model of Linear Elasticity}
        
\subsection{Geometric Visualization}

In the context of shape analysis, our primary objective is to study the geometric evolution of the domain $\Omega$ rather than the internal mechanical stress distribution. Consequently, in the following simulation examples, we visualize the initial configuration and the deformed configuration, where the coordinate transformation is defined by the mapping $\phi(\mathbf{x}) = \mathbf{x} + \mathbf{u}(\mathbf{x})$, with $\mathbf{u}$ the displacement field.

To ensure a more readable presentation, any vectors with coordinates denoted as $(x, y, z)$ in the text are treated as column vectors.

	\subsection{Example 1: Cantilever Beam with Gravity and Simulation Results}
	This is the starting model: a cantilever beam fixed at one end and subject to gravity.

	\begin{figure}[H]
		\centering
		\includegraphics[width=0.6\textwidth]{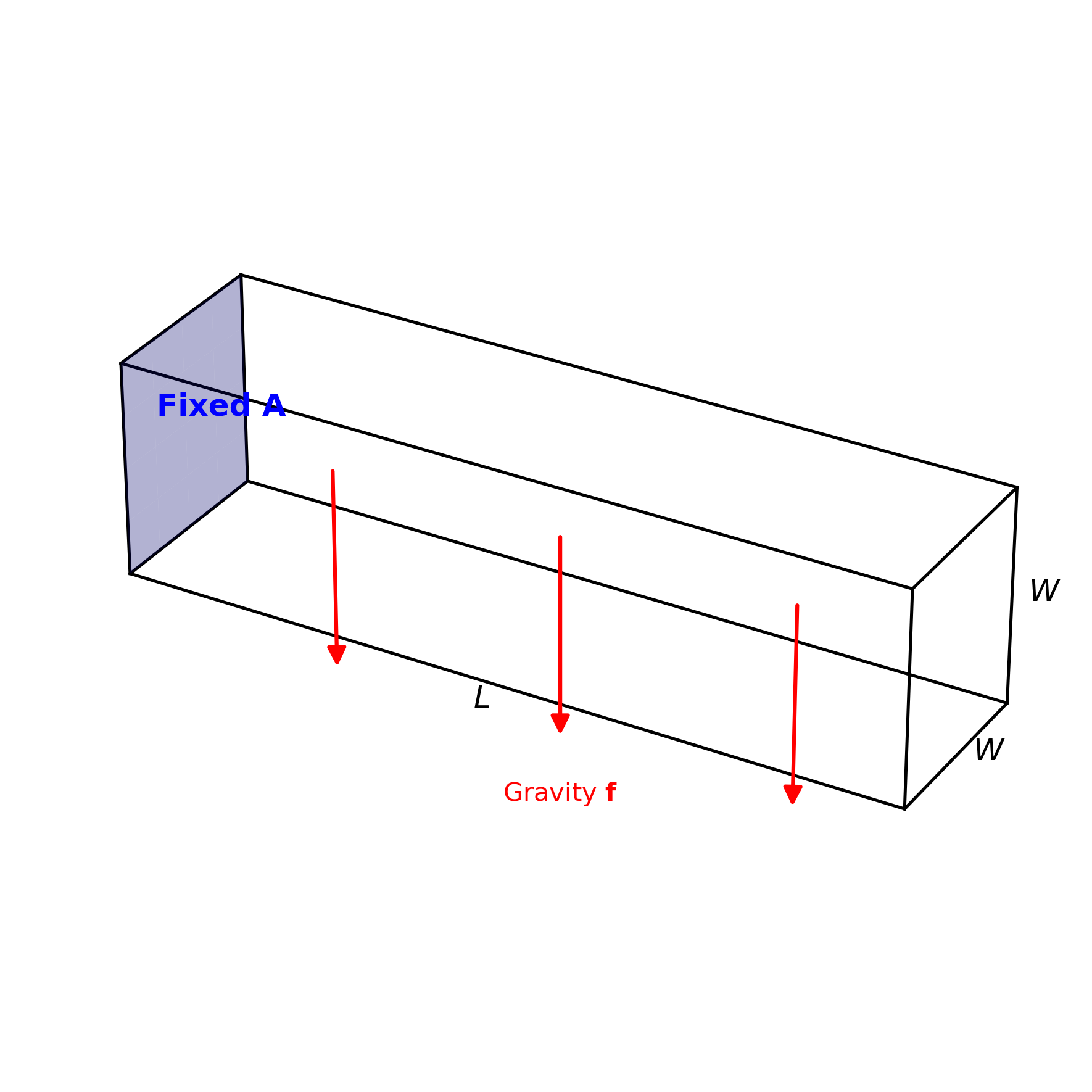}
		\caption{Schematic for Example 1.}
		\label{fig:sketch_EX1-0}
	\end{figure}
	
	\begin{itemize}
		\item \textbf{Geometry:} Rectangular beam with length $L=1$ and width/height $w=0.2$. 
		\item \textbf{Body Force:} Gravity acting downwards, $\vec{f} = (0, 0, -\rho g)$.
		\item \textbf{Dirichlet BC ($\Gamma_D$):} The left face, $A$, is fixed. $\vec{u} = \vec{0}$ on $A$.
		\item \textbf{Neumann BC ($\Gamma_N$):} The rest of the boundary, $\Gamma \setminus A$, is traction-free. $\vec{T} = \vec{0}$ on $\Gamma \setminus A$.
	\end{itemize}
\clearpage

\begin{figure}[H]
	\centering
	\includegraphics[width=0.8\textwidth]{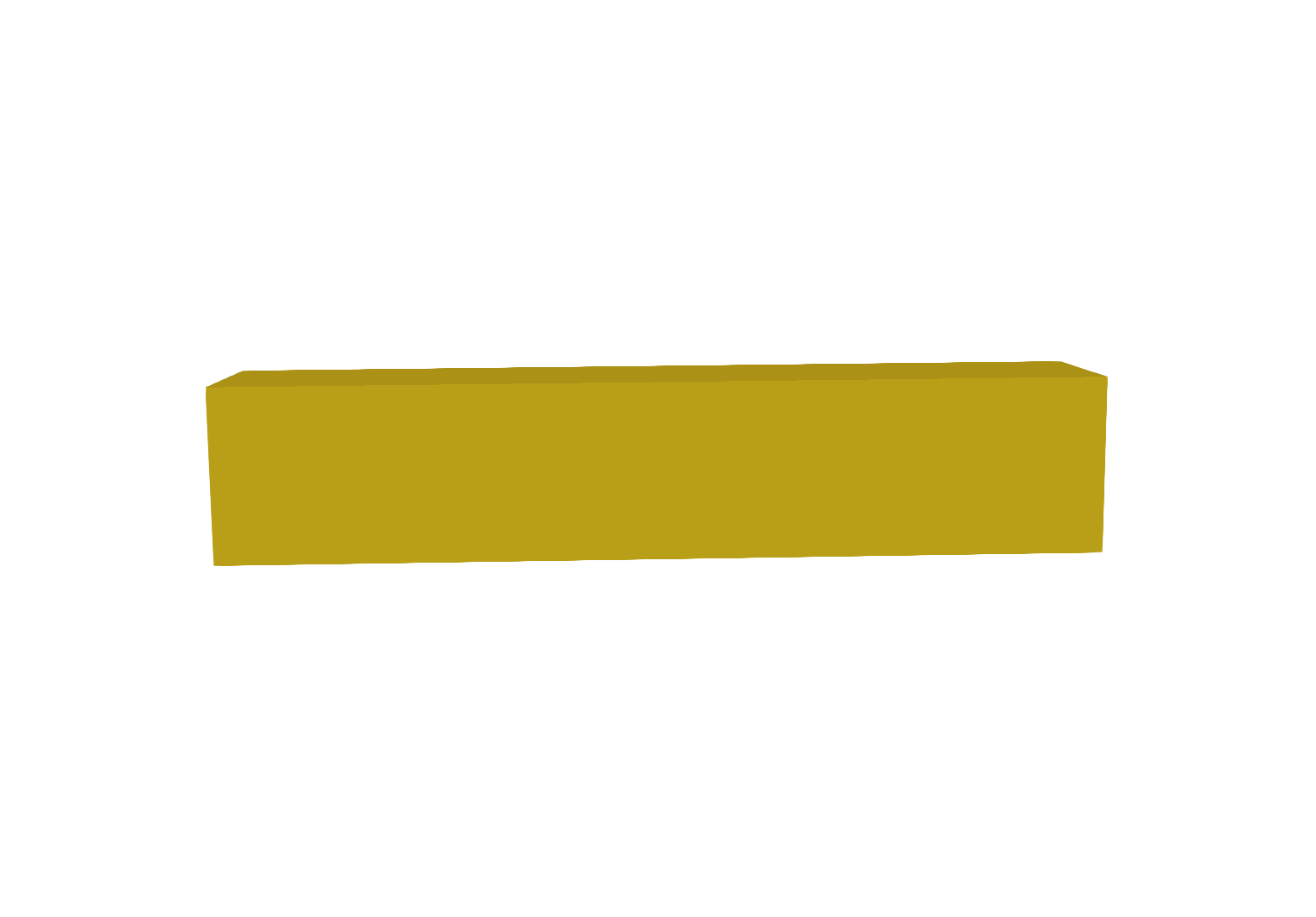}
	\caption[Reference configuration ($\Omega_0$) of the cantilever beam for Example 1]{Reference configuration ($\Omega_0$) of the cantilever beam used in Example 1. The geometry is a rectangular domain with length $L=1$ and cross-sectional width $w=0.2$. The mesh represents the discretized domain used for the Finite Element Method. The left face ($x=0$) is subject to a homogeneous Dirichlet boundary condition ($\mathbf{u}=\mathbf{0}$), serving as the initial state for the shape evolution process.}
	\label{fig:sketch_EX1-1}
\end{figure}

\begin{figure}[H]
	\centering
	\includegraphics[width=0.8\textwidth]{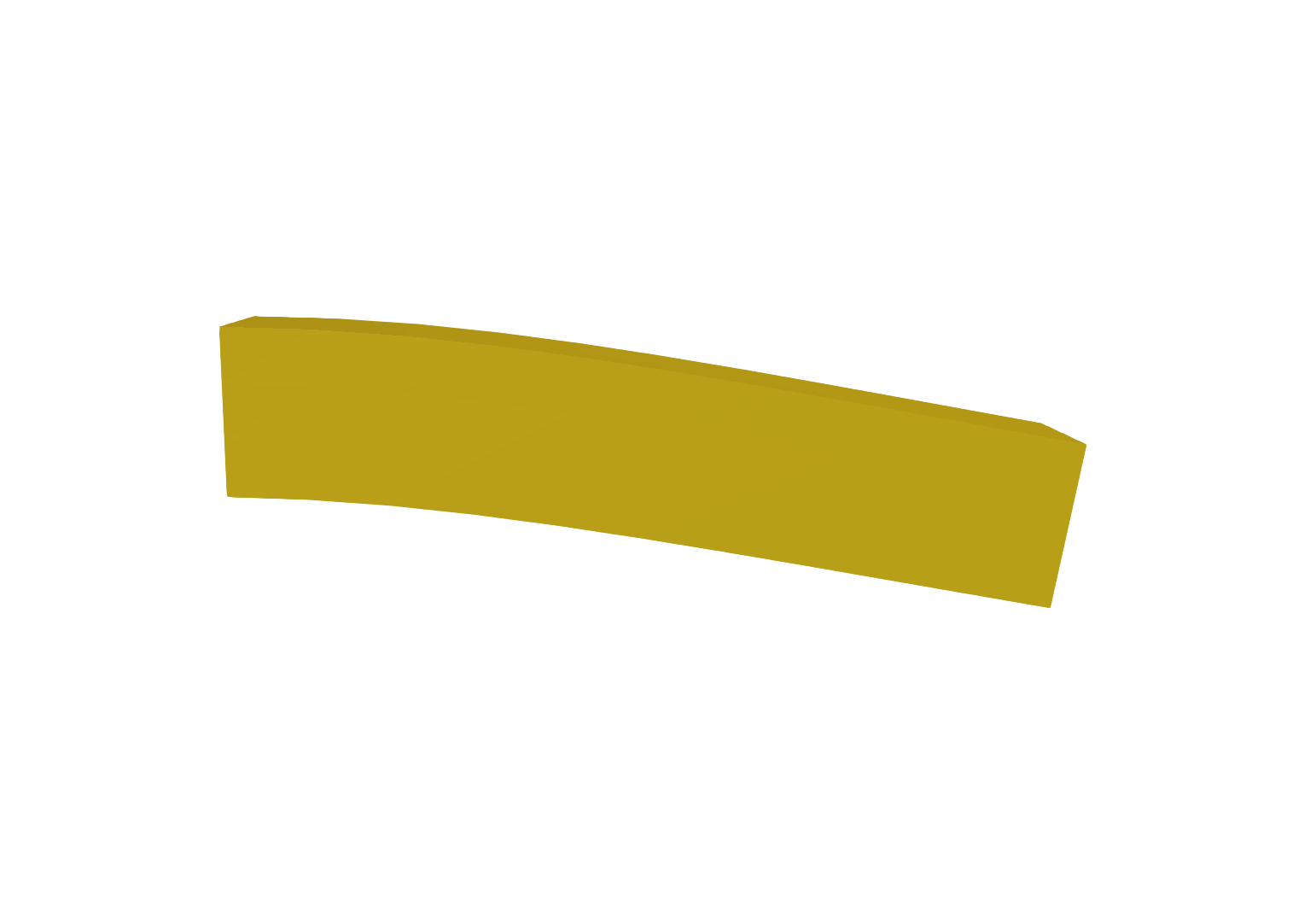}
	\caption[Deformed configuration ($\Omega_1$) of the cantilever beam under gravitational loading]{Deformed configuration ($\Omega_1$) of the cantilever beam under gravitational loading. This shape results from solving the standard linear elasticity equilibrium equations with a downward body force density $\mathbf{f}=(0, 0, -\rho g)$. The beam exhibits significant bending in the $-z$ direction.}
	\label{fig:sketch_model_1-2}
\end{figure}

	\clearpage

	\subsection{Example 2: Ball in a Hole}
	This model simulates a ball with the uniformly force on its upper hemisphere (denoted $A_2$) and a fixed lower hemisphere (denoted $A_1$).
	
	\begin{figure}[H]
		\centering
		\includegraphics[width=0.6\textwidth]{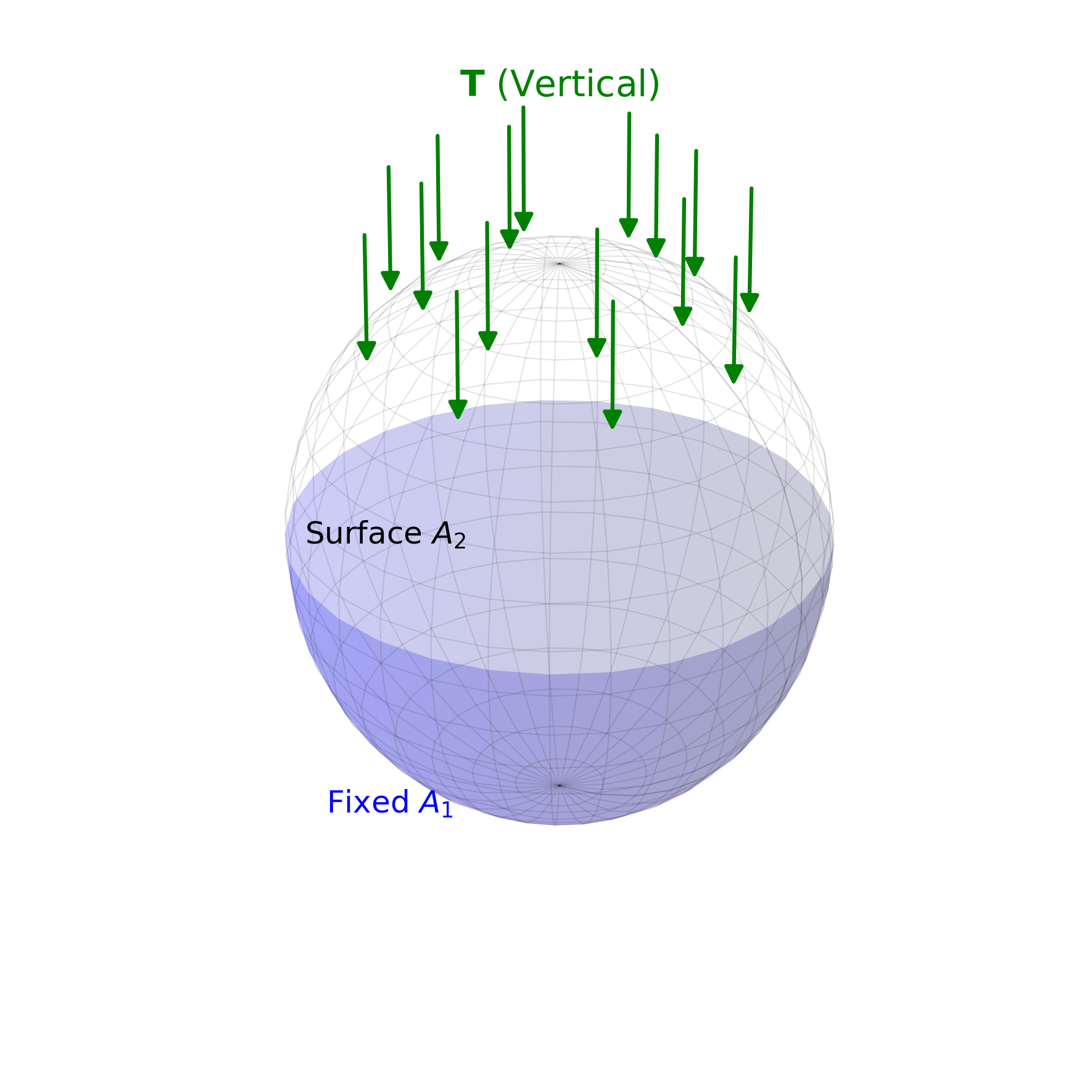}
		\caption{Schematic for Example 2}
		\label{fig:sketch_model_6}
	\end{figure}
	
	\begin{itemize}
		\item \textbf{Geometry:} A ball centered at $ (0,0,0) $ with radius 0.5.
		\item \textbf{Body Force:} Ignore gravity, so that $\vec{f} = (0, 0, 0)$.
		\item \textbf{Dirichlet BC ($\Gamma_D$):} The lower hemisphere is fixed:
		$\Gamma_D = A_1$, so that
		$\vec{u} = \vec{0}$ on $\Gamma_D$.
		
        \item \textbf{Neumann BC ($\Gamma_N$):} A uniform vertical downward traction $\vec{T}$ is applied to the upper hemisphere $A_2 = \Gamma \setminus \Gamma_D$. The magnitude of the traction is derived from a total vertical force $F = 0.4\,\text{N}$ (where $\text{N}$ denotes Newtons), such that:
        \begin{equation}
        \vec{T} = \left( 0, 0, -\frac{F}{\text{Area}(A_2)} \right) \quad \text{on } A_2.
        \end{equation}
        This formulation approximates the downward load as a constant vertical traction field rather than a normal pressure.
	\end{itemize}

    \clearpage
	
    \begin{figure}[H]
    \centering
    \includegraphics[width=0.7\textwidth]{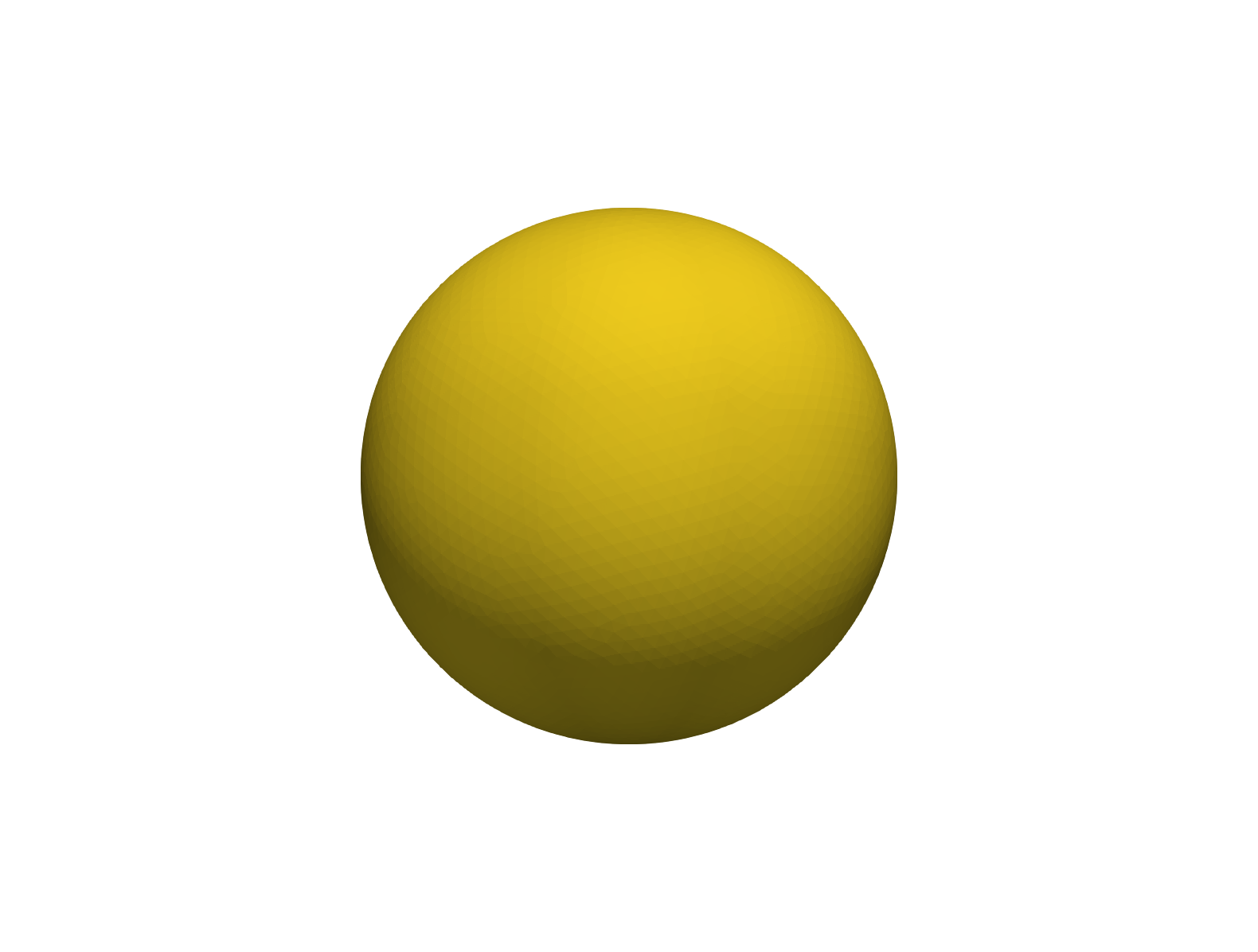}
    \caption[Reference configuration ($\Omega_0$) of the sphere for Example 2]{Reference configuration ($\Omega_0$) for Example 2. The geometry is a sphere of radius $R=0.5$ centered at the origin. The lower hemisphere ($z < 0$) is subject to a fixed Dirichlet boundary condition ($\mathbf{u}=\mathbf{0}$), acting as a rigid support for the structure.}
    \label{fig:sketch_EX2-0}
\end{figure}

\begin{figure}[H]
    \centering
    \includegraphics[width=0.7\textwidth]{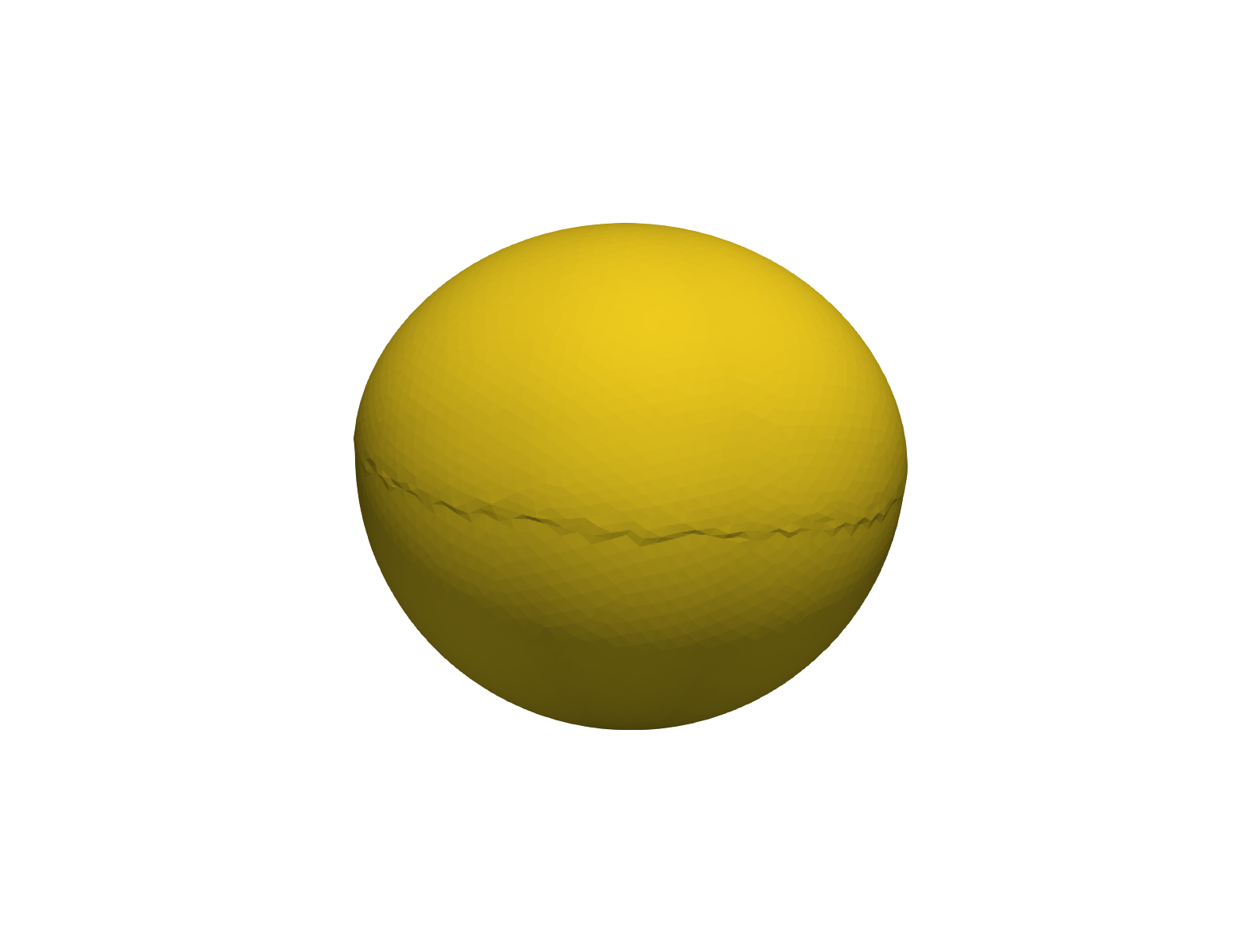}
    \caption[Deformed configuration ($\Omega_1$) of the sphere under surface traction]{Deformed configuration ($\Omega_1$) of the sphere under surface traction. A uniform downward pressure (corresponding to a total force $F=0.4\,\text{N}$) is applied to the upper hemisphere. In the absence of body forces, the deformation is driven entirely by the external load, resulting in a visible vertical compression against the fixed base.}
    \label{fig:sketch_EX2-1}
\end{figure}

\clearpage

\section{Variational Formulation with an Additive Growth Tensor}
\label{sec:growth_model}

We are now considering a modification to the linear elasticity model to account for volumetric growth, which refers to changes occurring within the interior of an object. There are numerous natural examples from biology, such as the addition of new material within biological tissues and the growth of plants or animals. This kind of behavior in the object can be modeled by $\textbf{morphoelasticity}$ as described in \cite{goriely2017mathematics}. In our construction of this model, we remain within the framework of linear elasticity, assuming that this growth does not alter the material's behavior.

\subsection{Morphoelastic Model Formulation} \label{Growthmodel}

% We introduce a \textbf{growth tensor} $\tens{G}(\vec{x})$, which is a symmetric second-order tensor field that describes the inherent, localized growth at each point $\vec{x} \in \Omega$. The key assumption is an additive decomposition of the total strain. And we assume $\tens{G}(\vec{x})$ is at least $C^1$.

We introduce a \textbf{growth tensor} $\tens{G}(\vec{x})$, which is a symmetric second-order tensor field that describes the inherent, localized growth at each point $\vec{x}\in\Omega$. More precisely,
\[
    \tens{G}:\Omega\to \mathbb{R}^{d\times d}_{\mathrm{sym}},
\]
so that for each $\vec{x}\in\Omega$, $\tens{G}(\vec{x})$ is a symmetric $d\times d$ matrix. This choice is consistent with the infinitesimal strain tensor
\[
    \tens{\varepsilon}(\vec{u}):\Omega\to \mathbb{R}^{d\times d}_{\mathrm{sym}},
\]
and therefore the difference
\[
    \tens{\varepsilon}(\vec{u})-\tens{G}
\]
is well-defined as a symmetric tensor field.

Physically, $\tens{G}$ represents the local growth strain, or the strain that would be induced by growth in the absence of elastic resistance. The elastic strain is then the part of the total strain that remains after subtracting this growth contribution. The key assumption is an additive decomposition of the total strain. For the strong-form derivation, we assume $\tens{G}$ is at least $C^1$, so that $\nabla\cdot\tens{\sigma}(\tens{G})$ is classically well-defined.

\begin{enumerate}
	\item \textbf{Strain Decomposition}: The total strain $\tens{\varepsilon}(\vec{u})$, derived from the displacement field, is the sum of the \textbf{elastic strain} $\tens{\varepsilon}_g$ and the growth strain $\tens{G}$. The material's stress is a response to the elastic strain only.
	\begin{equation}
		\tens{\varepsilon}_g = \tens{\varepsilon}(\vec{u}) - \tens{G}
	\end{equation}
	
	\item \textbf{Constitutive Law}: The stress tensor $\tens{\sigma}_g$ is related to the elastic strain $\tens{\varepsilon}_g$ by the standard linear isotropic constitutive law:
	\begin{equation}
		\tens{\sigma}_g = \tens{\sigma}(\tens{\varepsilon}_g) = \lambda \text{tr}(\tens{\varepsilon}_g)\mathbf{I} + 2\mu\tens{\varepsilon}_g
	\end{equation}
	Substituting the strain decomposition gives the stress in terms of displacement and growth:
	\begin{equation}
		\tens{\sigma}_g(\tens{\varepsilon}(\vec{u}), \tens{G}) = \lambda \text{tr}(\tens{\varepsilon}(\vec{u}) - \tens{G})\mathbf{I} + 2\mu(\tens{\varepsilon}(\vec{u}) - \tens{G})
	\end{equation}
	
	\item \textbf{Equilibrium}: The equilibrium equation states that the divergence of the internal stress must balance the external body forces $\vec{f}$:
	\begin{equation}
		-\nabla \cdot \tens{\sigma}_g = \vec{f}
	\end{equation}
	By substituting the constitutive law, we can see that the growth term acts as an effective body force:
	\begin{equation}
		-\nabla \cdot (\tens{\sigma}(\tens{\varepsilon}(\vec{u})) - \tens{\sigma}(\tens{G})) = \vec{f} \implies -\nabla \cdot \tens{\sigma}(\tens{\varepsilon}(\vec{u})) = \vec{f} - \nabla \cdot \tens{\sigma}(\tens{G})
	\end{equation}
	where we write $\tens{\sigma}(\tens{G}) = \lambda \text{tr}(\tens{G})\mathbf{I} + 2\mu\tens{G} $ for simplicity.

%     \item \textbf{A Specific Growth Tensor Definition}: To model localized volumetric growth, we define the growth tensor $\tens{G}$ using a Gaussian function. The spatial distribution is governed by:
% \begin{equation}
%     \tens{G}(\vec{x}) = g(\vec{x})\mathbf{I} = A_g \exp\left(-\frac{|\vec{x} - \vec{x}_0|^2}{2\tau^2}\right) \mathbf{I}
% \end{equation}
% Here, $\vec{x}_0$ represents the \textbf{center of the growth region}, determining the specific spatial location within the domain $\Omega$ where the growth intensity is maximal. The scalar $A_g$ denotes the peak growth amplitude, while $\tau$ controls the spatial spread of the growth effect.We take $\tau$ small enough so that $g\simeq 0$ on $\partial\Omega$, equivalently $\tens{G}\simeq \mathbf{0}$ on $\partial\Omega$.

\item \textbf{A Specific Growth Tensor Definition}: To model localized volumetric growth, we define the growth tensor $\tens{G}$ using a Gaussian function.

In this isotropic example, the growth tensor is taken to be a scalar multiple of the identity:
\[
    \tens{G}(\vec{x})=g(\vec{x})\mathbf{I}.
\]
Thus the growth is locally the same in every spatial direction. The scalar function
\[
    g:\Omega\to\mathbb{R}
\]
controls the magnitude of the growth, while the identity matrix $\mathbf{I}$ determines its isotropic tensor structure.

More specifically, the scalar spatial distribution is governed by the Gaussian function
\begin{equation}
    g(\vec{x})
    =
    A_g \exp\left(-\frac{|\vec{x} - \vec{x}_0|^2}{2\tau^2}\right),
\end{equation}
and hence
\begin{equation}
    \tens{G}(\vec{x})
    =
    g(\vec{x})\mathbf{I}
    =
    A_g \exp\left(-\frac{|\vec{x} - \vec{x}_0|^2}{2\tau^2}\right)\mathbf{I}.
\end{equation}

Here, $\vec{x}_0$ represents the \textbf{center of the growth region}, determining the specific spatial location within the domain $\Omega$ where the growth intensity is maximal. The scalar $A_g$ denotes the peak growth amplitude, while $\tau$ controls the spatial spread of the growth effect. We take $\tau$ small enough so that $g\simeq 0$ on $\partial\Omega$, equivalently $\tens{G}\simeq \mathbf{0}$ on $\partial\Omega$.

\end{enumerate}

From a biological perspective, it's noticeable that points closer to a growth center typically have a larger influence on growth, while points further away have less. This principle is straightforward to model for the simple, convex objects we are currently simulating. However, designing a biologically plausible model for more complex, non-convex geometries becomes significantly more challenging.

\subsection{Derivation of the Variational Formulation}
\label{Variational_with_Growth}

We derive the weak form using the \textbf{Principle of Minimum Potential Energy}, which states that a conservative system in stable equilibrium occupies a configuration that minimizes its total potential energy \cite{eslami2014finite}.

The total potential energy $E(\vec{u})$ is the sum of the internal strain energy stored in the elastic deformation and the potential energy of the external loads.

\begin{enumerate}
	\item \textbf{Energy Density}: The strain energy density $\Psi$ depends only on the elastic strain $\tens{\varepsilon}_g$. It is given by \cite{sadd2020elasticity}:
	\begin{equation}
		\Psi(\tens{\varepsilon}_g) = \frac{\lambda}{2}(\text{tr}(\tens{\varepsilon}_g))^2 + \mu\text{tr}(\tens{\varepsilon}_g^2)
	\end{equation}
	
	\item \textbf{Total Potential Energy}: The total energy functional $E(\vec{u})$ is:
	\begin{equation}
		E(\vec{u}) = \int_{\Omega} \Psi(\tens{\varepsilon}(\vec{u}) - \tens{G}) \,dV - \int_{\Omega} \vec{f} \cdot \vec{u} \,dV - \int_{\Gamma_N} \vec{T} \cdot \vec{u} \,dS
	\end{equation}
	The negative signs on the force terms arise because as forces do positive work, the potential energy of the system decreases ($W = -\Delta E_p$).
\end{enumerate}

To find the minimum, we take the first variation of the energy, $\delta E$, with respect to a small, arbitrary perturbation in displacement $\vec{u} \to \vec{u} + \delta \vec{v}$, and set it to zero.
\begin{equation}
	\delta E = \left[ \frac{d}{d\delta} E(\vec{u} + \delta \vec{v}) \right]\bigg|_{\delta=0} = 0
\end{equation}

Then, we can derive the new variational formulation using variational calculus. Details are in \cref{appendix:Variational_Morpho}

\subsection{Summary of the New Variational Problem}
\label{Summary_Var_morpho}

The modified problem can be stated concisely: find the displacement field $\vec{u} \in \mathcal{V}$ such that for all test functions $\vec{v} \in \hat{\mathcal{V}}$:
\begin{equation}
	a(\vec{u}, \vec{v}) = \mathcal{L}(\vec{v})
\end{equation}
The bilinear form $a(\vec{u},\vec{v})$ is identical to the standard linear elasticity problem:
\begin{equation}
	a(\vec{u}, \vec{v}) = \int_{\Omega} \tens{\sigma(\varepsilon(\vec{u})}) : \nabla \vec{v} \,dV
\end{equation}
The linear form $\mathcal{L}(\vec{v})$, however, now contains an additional term representing the effects of growth:
\begin{equation}
	\mathcal{L}(\vec{v}) = \int_{\Omega} \vec{f} \cdot \vec{v} \,dV + \int_{\Gamma_N} \vec{T} \cdot \vec{v} \,dS + \int_{\Omega} \tens{\sigma}(\tens{G}) : \nabla \vec{v} \,dV
	\label{eq:variation_Growth}
\end{equation}
Here, $\tens{\sigma}(\tens{G}) = \lambda \text{tr}(\tens{G})\mathbf{I} + 2\mu\tens{G} $.
Since $\tens{G}$ is known, by the conclusion in \cref{solution}, the equation has a unique solution.
\clearpage

\section{Simulation Examples based on the Model of Morphoelasticity}

\subsection{Model 3: Cantilever Beam with Gravity and Growth}

This model simulates a cantilever beam fixed at one end, subjected to both downward gravity and a localized internal volumetric growth.

\begin{figure}[H]
	\centering
	\includegraphics[width=0.6\textwidth]{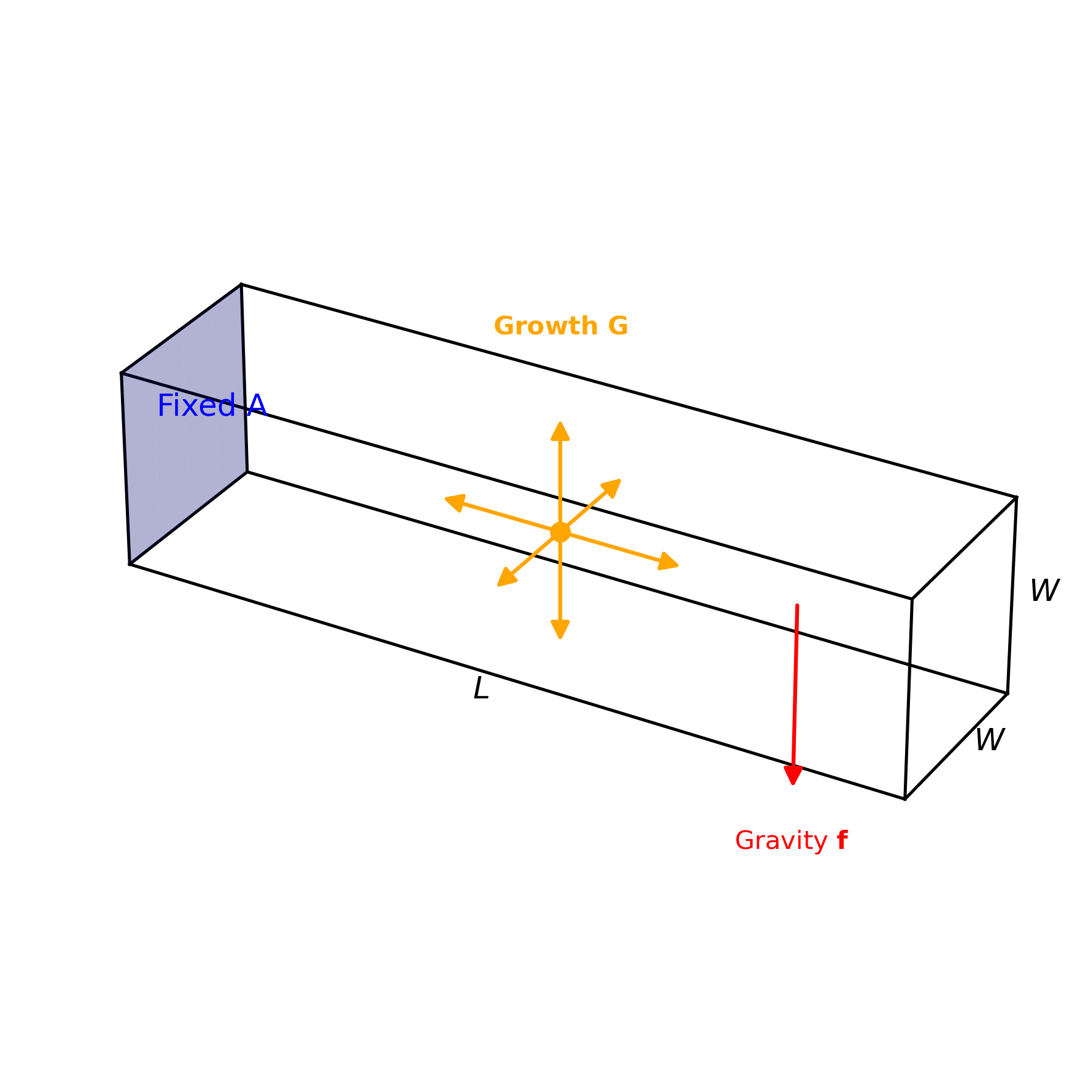}
	\caption[Schematic for Model 3]{Schematic for Model 3. The red arrows indicate the direction of gravity, the orange arrows represent the localized internal growth, and the blue region denotes the fixed boundary.}
	\label{fig:sketch_model_9}
\end{figure}

\begin{itemize}
	\item \textbf{Geometry:} A rectangular beam with length $L=1$ and width/height $w=0.2$.
    \item \textbf{Growth:} A localized growth region is centered at $(L/2, w/2, w/2)$, with a growth amplitude coefficient $A_g = 0.2, \tau = 0.1$. 
	\item \textbf{Body Force:} Gravity acts downwards, $\vec{f} = (0, 0, -\rho g)$.
	\item \textbf{Dirichlet BC ($\Gamma_D$):} The left face, denoted as $A$, is fixed. Thus, $\vec{u} = \vec{0}$ on $A$.
	\item \textbf{Neumann BC ($\Gamma_N$):} The remainder of the boundary, $\Gamma \setminus A$, is traction-free. $\vec{T} = \vec{0}$ on $\Gamma \setminus A$.
\end{itemize}

\begin{figure}[H]
	\centering
	\includegraphics[width=0.8\textwidth]{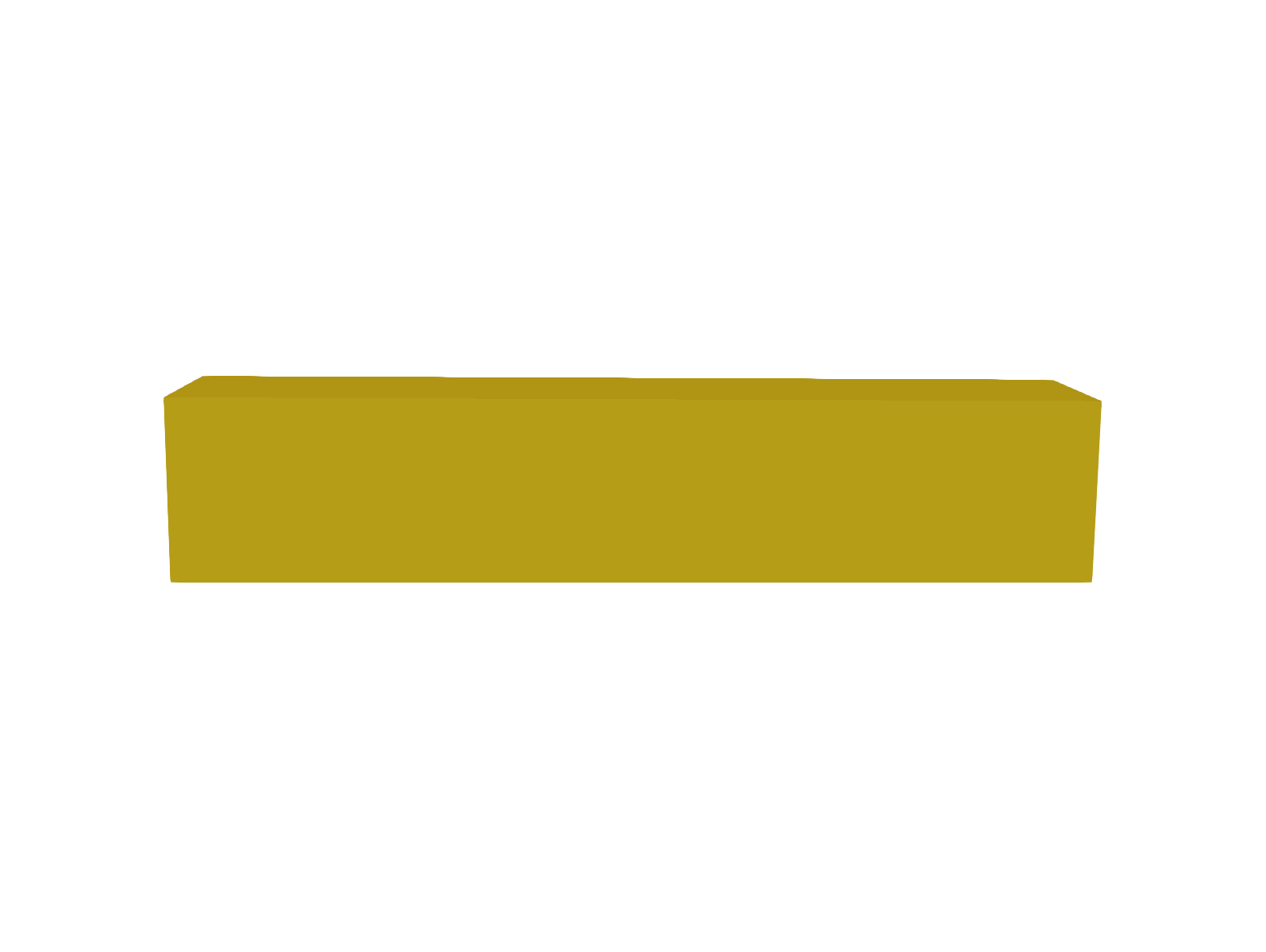}
	\caption[Reference configuration ($\Omega_0$) for Model 3]{Reference configuration ($\Omega_0$) for Model 3. A straight rectangular cantilever beam supported rigidly at the left boundary.}
	\label{fig:sketch_EX3-0}
\end{figure}
	
\begin{figure}[H]
	\centering
	\includegraphics[width=0.8\textwidth]{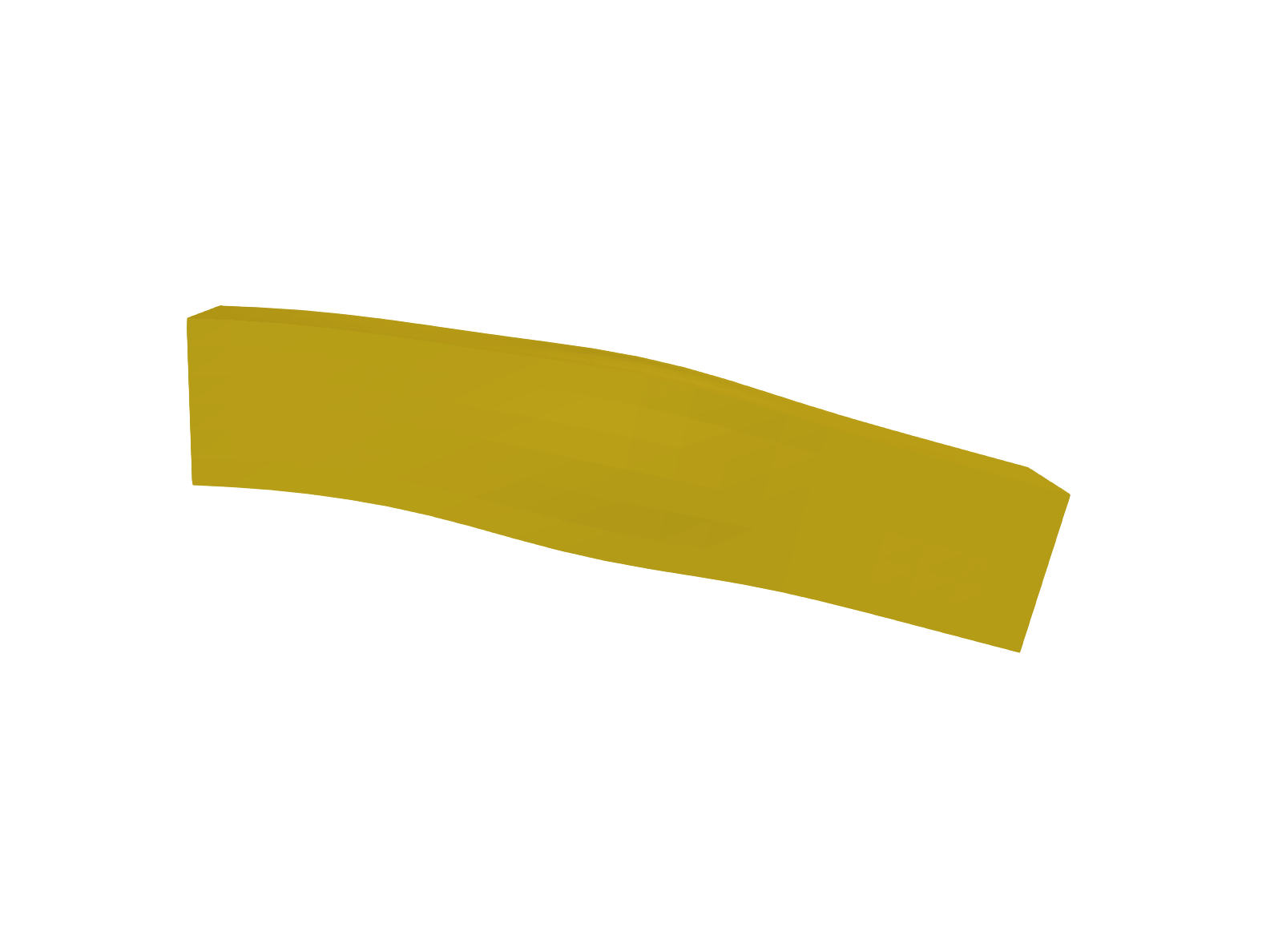}
	\caption[Deformed configuration ($\Omega_1$) of the beam with localized growth]{Deformed configuration ($\Omega_1$) of the beam. The final shape demonstrates the combined mechanical response: a downward bending induced by the external body force and a localized volumetric expansion driven by the internal growth tensor $\tens{G}$ at the center.}
	\label{fig:sketch_EX3-1}
\end{figure}

\clearpage

\subsection{Example 4: Ball with Pressure and Growth}

This model extends Example 2 by introducing an internal growth mechanism. It simulates a sphere subjected to a uniform downward surface traction on its upper hemisphere, a fixed lower hemisphere, and a localized growth region near its center.

\begin{figure}[H]
	\centering
	\includegraphics[width=0.6\textwidth]{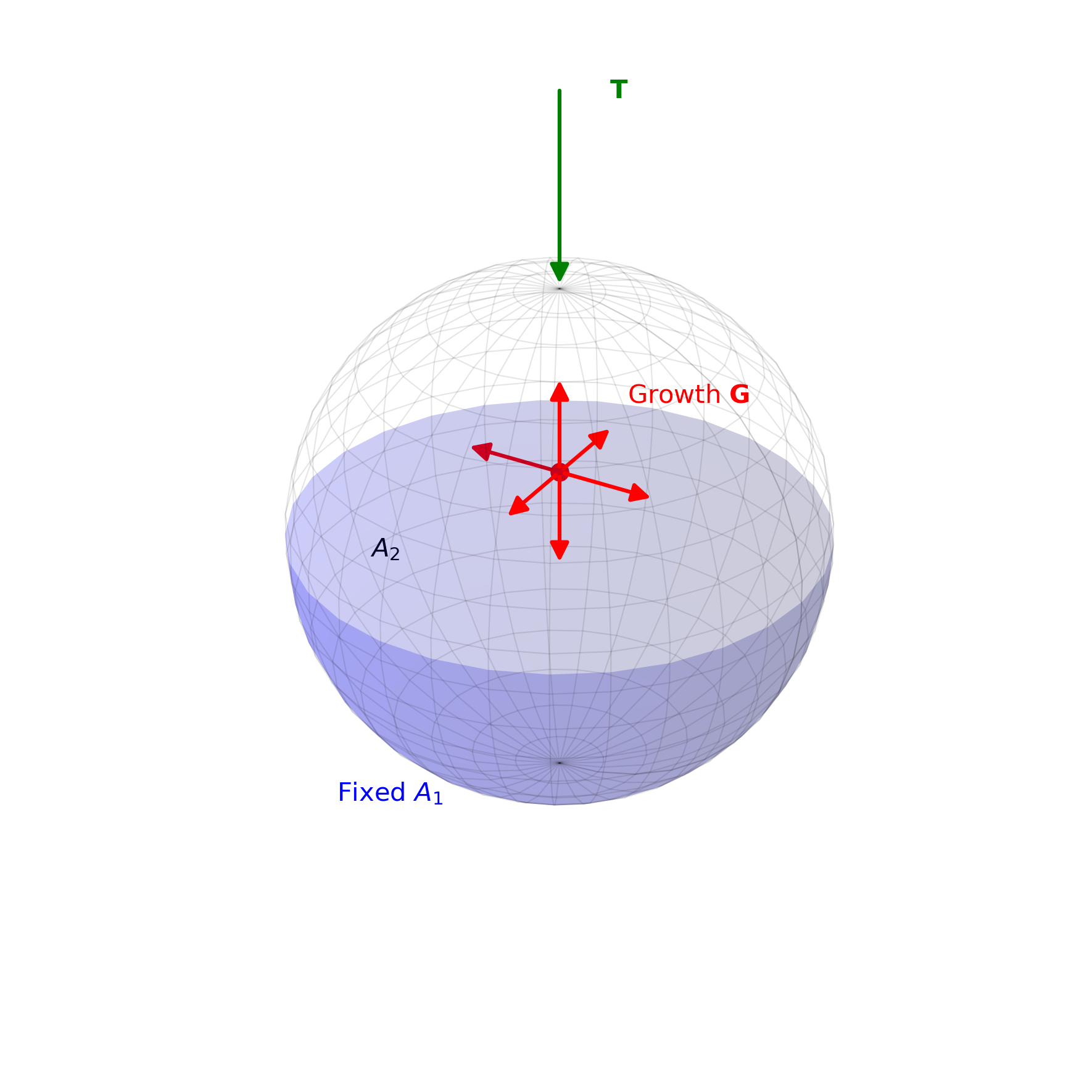}
	\caption[Schematic for Example 4]{Schematic for Example 4. Green arrows represent the external downward pressure, while red arrows indicate the internal growth.}
	\label{fig:sketch_model_7}
\end{figure}

\begin{itemize}
	\item \textbf{Geometry:} A sphere centered at $(0,0,0)$ with radius $R=0.5$.
    \item \textbf{Growth:} The localized growth seed is positioned at $(0, 0, 0.25)$, governed by the growth coefficient $A_g$ and $\tau = 0.1$.
	\item \textbf{Body Force:} Gravity is neglected, yielding $\vec{f} = (0, 0, 0)$.
	\item \textbf{Dirichlet BC ($\Gamma_D$):} The lower hemisphere is fixed, defined as $\Gamma_D = A_1$. Thus, $\vec{u} = \vec{0}$ on $\Gamma_D$.
	\item \textbf{Neumann BC ($\Gamma_N$):} A uniform downward traction $\vec{T} = \frac{\vec{F}}{\text{Area}(A_2)}$ is applied on the upper hemisphere $A_2 = \Gamma \setminus \Gamma_D$. Assuming a total downward force magnitude of $F=0.4\,\text{N}$, we approximate the pressure as a vertical downward traction.
\end{itemize}

\clearpage

\subsection{The Interplay of Pressure and Growth: A Critical Threshold}

Numerical simulations of Example 4 reveal a competitive interplay between the external pressure $\vec{T}$ and the internal growth tensor $\tens{G}(\vec{x})$. Specifically, the growth coefficient $A_g$ dictates whether the sphere exhibits a net upward expansion or downward compression. For values of $A_g$ below a certain critical threshold, the $z$-component of the displacement field $\vec{u}$ becomes predominantly negative.

\begin{figure}[H]
	\centering
	\includegraphics[width=0.7\textwidth]{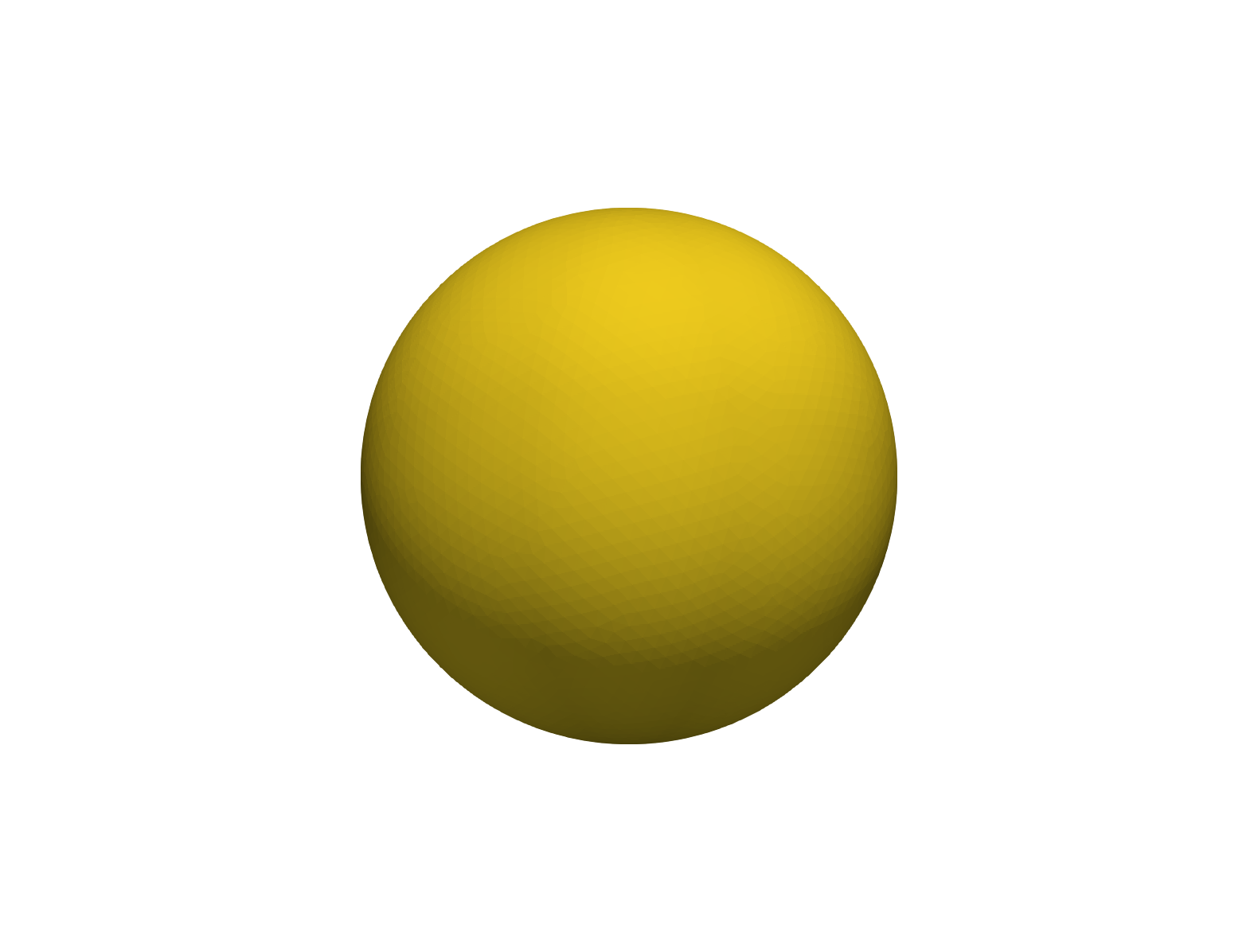}
	\caption{Reference configuration ($\Omega_0$) for Example 4.}
	\label{fig:sketch_EX4-0}
\end{figure}

\begin{figure}[H]
	\centering
	\includegraphics[width=0.7\textwidth]{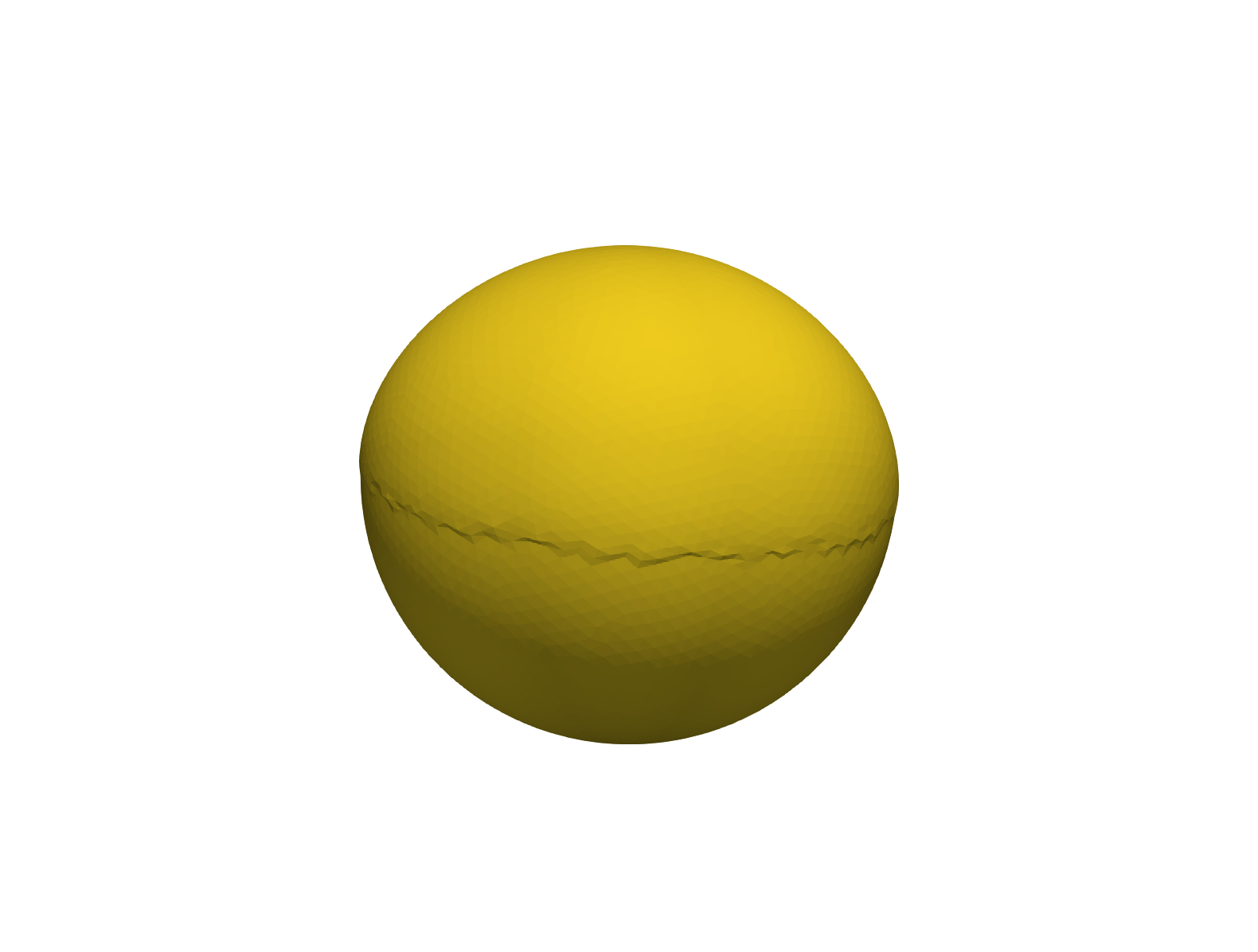}
	\caption[Deformed configuration with a low growth coefficient ($A_g = 0.1$)]{Deformed configuration with a low growth coefficient ($A_g = 0.1$). The external pressure dominates the internal growth, resulting in a net structural compression (negative vertical displacement).}
	\label{fig:sketch_EX4-1}
\end{figure}

\begin{figure}[H]
	\centering
	\includegraphics[width=0.7\textwidth]{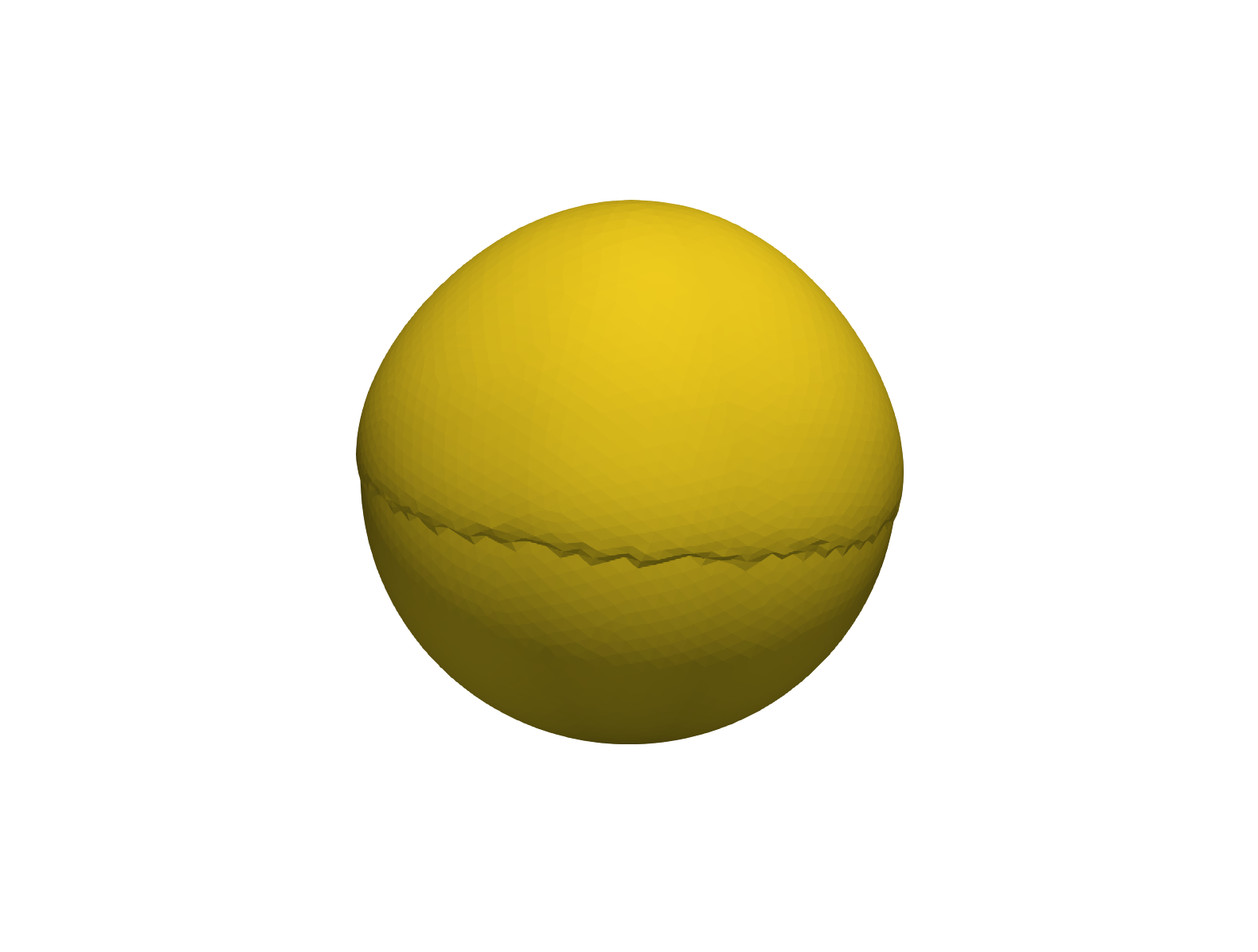}
	\caption[Deformed configuration with an intermediate growth coefficient ($A_g = 1$)]{Deformed configuration with an intermediate growth coefficient ($A_g = 1$). The internal growth begins to significantly counteract the external pressure.}
	\label{fig:sketch_EX4-2}
\end{figure}

\begin{figure}[H]
	\centering
	\includegraphics[width=0.8\textwidth]{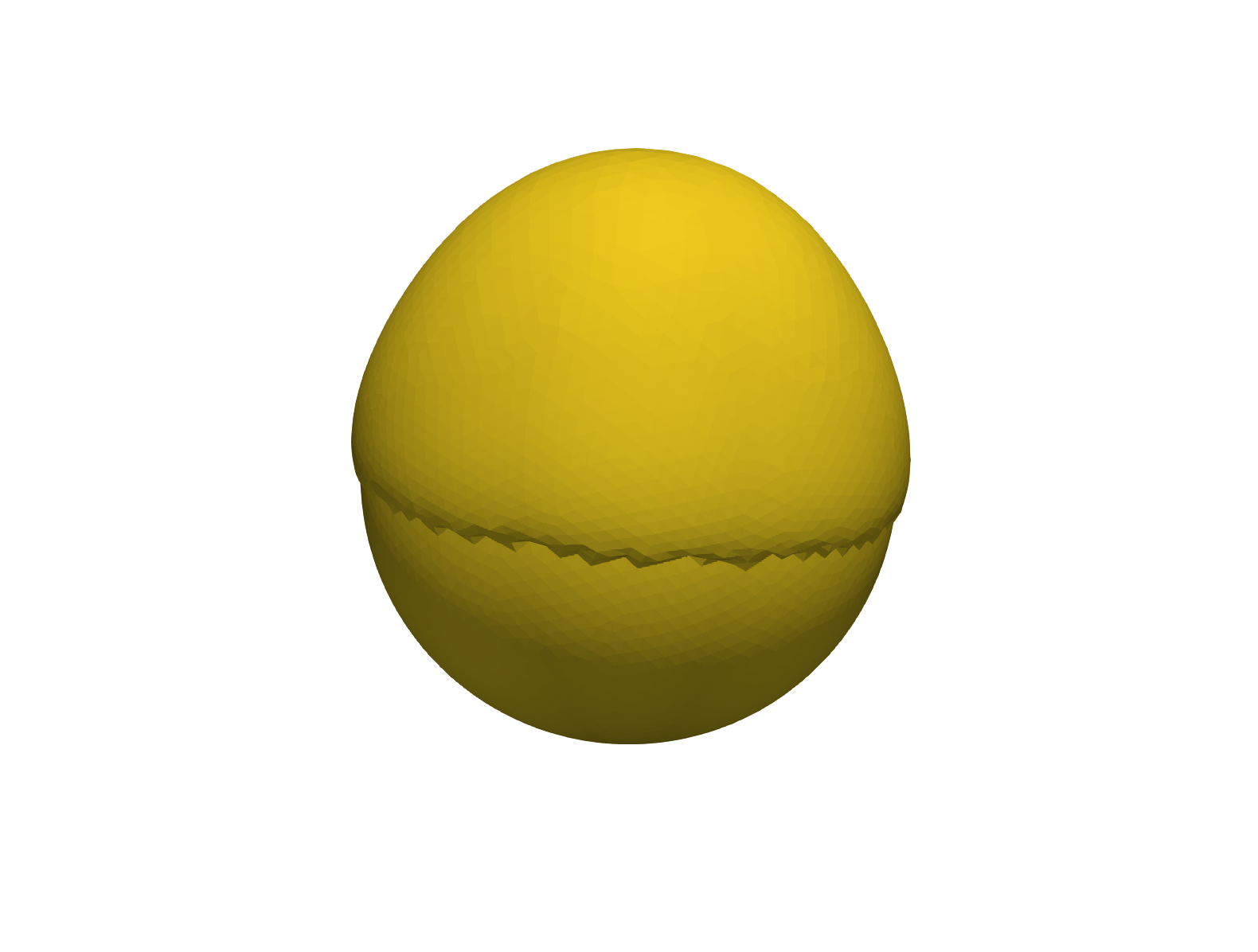}
	\caption[Deformed configuration with an intermediate growth coefficient ($A_g = 2$)]{Deformed configuration with a high growth coefficient ($A_g = 2$). The internal growth potential entirely overcomes the external compressive load, leading to an upward volumetric expansion.}
	\label{fig:sketch_EX4-3}
\end{figure}

While determining the exact analytical value of this critical threshold is mathematically challenging, these numerical results perfectly align with the underlying principles of morphoelasticity \cite{charon2023shape}. In this framework, growth acts as an internal driving force. When the growth potential ($A_g$) is insufficient to overcome the mechanical work exerted by the external traction, the structure is inevitably compressed rather than inflated.

\section{Preliminary Analysis of the Inverse Problem}

Thus far, much of the modeling and simulation effort has focused on the forward problem: solving for the displacement vector field $\bm{u}$ given a prescribed growth. This raises a natural question for an inverse problem: if the final displacement field $\bm{u}$ is known (e.g., from experimental or clinical measurements), can we determine the underlying growth tensor $\bm{G}$ that caused it? 

\subsection{Formulation and Ill-Posedness of the Inverse Problem}
By rearranging the equilibrium equation and utilizing the linearity of the stress operator $\bm{\sigma}(\cdot)$, we arrive at a linear partial differential equation for the unknown tensor field $\bm{G}$:
\begin{equation}
	\nabla \cdot \bm{\sigma}(\bm{G}) = \bm{h}(\bm{x}),
	\label{eq:inverse_pde}
\end{equation}
where
\begin{equation}
	\bm{h}(\bm{x}) = \bm{f}(\bm{x}) + \nabla \cdot \bm{\sigma}(\bm{\varepsilon}(\bm{u}(\bm{x})))
	\label{eq:h_definition}
\end{equation}
is a known vector field entirely determined by the observed data and external forces. 

The uniqueness of the solution depends on the properties of the linear operator $\mathcal{L}(\bm{G}) = \nabla \cdot \bm{\sigma}(\bm{G})$. A unique solution exists only if the null space (or kernel) of $\mathcal{L}$ is trivial, i.e., if $\mathcal{L}(\bm{G}) = \bm{0}$ implies $\bm{G} = \bm{0}$. 

However, if there are no other constraints on $\bm{G}$ apart from it being a $C^1$ symmetric tensor, non-trivial tensor fields $\bm{G}_\mathrm{h} \neq \bm{0}$ exist that satisfy the homogeneous equation:
\begin{equation}
	\nabla \cdot \bm{\sigma}(\bm{G}_\mathrm{h}) = \bm{0}.
	\label{eq:homogeneous}
\end{equation}
To illustrate this, consider a non-trivial example with Lamé parameters $\lambda = -3$ and $\mu = 6$ (which yields a positive bulk modulus $K > 0$). The following tensor field satisfies the homogeneous equation:
\begin{equation}
	\bm{G}_\mathrm{h}(\bm{x}) =
	\begin{pmatrix}
		x & 0 & 0 \\
		0 & -\left(\frac{\lambda + 2\mu}{2\lambda}\right)x & 0 \\
		0 & 0 & -\left(\frac{\lambda + 2\mu}{2\lambda}\right)x
	\end{pmatrix}
	\label{eq:Gh_example}
\end{equation}

Mathematically, we can apply the theory of Beltrami representation theory introduced in \cite{sadd2020elasticity} to represent any stress field with zero divergence as $\Sigma = \nabla \times (\nabla \times \boldsymbol{\Phi})^T$ for an arbitrary smooth symmetric tensor field $\boldsymbol{\Phi}$, and since the linear constitutive operator $\sigma(\mathbf{G})$ is a bijection (as soon as $\mu\neq 0$ and $3\lambda + 2\mu \neq 0$), every such $\boldsymbol{\Phi}$ maps to a non-trivial growth tensor $\mathbf{G}_h \neq \mathbf{0}$ in the null space, which rigorously demonstrates that the inverse problem is fundamentally ill-posed with infinitely many solutions.

\subsection{Identifiable Growth Models as a Partial Solution}
To resolve this ambiguity, it is necessary to provide \textbf{identifiable growth models} by restricting the allowable form of $\bm{G}$. For instance, consider an isotropic growth model defined by a scalar field $g(\bm{x})$, such that $\bm{G} = g\mathbf{I}$. 

Substituting this into the stress operator gives $\bm{\sigma}(g\mathbf{I}) = \lambda \text{tr}(g\mathbf{I})\mathbf{I} + 2\mu (g\mathbf{I}) = (3\lambda + 2\mu)g\mathbf{I}$. The condition for falling into the null space, $\nabla \cdot \bm{\sigma}(\bm{G}) = \bm{0}$, then requires:
\begin{equation}
    (3\lambda + 2\mu) \nabla g = \bm{0} \implies (3\lambda + 2\mu) \partial_i g = 0 \quad \text{for } i=1, 2, 3.
\end{equation}
As long as $3\lambda + 2\mu \neq 0$, which is exactly consistent with our setting where the bulk modulus is strictly positive (\cref{appendix:Lamé–Navier}), this equation requires $\nabla g = \bm{0}$, meaning $g$ must be a constant on the domain $\Omega$. Consequently, any spatially varying parametric model, such as one where $g(\bm{x})$ is defined by a \textbf{Gaussian function}, cannot belong to the null space. Thus, by adopting such a parametric representation, the growth model becomes strictly identifiable.

\subsection{Argument of Defining Distance via Minimum Potential Energy}

In the shape space framework, determining the ``true'' growth mechanism equates to finding an optimal trajectory between an initial configuration, $\Omega_A$, and a final deformed configuration, $\Omega_B$. This perspective fundamentally shifts our focus from asking ``what is the single static cause?'' to ``what is the most physically plausible path connecting these two shapes?''.

To quantify the notion of a ``plausible path,'' a natural starting point is to define a cost function or distance metric based on the principle of minimum potential energy. An initial, intuitive approach might consider the total stored elastic energy required for a single-step deformation mapping the reference configuration $\Omega_A$ to the final state $\Omega_B$ via a displacement field $\bm{u}^*$:
\begin{equation}
	E(\bm{u}^*) = \int_{\Omega_A} \Psi(\bm{\varepsilon}(\bm{u}^*)) \, dX
\end{equation}
where $\Psi(\bm{\varepsilon}) = \frac{\lambda}{2}(\text{tr}(\bm{\varepsilon}))^2 + \mu \text{tr}(\bm{\varepsilon}^2)$ is the strain energy density, and $\bm{u}^*$ is the specific mapping that minimizes this energy among all valid transformations between the two shapes.

While conceptually straightforward, this static definition exhibits fundamental limitations. Mathematically, it cannot serve as a valid Riemannian distance metric on the shape manifold because it inherently violates symmetry---the elastic energy required to deform $\Omega_A$ into $\Omega_B$ is generally not equal to the energy required to deform $\Omega_B$ back into $\Omega_A$ due to the change in the reference integration domain. Furthermore, it fails to satisfy the triangle inequality. 

More importantly, from a mechanical standpoint, this one-step approach is fundamentally incompatible with our underlying physical model. The strain energy density $\Psi(\bm{\varepsilon})$ relies on linear elasticity, which is strictly valid only under the assumption of infinitesimally small deformations. Applying this static model to evaluate large, finite shape differences inherently violates this assumption, leading to physically inaccurate stress-strain approximations. 

Coupled with the biological reality that shape evolution is an accumulative, path-dependent process rather than a sudden, one-step elastic distortion, the static energy functional merely evaluates the energetic cost of a final state. It completely ignores the continuous structural remodeling, localized growth, and intermediate energy dissipation that characterize natural biological development. Therefore, to construct a rigorous and physically meaningful geodesic distance, a transition to a time-dependent, continuous evolutionary framework is strictly necessary. This will be introduced in \cref{sec:continuous_framework}.

\subsection{Argument of Change of Growth Tensor in Continuous Evolution}
In the preceding simulations (Sections 3 and 4), the models were executed as one-step static deformations. To transition to a shape space trajectory, we must now consider growth as a continuous process occurring over multiple steps.

Previously, we defined the isotropic growth tensor as:
\begin{equation}
	\tens{G}(\vec{x}) = g(\vec{x})\mathbf{I} = A_g \exp\left(-\frac{|\vec{x} - \vec{x}_0|^2}{2\tau^2}\right) \mathbf{I}
\end{equation}
Although this formulation produces reasonable static results, it presents a fundamental issue for continuous evolution. The definition $\tens{G}(\vec{x})$ relies on spatial coordinates ($\vec{x}$), representing an \textbf{Eulerian description}. As highlighted in modern geometric data analysis \cite{feydy2020geometric}, an Eulerian perspective is often ill-suited for tracking large, free-form anatomical deformations because the material points themselves move as the body deforms. Therefore, to ensure physical plausibility and accurately track the growth of specific tissue regions, the growth tensor must be formulated in a \textbf{Lagrangian description}, defined on the reference material coordinates $\mathbf{X}$, which is introduced in the next section \cref{eq:lagrangian_ref}.

This ensures that the growth source moves correctly with the deforming material during continuous simulation.

	\section{Continuous Formulation of the Growth Model}

The iterative simulation, where growth occurs in discrete steps, is a numerical approximation of an underlying continuous time evolution. To establish the theoretical foundation, we now formulate this process as a continuous model, assuming a quasi-static evolution. The quasi-static assumption implies that the growth process is slow enough for the body to be in mechanical equilibrium at every instant, allowing us to neglect inertial effects (i.e., terms involving acceleration $\rho \ddot{\vec{u}}$) \cite{charon2023shape}. 

\subsection{From Discrete Steps to Continuous Time}

We associate the iteration number $n$ with a continuous time variable $t$. The incremental displacement $\vec{u}_{\text{sol}}$ computed in each step corresponds to the displacement occurring over a small time interval $\Delta t$. This naturally introduces the concept of a ``velocity field," $\vec{v}(t,\vec{x})$, which describes the instantaneous motion of the material. The relationship is given by:
\begin{equation}
	\vec{u}_{\text{sol}} \approx \vec{v} \cdot \Delta t \quad \text{or equivalently,} \quad \vec{v}(t,\vec{x}) = \frac{\partial \vec{u}(t,\vec{x})}{\partial t}.
\end{equation}
All fields, including displacement, stress, and the domain itself, are now considered functions of time.

\subsection{From Static to Continuous Evolution Framework}
\label{sec:continuous_framework}

In classical elasticity, a physical body possesses a permanent memory of its initial, stress-free reference configuration. Consequently, one might intuitively attempt to formulate a continuous growth model by directly taking the time derivative of the static equilibrium equations. However, biological tissues are living materials that continuously undergo \textit{remodeling}. In every infinitesimally small time step $dt$, the tissue accommodates internal growth, undergoes slight elastic deformation, and rapidly remodels to dissipate internal stresses. In this process, it essentially ``forgets'' its previous configuration and treats its current, updated shape as a new stress-free reference state \cite{goriely2017mathematics}. 

Therefore, a rigorous continuous formulation cannot rely on the time derivative of a global elastic state. Instead, it must be constructed by accumulating the incremental energy dissipated during each infinitesimal remodeling step. Assume a homogeneous hyperelastic model where the elastic strain energy of a deformation mapping $\phi$ defined on a domain $\Omega$ is given by:
\begin{equation}
	\boldsymbol{W}(\Omega, \phi) = \int_{\Omega} W\big(\nabla\phi(\vec{x})^T\nabla\phi(\vec{x})\big) \, d\vec{x}
\end{equation}
where $W$ is defined on the set of symmetric matrices and takes values in $[0, +\infty)$, with the stress-free state yielding $W(I)=0$.

Assuming that $W$ is at least $C^2$ and using the fact that the identity matrix $I$ is a global minimum, the Taylor expansion around a purely small perturbation $\delta D$ at $\delta=0$ yields:
\begin{equation}
	W(I + \delta D) = W(I) + \nabla W(I)\cdot \delta D + \frac{1}{2} \nabla^2 W(I) (\delta D, \delta D) + o(\delta^2) = \frac{1}{2} \nabla^2 W(I) (\delta D, \delta D) + o(\delta^2)
\end{equation}
Here, the 4th-order tensor $\nabla^2 W(I)$ represents the initial stiffness tensor of the material. Its double contraction with the perturbation forms a positive-definite quadratic form on the set of symmetric matrices, which we define for short as $w(D) \doteq \nabla^2 W(I)(D, D)$. 

Now, assume that the domain $\Omega$ and mapping $\phi$ depend on time, with the current domain $\Omega_t = \phi(t,\Omega_0)$ acting as the updated reference. The spatial velocity field is given by $\vec{v}(t,\vec{x}) = \partial_t \phi(t,\vec{X})$. Taking an infinitesimal time $dt$ and using the quasi-static assumption, the incremental deformation mapping from time $t_k$ to $t_k + dt$ is defined as:
\begin{equation}
	\psi(t_k) \doteq \phi(t_k+dt)\circ \phi(t_k)^{-1} \approx \mathrm{id} + dt \cdot \vec{v}(t_k)
\end{equation}
The spatial gradient of this incremental mapping is $\nabla\psi = I + dt \nabla \vec{v}$. The corresponding right Cauchy-Green deformation tensor, which measures the pure strain, is calculated as:
\begin{align}
	\nabla\psi^T\nabla\psi &= (I + dt \nabla \vec{v}^T)(I + dt \nabla \vec{v}) \nonumber \\
	&= I + dt(\nabla \vec{v} + \nabla \vec{v}^T) + dt^2 (\nabla \vec{v}^T \nabla \vec{v}) \nonumber \\
	&= I + 2dt \cdot \tens{\varepsilon}(\vec{v}) + o(dt)
\end{align}
where $\tens{\varepsilon}(\vec{v}) = \frac{1}{2}(\nabla \vec{v} + \nabla \vec{v}^T)$. 

To incorporate volumetric growth, we model the instantaneous growth tensor between time $t$ and $t+dt$ as $\tilde{G}(t) = I + dt \cdot \mathfrak{g}(t)$, where $\mathfrak{g}(t)$ is the instantaneous growth rate tensor. According to the multiplicative decomposition of morphoelasticity ($F = F_e F_g$) \cite{goriely2017mathematics}, the effective elastic deformation is obtained by pulling back the total deformation by the inverse of the growth tensor. Using the approximation $\tilde{G}(t)^{-1} \approx I - dt \cdot \mathfrak{g}(t)$, since $dt$ small enough, the effective elastic strain metric becomes:
\begin{align}
	C_e &= \tilde{G}(t)^{-T} \big(\nabla\psi^T \nabla\psi\big) \tilde{G}(t)^{-1} \nonumber \\
	&\approx (I - dt \cdot \mathfrak{g}^T) \big(I + 2dt \cdot \tens{\varepsilon}(\vec{v})\big) (I - dt \cdot \mathfrak{g}) \nonumber \\
	&= I + 2dt \cdot \tens{\varepsilon}(\vec{v}) - dt \cdot \mathfrak{g} - dt \cdot \mathfrak{g}^T + o(dt)
\end{align}
Assuming isotropic or symmetric growth ($\mathfrak{g} = \mathfrak{g}^T$), this simplifies to:
\begin{equation}
	C_e = I + 2dt \big(\tens{\varepsilon}(\vec{v}) - \mathfrak{g}(t)\big) + o(dt)
\end{equation}

A critical mathematical observation arises here regarding the accumulation of energy over time. By substituting the effective strain perturbation $\delta D = 2dt (\tens{\varepsilon}(\vec{v}) - \mathfrak{g})$ into the second-order Taylor expansion $W(I+\delta D) \approx \frac{1}{2}w(\delta D)$, and utilizing the property of quadratic forms where $w(cX) = c^2 w(X)$, the energy dissipated in a single infinitesimal step is:
\begin{equation}
    W_{\text{step}} = W(I + \delta D) \approx \frac{1}{2} w\big( 2dt (\tens{\varepsilon}(\vec{v}) - \mathfrak{g}) \big) = 2dt^2 w\big( \tens{\varepsilon}(\vec{v}) - \mathfrak{g} \big) + o(dt^2).
\end{equation}
If we were to simply sum the total energy over the entire evolution $T$ (partitioned into $N = T/dt$ steps), the total energy $\sum W_{\text{step}} \approx dt \sum (2dt \cdot w)$ would scale with $O(dt)$ and approach zero as $dt \to 0$. To avoid this, we evaluate the limit of the energy divided by the time step $dt$. Thus, the total continuous action over time $[0, T]$ is defined as:
\begin{equation}
	\mathcal{J} = \lim_{dt \to 0} \frac{1}{dt} \sum_{k=0}^{T/dt} \boldsymbol{W} (\Omega_{t_k}, \psi(t_k)) \simeq 2\int_0^T\int_{\Omega_t} w\big(\tens{\varepsilon}(\vec{v}(t,\vec{x})) - \mathfrak{g}(t,\vec{x})\big) \, d\vec{x} \, dt
\end{equation}

For a linear isotropic elastic material, taking the standard strain energy density:
\begin{equation}
	W(D)= \frac{\lambda}{2} \mathrm{tr}(D)^2 + \mu \mathrm{tr}(DD^T) 
\end{equation}
$W$ is inherently a quadratic form, meaning $w = 2W$. Disregarding multiplicative constants, this exactly recovers our regularized variational model. Crucially, the functional relies on a positive-definite quadratic form, providing a rigorous Riemannian distance metric for the shape space.

\subsection{Connection to the Numerical Implementation}
\label{Numerical_imp}

The continuous action functional derived above  justifies our iterative numerical scheme in which the evolution is discretized into a sequence of instantaneous energy-minimization problems. 

\begin{enumerate}
	\item Time is discretized into small steps $\Delta t$, corresponding to the incremental intervals in the continuous theory.
	\item At each time step $t_n \to t_{n+1}$, the instantaneous growth rate $\mathfrak{g}(t_n)$ acts as the active source. We do not use finite differences of a global tensor, but instead evaluate the local instantaneous growth generation.
	\item The variational problem solved in each iteration:
	\begin{equation}
		a(\vec{u}_{\text{sol}}, \vec{v}) = L(\vec{v})
	\end{equation}
	is mathematically equivalent to minimizing the spatial integral of the instantaneous dissipation rate $\int_{\Omega_t} w\big(\tens{\varepsilon}(\vec{v}) - \mathfrak{g}\big) d\vec{x}$. We solve for the optimal velocity field $\vec{v}(t_n)$ (represented by the incremental displacement $\vec{u}_{\text{sol}} \approx \vec{v}(t_n) \Delta t$) that restores mechanical equilibrium for that specific increment.
	\item The final step updates the domain $\Omega_{t_n} \to \Omega_{t_{n+1}}$, naturally resetting the reference configuration and fulfilling the memory-less remodeling assumption of the morphoelastic framework.
\end{enumerate}

As formalized above, the core of the numerical simulation is to evaluate the instantaneous growth rate $\mathfrak{g}(t,\vec{x})$ and compute the optimal velocity field $\vec{v}$ to update the mesh. 

To ensure physical plausibility, the instantaneous growth rate must be tied to the material points, even as they move. The simplest choice is to define the growth generation in the initial, undeformed reference space $\Omega_0$, and push it forward to the current Eulerian space $\Omega_t$. 

Let the reference (initial) position of the growth center be denoted as $\vec{X}_0$. We prescribe an initial scalar growth rate distribution $\mathfrak{g}_0(\vec{X})$ based on the material distance to $\vec{X}_0$:
\begin{equation}
	\mathfrak{g}_0(\vec{X}) = A_g \exp\left(-\frac{\|\vec{X} - \vec{X}_0\|^2}{2\tau^2}\right)
	\label{eq:lagrangian_ref}
\end{equation}

To implement this in our simulation at any current time $t$, we utilize the space-time diffeomorphism $\phi(t,\vec{X}) = \vec{x}$. The material point $\vec{X}$ currently located at spatial position $\vec{x}$ is retrieved via the inverse mapping $\vec{X} = \phi(t)^{-1}(\vec{x})$. Thus, the instantaneous growth rate tensor evaluated in the current configuration becomes:
\begin{equation}
	\mathfrak{g}(\vec{x}, t) = \Big(\mathfrak{g}_0 \circ \phi(t)^{-1}(\vec{x})\Big) \tens{I} = A_g \exp\left(-\frac{\|\phi(t)^{-1}(\vec{x}) - \vec{X}_0\|^2}{2\tau^2}\right) \tens{I}
	\label{eq:lagrangian_current}
\end{equation}

In the discrete computational setting, the inverse mapping $\phi(t)^{-1}(\vec{x})$ is continuously tracked by accumulating the inverted displacement field $\vec{u}$ over time. 

	\section{From Inverse Problem to Shape Space and Geodesic Distance}
	
	As discussed, the static models developed previously allow us to solve the forward problem. However, a more challenging and often more relevant question in biology is the \textit{inverse problem}: given an observed shape change, what is the underlying growth pattern that caused it? Since the inverse problem is ill-posed, to address this ambiguity, we must move beyond a simple static equilibrium view and consider the entire deformation as a continuous path within the shape space. This perspective shifts the focus from finding a single cause for a final state to identifying the most plausible or "energy-efficient" trajectory between two shapes. This leads us to the concept of the length of shortest path between two shapes, measured by a physically meaningful metric.
	
	This section outlines the transition from our static framework to a dynamic one and connected it to the diffeomorphic mapping, culminating in an optimal control problem that defines such a geodesic path, inspired by the principles of diffeomorphic shape analysis \cite{charon2023shape}.

\subsection{Recap: Transition to a Dynamic Framework}

To model a continuous growth process, we introduce time, $t \in[0, 1]$. The shape of the object is no longer static but evolves continuously.
\begin{itemize}
	\item \textbf{Deformation Map:} The position of a material point, initially at $\vec{X} \in \Omega_0$, is given by a time-parameterized spatial diffeomorphism $\phi(t, \cdot)$ such that $\vec{x} = \phi(t, \vec{X})$ at time $t$.
	\item \textbf{Velocity Field:} The evolution is driven by a Eulerian velocity field $\vec{v}(t, \vec{x})$, representing the instantaneous speed of the material at spatial position $\vec{x}$. It governs the deformation map via the flow equation:
	\begin{equation}
		\partial_t \phi(t, \vec{X}) = \vec{v}\big(t, \phi(t, \vec{X})\big) \quad \text{or concisely} \quad \partial_t \phi(t) = \vec{v}(t) \circ \phi(t).
	\end{equation}
	This aligns with our discrete numerical steps where the velocity is evaluated as $\vec{v} \approx \frac{\vec{u}(\vec{X}, t + \Delta t ) - \vec{u}(\vec{X}, t)}{\Delta t}$.
	\item \textbf{Evolving Domain:} The current domain occupied by the object is simply $\Omega_t = \phi(t, \Omega_0)$.
	\item \textbf{Symmetric Velocity Gradient (Strain Rate Equivalent):} The instantaneous rate of elastic deformation is described by the tensor $\tens{\varepsilon}(\vec{v})$, defined as the symmetric part of the spatial velocity gradient:
	\begin{equation}
		\tens{\varepsilon}(\vec{v}) = \frac{1}{2} \left( \nabla\vec{v} + (\nabla\vec{v})^T \right).
	\end{equation}
\end{itemize}
In this dynamic context, our goal is to find an entire evolution path $\{\phi(t, \cdot)\}_{t \in[0,1]}$ that is optimal in some sense.

    \subsection{Defining a Norm for Growth Rate Tensors}

Following the Riemannian viewpoint in \cite{charon2023shape}, we define the ``cost'' of a continuous growth pattern. 

\textit{Note on notation: In the preceding continuous derivation (\cref{sec:continuous_framework}), we utilized $\mathfrak{g}$ to denote the instantaneous growth rate to strictly distinguish it from the finite static growth step. However, as we now formulate the continuous optimal control problem, we will use $\tens{G}$ to denote this \textbf{instantaneous growth tensor field} (i.e., replacing $\mathfrak{g}$ with $\tens{G}$). This overloads the symbol slightly but maintains notational consistency with our earlier variational forms.}

The norm-squared of a growth rate tensor field $\tens{G}$ on a domain $\Omega$ is defined as the minimum energy required to accommodate it:
\begin{equation}
	\|\tens{G}\|^2_{[\Omega]} = \min_{\vec{v} \in V} \left( \kappa \|\vec{v}\|_V^2 + \int_{\Omega} \Psi(\tens{\varepsilon}(\vec{v}) - \tens{G}) \, dV \right)
	\label{eq:growth_norm}
\end{equation}
Here:
\begin{itemize}
	\item $\tens{G}$ is the instantaneous growth tensor (symmetric), corresponding to the control field $g$ in \cite{charon2023shape}.
	\item $\Psi(\cdot)$ is the quadratic strain energy density form. As derived in \cref{sec:continuous_framework}, it now measures the instantaneous elastic power dissipation rate.
	\item The term $\tens{\varepsilon}(\vec{v}) - \tens{G}$ represents the instantaneous \textit{elastic} deformation rate metric.
	\item $\kappa \|\vec{v}\|_V^2$ is a regularization term ensuring the well-posedness of the minimization problem, corresponding to the reproducing kernel Hilbert space (RKHS) norm in \cite{charon2023shape}, which will be introduced in \cref{sec:variational_regularized}.
\end{itemize}
For each given growth rate $\tens{G}$, a unique velocity field $\vec{v}$ exists that minimizes \cref{eq:growth_norm} under suitable boundary conditions \cite{charon2023shape}.

\subsection{The Optimal Control Problem for Geodesic Paths}

With the growth rate norm defined, we can now formulate the problem of finding the most ``efficient'' trajectory—the geodesic path—between an initial shape $\Omega_0$ and a final shape $\Omega_1$. This is framed as an optimal control problem, where we seek the time-dependent growth field $\tens{G}(t)$ that drives the deformation while minimizing the total integrated cost.

\textbf{Objective:} Find the growth evolution $\tens{G}(t)$ that minimizes the total cost functional (the squared geodesic distance):
\begin{equation}
	J(\tens{G}(\cdot)) = \int_0^1 \|\tens{G}(t)\|^2_{[\Omega(t)]} \, dt
	\label{eq:cost_functional}
\end{equation}	
\textbf{Subject to the following constraints:} 
\begin{equation} \label{eq:constraints}
\begin{aligned}
    \partial_t \phi(t) &= \vec{v}(t) \circ \phi(t), \\
    \vec{v}(t) &= \arg\min_{\vec{v} \in V} \left( \kappa \|\vec{v}\|_V^2 + \int_{\Omega(t)} \Psi(\tens{\varepsilon}(\vec{v}) - \tens{G}(t)) \, dV \right).
\end{aligned}
\end{equation}
     And the path must accurately connect the start and end shapes:
	\begin{align}
		\phi(0, \Omega_0) &= \Omega_0 \\
		\phi(1, \Omega_0) &= \Omega_1
	\end{align}

As shown in \cite{charon2023shape}, this problem is mathematically equivalent to simultaneously finding the optimal pair of controls $(\vec{v}(t), \tens{G}(t))$ and the optimal state trajectory (driven by $\vec{v}(t)$) that minimize the combined action functional:
\begin{equation}
	\min_{\vec{v}(\cdot), \tens{G}(\cdot)} \int_0^1 \left( \kappa \|\vec{v}(t)\|_V^2 + \int_{\Omega(t)} \Psi(\tens{\varepsilon}(\vec{v}(t)) - \tens{G}(t)) \, dV \right) dt
	\label{eq:action_functional}
\end{equation}
This formulation (corresponding to Eq. 16 in \cite{charon2023shape}) defines the total action for a growth path. Its minimizer represents the geodesic path in the morphoelastic shape space. One can verify that $\|\tens{G}\|_{[\Omega]}$ defines a rigorous norm, and the above minimization problem defines a valid geodesic distance $\mathcal{D}$ (see \cref{appendix:proofs}), with the solution guaranteed to exist under suitable topological assumptions \cite{charon2023shape}.

\paragraph{Numerical Equivalence to the Forward Variational Problem}

Leaving the complete solution of the full inverse optimal control problem for future work, we address here the specific subproblem that specifies a Lagrangian model for the evolution of the growth tensor $\tens{G}$, and extract the resulting forward shape evolution. Given that the velocity field solves the instantaneous variational problem in \cref{eq:growth_norm}, we precisely let $\tens{G}(t) = \tens{G}_0 \circ \phi(t)^{-1}$ (which maps exactly to the tracking of $\mathfrak{g}(\vec{x},t)$ defined in \cref{eq:lagrangian_current}), where $\partial_t \phi(t) = \vec{v}(t)\circ \phi(t)$ and $\vec{v}(t)$ solves \cref{eq:growth_norm}.

\section{Variational Formulation with Smoothing Regularization}
\label{sec:variational_regularized}

To solve for the displacement field $\vec{u}$ resulting from the growth tensor $\tens{G}$, we formulate the problem within the framework of the Principle of Minimum Potential Energy. In addition to the elastic energy arising from the material's response to growth, we introduce a regularization term to ensure the smoothness of the resulting displacement field. This is crucial for obtaining physically plausible and well-behaved solutions, especially in numerical implementations.

This approach can be seen as a simplified, static version of more comprehensive dynamic models in computational anatomy, where smoothness is enforced on a time-dependent velocity field to generate geodesic flows of diffeomorphisms \cite{beg2005computing}. In our static context, penalizing the spatial derivatives of the displacement field $\vec{u}$ serves a similar purpose: it ensures that we have a smooth evolution.

\subsection{The Regularized Energy Functional}

The total potential energy functional, $E_{\text{reg}}(\vec{u})$, is composed of the elastic strain energy and a smoothing regularization term.

We first analyze the elastic strain term. The internal energy stored in the body is due to the elastic strain, $\tens{\varepsilon}_g = \tens{\varepsilon}(\vec{u}) - \tens{G}$. The strain energy density, $\Psi$, for an isotropic linear elastic material is given by \cite{sadd2020elasticity}:
\begin{equation}
	\Psi(\tens{\varepsilon}_g) = \frac{\lambda}{2} \text{tr}(\tens{\varepsilon}_g)^2 + \mu \text{tr}(\tens{\varepsilon}_g^2)
\end{equation}
where $\lambda$ and $\mu$ are the Lamé constants.

 To enforce smoothness, we add a penalty term proportional to the squared $L^2$-norm of the displacement field with an operator. This term penalizes sharp curvatures in $\vec{u}$ and yields a smoother transformation. The regularization energy is:
\begin{equation}
	E_{\text{reg}}(\vec{u}) = \frac{\kappa}{2} \int_{\Omega} \| L \vec{u}\|^2 \, dV
\end{equation}
where $\kappa > 0$ is a regularization parameter that controls the trade-off between accommodating the growth and maintaining the smoothness of the displacement, and $L$ is a linear operator acting on $\vec{u}$.

\paragraph{Conclusion of Total Potential Energy.} The total energy functional to be minimized is given as follows (external loading terms set to zero, as we focus on the effect of internal growth):
\begin{equation}
	E_{\text{total}}(\vec{u}) = \int_{\Omega} \Psi(\tens{\varepsilon}(\vec{u}) - \tens{G}) \, dV + \frac{\kappa}{2} \int_{\Omega} \|L \vec{u}\|^2 \, dV
	\label{eq:total_energy_reg_new}
\end{equation}

\subsection{Applying the Sobolev Embedding Theorem}

In the context of three-dimensional shape analysis $(d=3)$, it is desirable for the displacement field $\vec{u}$ to have at least $C^1$ regularity. This is because the discrete deformation map is written as
\[
    \vec{\phi}(\vec{x})=\vec{x}+\vec{u}(\vec{x}),
\]
and a $C^1$ displacement field ensures that the deformation gradient
\[
    \nabla \vec{\phi}(\vec{x})
    =
    \mathbf{I}+\nabla \vec{u}(\vec{x})
\]
is pointwise well-defined and continuous. This regularity is important for obtaining a smooth and physically meaningful deformation. Strictly speaking, $C^1$ regularity alone does not guarantee that $\vec{\phi}$ is a global diffeomorphism; additional conditions, such as non-degeneracy of the Jacobian determinant and global injectivity, are needed. Nevertheless, $C^1$ regularity is a natural minimal smoothness requirement in this setting.

According to the Sobolev Embedding Theorem \cite{arbogast2025functional}, the embedding into H\"older or continuous spaces depends on the relation $mp>d$. Specifically, for a bounded domain $\Omega\subset\mathbb{R}^3$ with a bounded extension operator, Part (d) of Theorem 7.22 in \cite{arbogast2025functional} states that
\begin{equation}
    W^{j+m,p}(\Omega)\hookrightarrow C^j(\Omega)
    \qquad \text{if } mp>d.
\end{equation}
In our implementation, we consider the Hilbert case $p=2$. To obtain $C^1$ regularity, we set $j=1$ and require
\begin{equation}
    2m>3.
\end{equation}
Thus $m>3/2$, and since $m$ must be an integer in the theorem, we take $m\ge 2$. Consequently,
\begin{equation}
    \vec{u}\in W^{1+2,2}(\Omega)=H^3(\Omega)
\end{equation}
is sufficient to guarantee the desired $C^1$ regularity. This motivates the use of a higher-order regularization operator in the variational formulation, since insufficient regularity may lead to non-smooth or potentially overlapping transformations in large deformation settings.

One point worth mentioning is the relation between Laplacian-type penalties and Sobolev norms. If the regularization operator is chosen formally as
\[
    L=\Delta^\beta,
\]
then the penalty
\[
    \|\Delta^\beta \vec{u}\|_{L^2(\Omega)}^2
\]
controls high-order derivatives of the displacement field. Under suitable elliptic boundary conditions and elliptic regularity assumptions, such a penalty is equivalent, up to lower-order terms, to controlling the full $H^{2\beta}(\Omega)$ norm. This equivalence is commonly justified by integration by parts, Poincar\'e-type inequalities, and elliptic regularity estimates \cite{evans2022partial}.

However, in the present elasticity model, the displacement field is not assumed to vanish on the entire boundary $\partial\Omega$. Instead, the boundary is decomposed as
\[
    \partial\Omega=\Gamma_D\cup\Gamma_N,
\]
where the homogeneous Dirichlet condition
\[
    \vec{u}=\vec{0}
    \qquad \text{on } \Gamma_D
\]
is imposed only on the fixed part of the boundary, while the remaining part $\Gamma_N$ is governed by natural Neumann-type boundary conditions. Therefore, the norm-equivalence statement should be understood under the corresponding homogeneous displacement space or under boundary conditions that make the elliptic operator well-posed. The main purpose of the operator penalty is to impose Sobolev-type smoothness on the displacement field and to suppress irregular or highly oscillatory deformations.

\subsection{Rigorous Operator Selection in LDDMM}

For a more mathematically rigorous treatment of the regularization, we refer to the Large Deformation Diffeomorphic Metric Mapping (LDDMM) framework as detailed in \cite{beg2005computing}. To guarantee that the flow of the velocity field results in a valid diffeomorphism, the differential operator $L$ is typically chosen to be of the Cauchy-Navier type:
\begin{equation}
	L = (-\alpha \Delta + \gamma)^\beta I
\end{equation}
where $\alpha, \gamma > 0$ and $I$ is the identity operator. In three-dimensional space ($d=3$), \cite{beg2005computing} specifies that the power $\beta$ must satisfy $\beta > 1.5$ to ensure sufficient Sobolev smoothness. This choice ensures that the operator $L^\dagger L$ (which appears in the Euler-Lagrange equations, $L^\dagger$ is called the adjoint operator) is of order $2\beta > 3$. In numerical implementation the lowest order we can choose is $\beta = 2$.

\subsubsection{General Operator and Strong Form}

To establish a generalized framework, we consider a regularization energy
defined by the squared \(L^2\)-norm of a differential operator applied to the
displacement field \(\mathbf u\). Let
\[
    A=-\alpha\Delta+\gamma I,
    \qquad
    L=A^\beta,
\]
where \(\alpha>0\), \(\gamma\ge 0\), \(\beta\ge 1\) is an integer, and \(I\)
denotes the identity operator acting componentwise on vector fields. The
regularization energy is
\begin{equation}
    E_{\mathrm{reg}}(\mathbf u)
    =
    \frac{\kappa}{2}
    \int_\Omega
    \|L\mathbf u\|^2\,dV
    =
    \frac{\kappa}{2}
    \int_\Omega
    (A^\beta\mathbf u)\cdot(A^\beta\mathbf u)\,dV.
\end{equation}

To derive the associated Euler--Lagrange equation, we compute the first
variation with respect to a test function \(\mathbf h\in\hat V\). Let
\(\epsilon\) be a scalar perturbation. Then
\begin{align}
    \delta E_{\mathrm{reg}}(\mathbf u;\mathbf h)
    &=
    \frac{d}{d\epsilon}
    E_{\mathrm{reg}}(\mathbf u+\epsilon\mathbf h)
    \bigg|_{\epsilon=0} \nonumber\\
    &=
    \frac{\kappa}{2}
    \int_\Omega
    \frac{d}{d\epsilon}
    \left[
    L(\mathbf u+\epsilon\mathbf h)
    \cdot
    L(\mathbf u+\epsilon\mathbf h)
    \right]
    \bigg|_{\epsilon=0}
    \,dV \nonumber\\
    &=
    \kappa
    \int_\Omega
    (L\mathbf u)\cdot(L\mathbf h)\,dV \nonumber\\
    &=
    \kappa
    \int_\Omega
    (A^\beta\mathbf u)\cdot(A^\beta\mathbf h)\,dV.
\end{align}

To isolate the test function \(\mathbf h\) and identify the interior strong
form, we formally use the adjoint of the regularization operator. If the
boundary conditions are chosen so that all boundary terms generated by
integration by parts vanish, then \(A\) is self-adjoint with respect to the
\(L^2\)-inner product. Consequently,
\[
    L^\dagger L
    =
    (A^\beta)^\dagger A^\beta
    =
    A^\beta A^\beta
    =
    A^{2\beta}.
\]
Therefore, formally,
\begin{equation}
    \kappa
    \int_\Omega
    (A^\beta\mathbf u)\cdot(A^\beta\mathbf h)\,dV
    =
    \kappa
    \int_\Omega
    (A^{2\beta}\mathbf u)\cdot\mathbf h\,dV.
\end{equation}
Thus the corresponding interior strong-form regularization term is
\begin{equation}
    \kappa L^\dagger L\mathbf u
    =
    \kappa A^{2\beta}\mathbf u
    =
    \kappa(-\alpha\Delta+\gamma I)^{2\beta}\mathbf u.
\end{equation}
This is a differential operator of order \(4\beta\).

This strong-form expression should be understood as an interior or formal strong
form. In the actual variational problem, the high-order boundary conditions are
not imposed separately. Instead, they are encoded weakly by the choice of the
admissible displacement space, the test space, and the natural boundary
conditions arising from the variational principle. In particular, the physical
Dirichlet condition is imposed only on the fixed part \(\Gamma_D\), while the
remaining boundary \(\Gamma_N\) carries the corresponding Neumann-type natural
boundary conditions.

% To isolate the test function $\mathbf{h}$ and identify the strong form, we rely on the properties of the base operator $A = -\alpha \Delta + \gamma I$. Assuming appropriate boundary conditions where boundary integrals vanish (e.g., $\mathbf{h}$ and its derivatives vanish on $\partial \Omega$), the Laplacian $\Delta$ is a self-adjoint operator. Consequently, the operator $A$ is self-adjoint:
% \begin{equation}
% 	\int_{\Omega} (A \mathbf{v}) \cdot \mathbf{w} \, dV = \int_{\Omega} \mathbf{v} \cdot (A \mathbf{w}) \, dV.
% \end{equation}
% Since the operator $L$ is a power of a self-adjoint operator ($L = A^\beta$), it is also self-adjoint. Specifically, applying integration by parts $\beta$ times transfers the operator $L$ from $\mathbf{h}$ to $\mathbf{u}$:
% \begin{align}
% 	\kappa \int_{\Omega} ( A^\beta \mathbf{u}) \cdot ( A^\beta \mathbf{h}) \, dV &= \kappa \int_{\Omega} ( A^\beta ( A^\beta \mathbf{u}))\cdot \mathbf{h}  \, dV \quad (\text{Self-adjointness}) \nonumber \\
% 	&= \int_{\Omega} (\kappa A^{2\beta} \mathbf{u}) \cdot \mathbf{h} \, dV
% \end{align}

% Thus, the variation of the regularization energy leads to the following term in the strong form of the equilibrium equation:
% \begin{equation}
% 	\kappa L^{\dagger}L \mathbf{u} = \kappa L^{2} \mathbf{u} = \kappa (-\alpha \Delta + \gamma)^{2\beta} \mathbf{u}
% \end{equation}
% This term represents a differential operator of order $4\beta$.

\subsection{Summary of the New Variational Problem}
The modified problem is to find the displacement field $\vec{u} \in V$ such that for all test functions $\vec{h} \in \hat{V}$:
\begin{equation}
	a_{\text{reg}}(\vec{u}, \vec{h}) = L_{\text{growth}}(\vec{h})
\end{equation}
where the bilinear form $a_{\text{reg}}(\cdot, \cdot)$ and the linear form $L_{\text{growth}}(\cdot)$ are defined as:
\begin{align}
	a_{\text{reg}}(\vec{u}, \vec{h}) &= \int_{\Omega} \tens{\sigma(\varepsilon(\vec{u})}) : \nabla\vec{h} \, dV + \kappa \int_{\Omega} \left[ (-\alpha \Delta + \gamma)^\beta \vec{u} \right] \cdot \left[ (-\alpha \Delta + \gamma)^\beta \vec{h} \right] \, dV \\
	L_{\text{growth}}(\vec{h}) &= \int_{\Omega} \tens{\sigma}(\tens{G}) : \nabla\vec{h} \, dV
\end{align}
The addition of the second term in $a_{\text{reg}}(\vec{u}, \vec{h})$ transforms the problem from a standard second-order elliptic system to a high-order system of order $4\beta$ (\cref{appendix:general_strong_form}). This formulation naturally yields smoother solutions for the displacement field $\vec{u}$, with the degree of regularity explicitly controlled by the parameter $\beta$. In comparison with the purely elastic formulation in \cref{sec:growth_model} (\cref{eq:variation_Growth}), this generalized variational principle ensures the smoothness required for large diffeomorphic deformations.

\subsubsection{Numerical Setting of Generalized Mixed Variational Formulation}

A significant numerical challenge arises when discretizing the generalized high-order regularization term derived in the previous subsection. Directly solving a partial differential equation of order \(4\beta\) would require finite element spaces with \(C^{2\beta-1}\) global continuity, which is computationally prohibitive for standard \(C^0\) Lagrange finite elements.

To circumvent this restriction, we employ a cascadic mixed finite element formulation. The main idea is to introduce auxiliary fields representing successive applications of the second-order elliptic operator, thereby decomposing the original high-order problem into a coupled cascade of second-order weak equations. The fundamental well-posedness, uniqueness, and stability of mixed formulations are established by Brezzi in \cite{M2AN_1974__8_2_129_0}. A classical foundational application is the mixed formulation for the fourth-order biharmonic equation, as developed by Ciarlet and Raviart in \cite{CIARLET1974125}. Further studies for more general high-order elliptic problems are provided in \cite{Bramble1985}.

Let
\[
    A=-\alpha\Delta+\gamma I,
\]
where \(I\) denotes the identity operator acting componentwise on vector fields. The strong-form regularization term derived formally in the previous subsection is
\[
    \kappa A^{2\beta}\mathbf u.
\]
To reduce this high-order operator to a sequence of second-order problems, we introduce auxiliary variables \(\mathbf u_k\), for \(k=0,1,\dots,2\beta-1\), with
\[
    \mathbf u_0=\mathbf u,
\]
where \(\mathbf u_0\) is the primary physical displacement field. The auxiliary variables are defined recursively by
\begin{align}
    \mathbf u_0 &= \mathbf u,\\
    \mathbf u_k &= A\mathbf u_{k-1}
    =
    (-\alpha\Delta+\gamma I)\mathbf u_{k-1},
    \qquad k=1,\dots,2\beta-1.
\end{align}
By substitution,
\[
    \mathbf u_k=A^k\mathbf u.
\]
Thus, setting \(N=2\beta\), the full regularization term can be recovered through the highest-order auxiliary field:
\begin{equation}
    \kappa A^{2\beta}\mathbf u
    =
    \kappa A\mathbf u_{N-1}.
\end{equation}

The primary displacement field satisfies the physical Dirichlet boundary condition on the fixed boundary, while the auxiliary variables are introduced only to reduce the differential order of the problem. Therefore, we distinguish the displacement space from the auxiliary spaces. Let
\[
    \mathcal V_D
    =
    \left\{
    \mathbf v\in [H^1(\Omega)]^d
    \mid
    \mathbf v=\mathbf 0 \text{ on } \Gamma_D
    \right\},
\]
and let
\[
    \mathcal W=[H^1(\Omega)]^d.
\]
Then the mixed space is taken to be
\begin{equation}
    \mathcal V_{\mathrm{mixed}}
    =
    \mathcal V_D\times \mathcal W^{N-1}.
\end{equation}
This choice imposes the physical essential boundary condition only on the primary displacement \(\mathbf u_0\), while avoiding artificial essential boundary conditions on the auxiliary fields \(\mathbf u_1,\dots,\mathbf u_{N-1}\).

The mixed problem is: find
\[
    (\mathbf u_0,\mathbf u_1,\dots,\mathbf u_{N-1})
    \in
    \mathcal V_{\mathrm{mixed}}
\]
such that for all test functions
\[
    (\boldsymbol{\psi}_0,\boldsymbol{\psi}_1,\dots,\boldsymbol{\psi}_{N-1})
    \in
    \mathcal V_{\mathrm{mixed}},
\]
the following coupled system is satisfied.

\begin{enumerate}
    \item \textbf{Mechanical Equilibrium.}
    
    The first equation balances the internal elastic stress, the regularization penalty expressed through the auxiliary variable \(\mathbf u_{N-1}\), and the growth-induced stress:

    \begin{equation}
		\int_\Omega \sigma(\mathbf{u}_0) : \nabla \boldsymbol{\psi}_0 \, dV 
		+ \kappa \int_\Omega \left( \alpha \nabla \mathbf{u}_{N-1} : \nabla \boldsymbol{\psi}_0 + \gamma \mathbf{u}_{N-1} \cdot \boldsymbol{\psi}_0 \right) dV 
		= \int_\Omega \sigma(G) : \nabla \boldsymbol{\psi}_0 \, dV
	\end{equation}
    
    \item \textbf{Recursive Constraints.}
    
    For each \(k=1,\dots,N-1\), we weakly enforce the hierarchical relation
    \[
        \mathbf u_k=A\mathbf u_{k-1}.
    \]
    Multiplying by the test function \(\boldsymbol{\psi}_k\in\mathcal W\) and applying integration by parts gives
    \begin{equation}
        \int_\Omega
        \mathbf u_k\cdot \boldsymbol{\psi}_k
        \,dV
        -
        \int_\Omega
        \left(
        \alpha\nabla \mathbf u_{k-1}:\nabla \boldsymbol{\psi}_k
        +
        \gamma \mathbf u_{k-1}\cdot \boldsymbol{\psi}_k
        \right)
        dV
        =
        0.
    \end{equation}
    
    The boundary terms generated in this integration-by-parts step are interpreted as the natural boundary conditions associated with the auxiliary equations and are handled weakly in the mixed variational formulation. Thus, no additional artificial essential boundary conditions are imposed on the auxiliary variables.
\end{enumerate}

\section{Simulation Examples}
\subsection{For $\beta = 2$}
For the lowest order case $\beta=2$, set $\kappa$ = 1. We have $L = (-\alpha \Delta + \gamma)^2$. Consequently, the operator acting on $\mathbf{u}$ in the strong form is $L^2 = (-\alpha \Delta + \gamma)^4$. Expanding this term using the binomial theorem reveals the high-order nature of the regularization:
\begin{equation}
	L^2 = \alpha^4 \Delta^4 - 4\alpha^3 \gamma \Delta^3 + 6\alpha^2 \gamma^2 \Delta^2 - 4\alpha \gamma^3 \Delta + \gamma^4 I
\end{equation}

To implement the rigorous LDDMM regularization derived in the previous section, where the operator is defined as $L = (-\alpha\Delta + \gamma)^2 I$, we must address the resulting eighth-order differential term $L^2 u = (-\alpha\Delta + \gamma)^4 u$. To solve this using standard $C^0$ Lagrangian finite elements, we apply the Mixed Finite Element Method, by introducing three auxiliary variables $\{w, z, p\}$. These variables represent successive applications of the second-order operator $(-\alpha\Delta + \gamma)$:

\begin{align}
	w &= (-\alpha\Delta + \gamma)u \\
	z &= (-\alpha\Delta + \gamma)w \\
	p &= (-\alpha\Delta + \gamma)z 
\end{align}

Under this decomposition, the eighth-order term is recovered as $(-\alpha\Delta + \gamma)p$. The resulting coupled variational problem seeks $(u, w, z, p) \in \mathcal{V}_{mixed}$ such that for all test functions $(h, \eta, \xi, q) \in \mathcal{V}_{mixed}$:

\begin{subequations}
	\label{eq:mixed_lddmm}
	\begin{align}
		\int_\Omega \sigma(u) : \nabla h \, dV + \alpha \int_\Omega \nabla p \cdot \nabla h \, dV + \gamma \int_\Omega p \cdot h \, dV &= \int_\Omega \sigma(G) : \nabla h \, dV \\
		\int_\Omega p \cdot q \, dV - \alpha \int_\Omega \nabla z \cdot \nabla q \, dV - \gamma \int_\Omega z \cdot q \, dV &= 0 \\
		\int_\Omega z \cdot \xi \, dV - \alpha \int_\Omega \nabla w \cdot \nabla \xi \, dV - \gamma \int_\Omega w \cdot \xi \, dV &= 0 \\
		\int_\Omega w \cdot \eta \, dV - \alpha \int_\Omega \nabla u \cdot \nabla \eta \, dV - \gamma \int_\Omega u \cdot \eta \, dV &= 0
	\end{align}
\end{subequations}

In this system:
	\cref{eq:mixed_lddmm}a enforces the mechanical equilibrium, where the LDDMM smoothing term is used to ensure high-order regularity of the displacement field. Equations \cref{eq:mixed_lddmm}b--\cref{eq:mixed_lddmm}d weakly enforce the hierarchical constraints that link the auxiliary variables back to the displacement $u$. By applying integration by parts, we transfer the Laplacian operators onto the test functions, allowing the use of $H^1$-conforming elements (such as $P_2$ elements) for all fields. This avoids the need for complex $C^3$ continuous elements.

As for the numerical implementation, we employ the direct solver MUMPS, available through FEniCSx; further details can be found in \cite{doi:10.1137/S0895479899358194}.

\subsection{Example 5}

\begin{center}
    \includegraphics[width=0.6\textwidth]{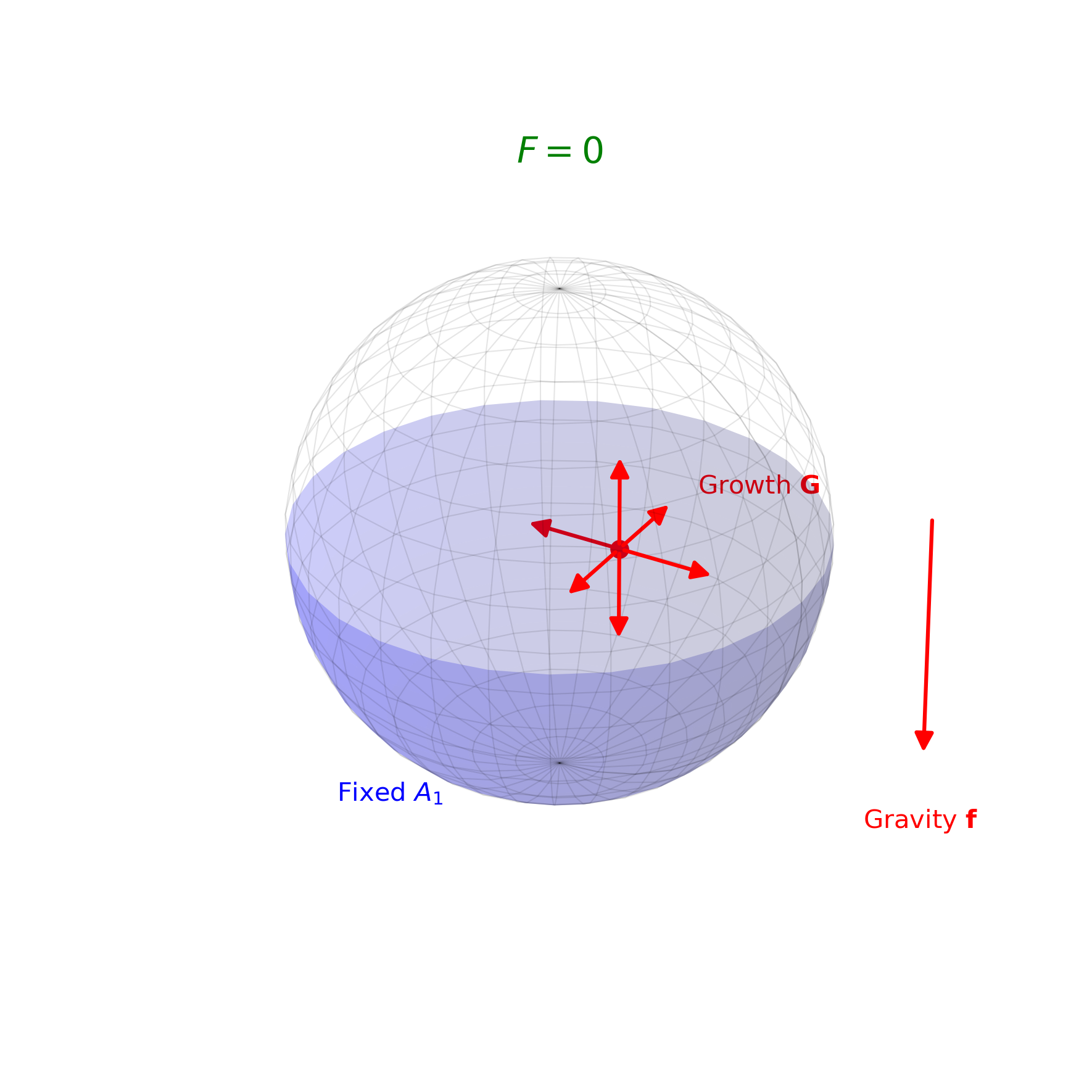}
    \captionof{figure}{Schematic for Example 5}
    \label{fig:sketch_model_8}
\end{center}

\begin{itemize}
	\item \textbf{Geometry:} A ball centered at $ (0,0,0) $ with radius 0.5.
    \item \textbf {Growth:} The growth center is $(0.25,0,0)$, with growth coefficient $A_g$ = 0.5 and $\tau$ = 0.1. 
	\item \textbf{Body Force:} $\vec{f} = (0, 0, 0)$.
	\item \textbf{Dirichlet BC ($\Gamma_D$):} The lower hemisphere is fixed.
	\item \textbf{Neumann BC ($\Gamma_N$):} $\vec{T} = 0$ on $A_2$, where $A_2$ = $\Gamma_N = \Gamma \setminus \Gamma_D$.
\end{itemize}

\begin{figure}[H] 
	\centering
	
	\setlength{\abovecaptionskip}{3pt} 
	\setlength{\belowcaptionskip}{0pt}
	
	\begin{subfigure}[b]{0.48\textwidth}
		\centering
		\includegraphics[width=\textwidth]{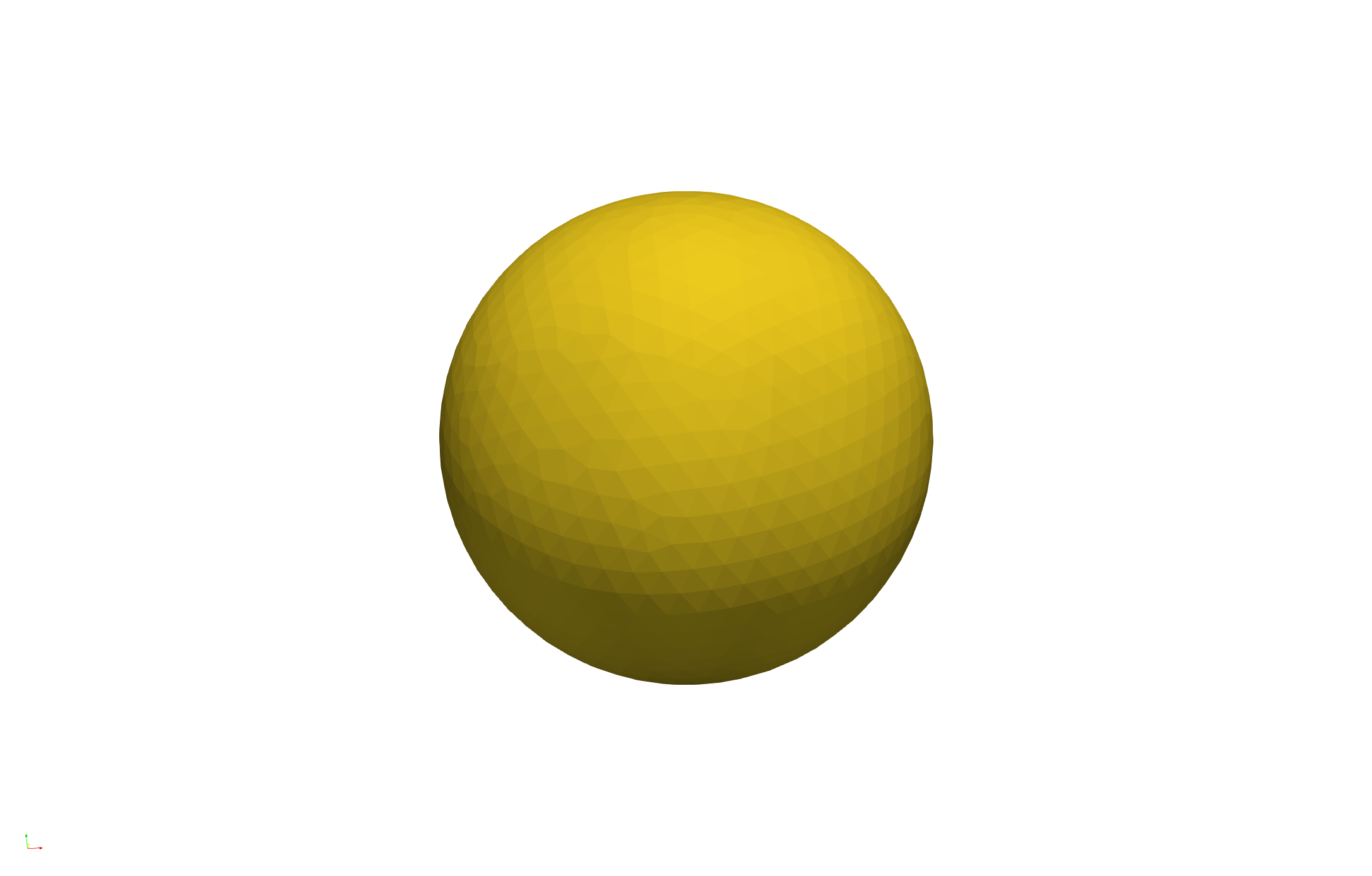} 
		\caption{Initial: Iteration 0}
	\end{subfigure}
	\hfill 
	\begin{subfigure}[b]{0.48\textwidth}
		\centering
		\includegraphics[width=\textwidth]{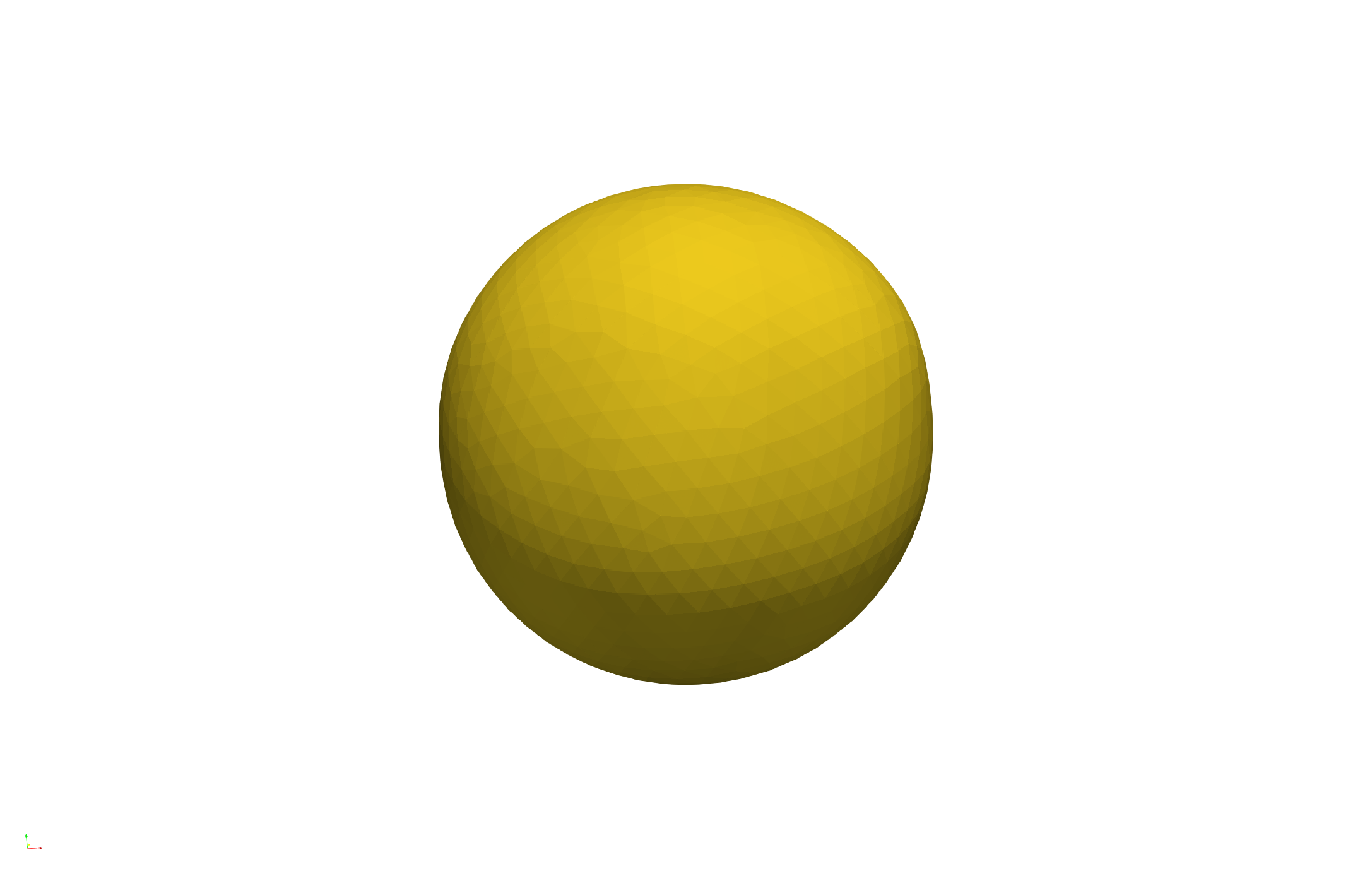}
		\caption{Iteration 1}
	\end{subfigure}
	
	\vspace{0.3em} 
	
	\begin{subfigure}[b]{0.48\textwidth}
		\centering
		\includegraphics[width=\textwidth]{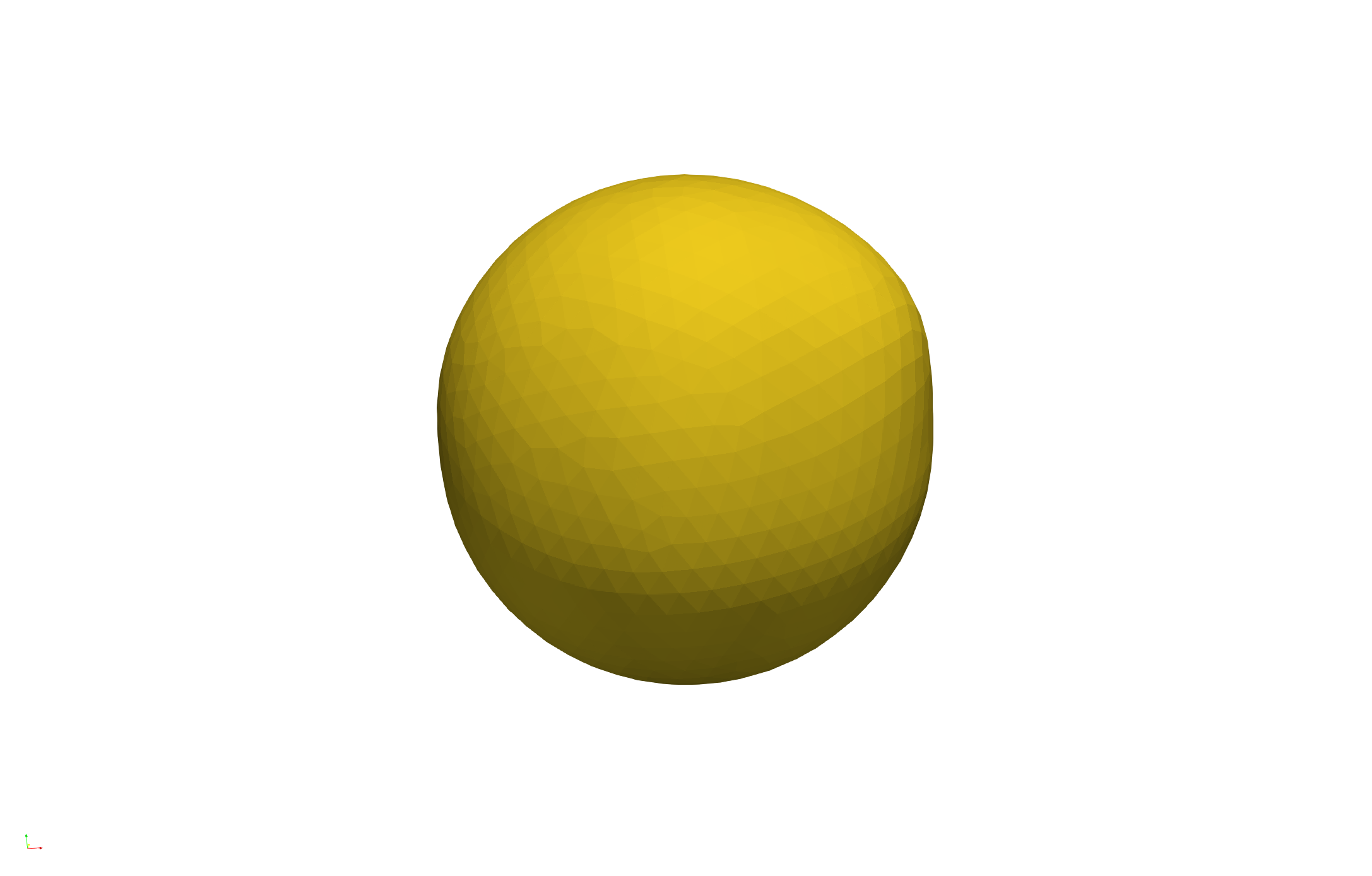} 
		\caption{Iteration 2}
	\end{subfigure}
	\hfill 
	\begin{subfigure}[b]{0.48\textwidth}
		\centering
		\includegraphics[width=\textwidth]{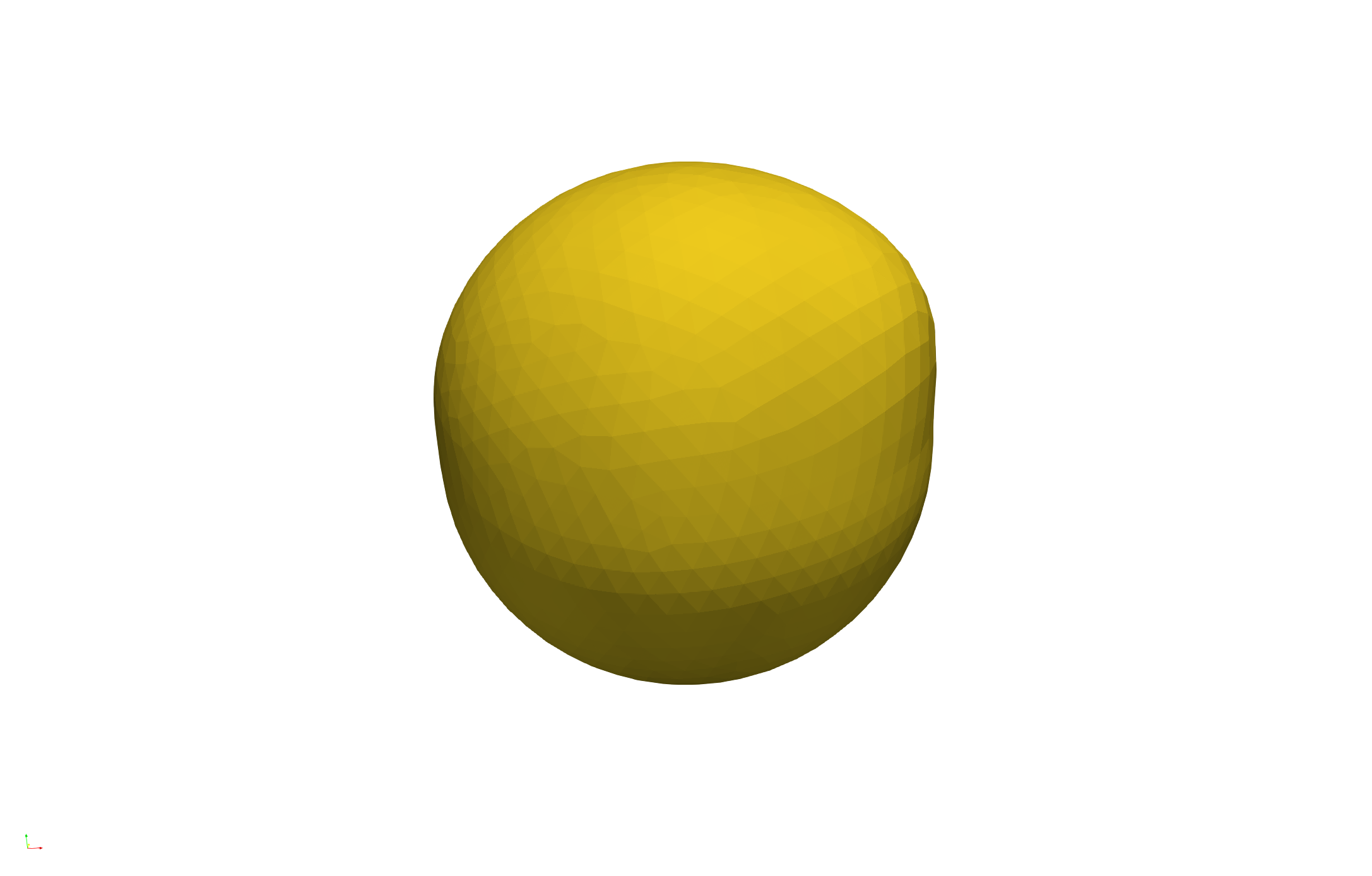}
		\caption{Iteration 3}
	\end{subfigure}
	
	\vspace{0.3em}
	
	\begin{subfigure}[b]{0.48\textwidth}
		\centering
		\includegraphics[width=\textwidth]{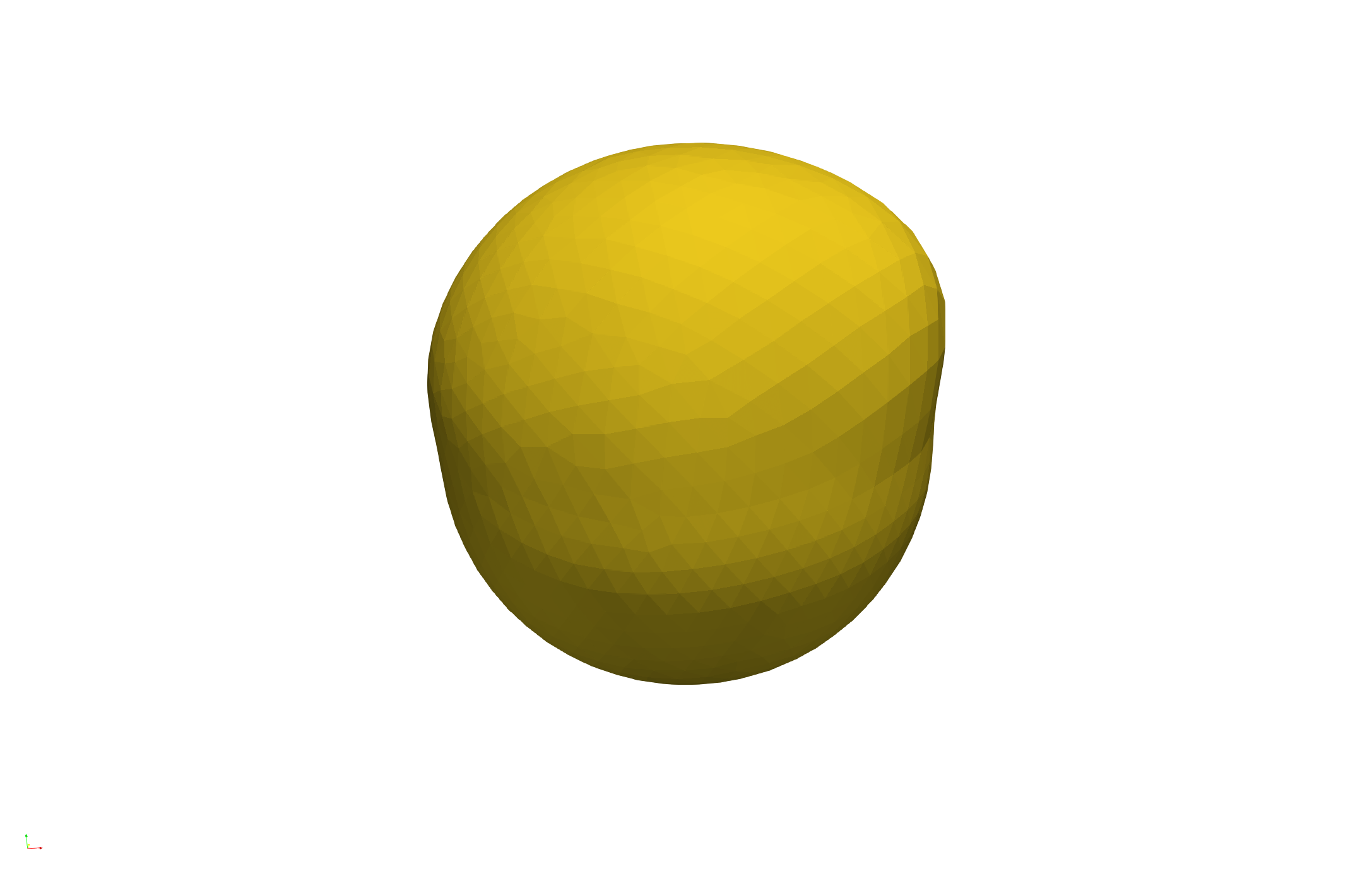} 
		\caption{Iteration 4}
	\end{subfigure}
	
	\caption[Numerical results of the sphere across different iterations]{ Numerical results across integration steps. Here we slightly increased the mesh size for computational reasons, using the LDDMM operator with $\beta = 2$ and set the parameters $\alpha$ and $\gamma$ equal to 0.01 and 0.001, $\kappa$ = 1.}
	\label{fig:my_five_plots}
\end{figure}

\clearpage

\subsection{Comparison between No Regularization and An LDDMM Regularization}

\begin{figure} [H]
	\centering
	
	\setlength{\abovecaptionskip}{5pt} 
	\setlength{\belowcaptionskip}{0pt}
	
	\begin{subfigure}[b]{0.45\textwidth}
		\centering
		\includegraphics[width=\textwidth]{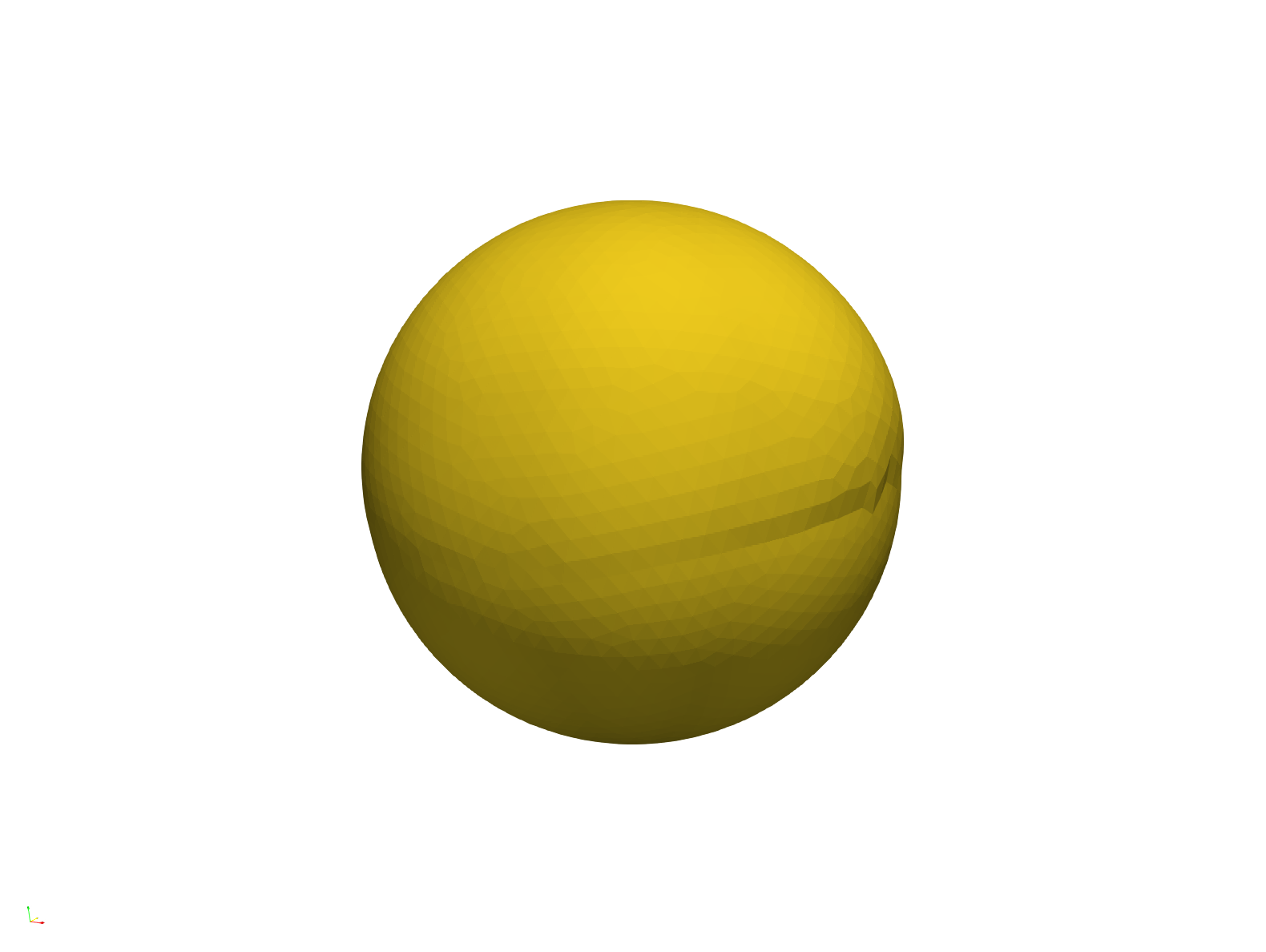} 
		\caption{No Reg: Iteration 1}
	\end{subfigure}
	\hfill 
	\begin{subfigure}[b]{0.45\textwidth}
		\centering
		\includegraphics[width=\textwidth]{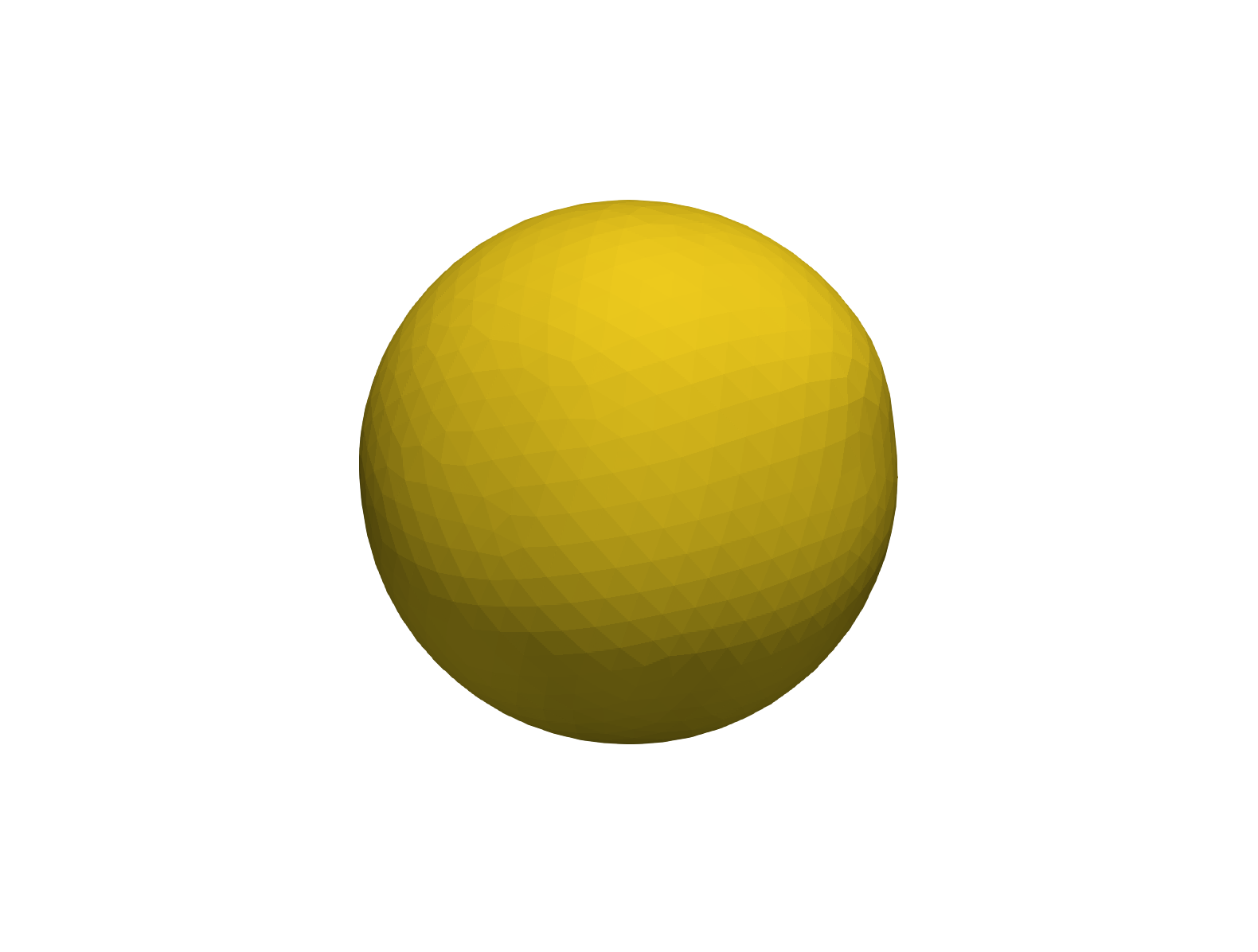}
		\caption{With Reg: Iteration 1}
	\end{subfigure}
	
	\vspace{0.3em} 
	
	\begin{subfigure}[b]{0.45\textwidth}
		\centering
		\includegraphics[width=\textwidth]{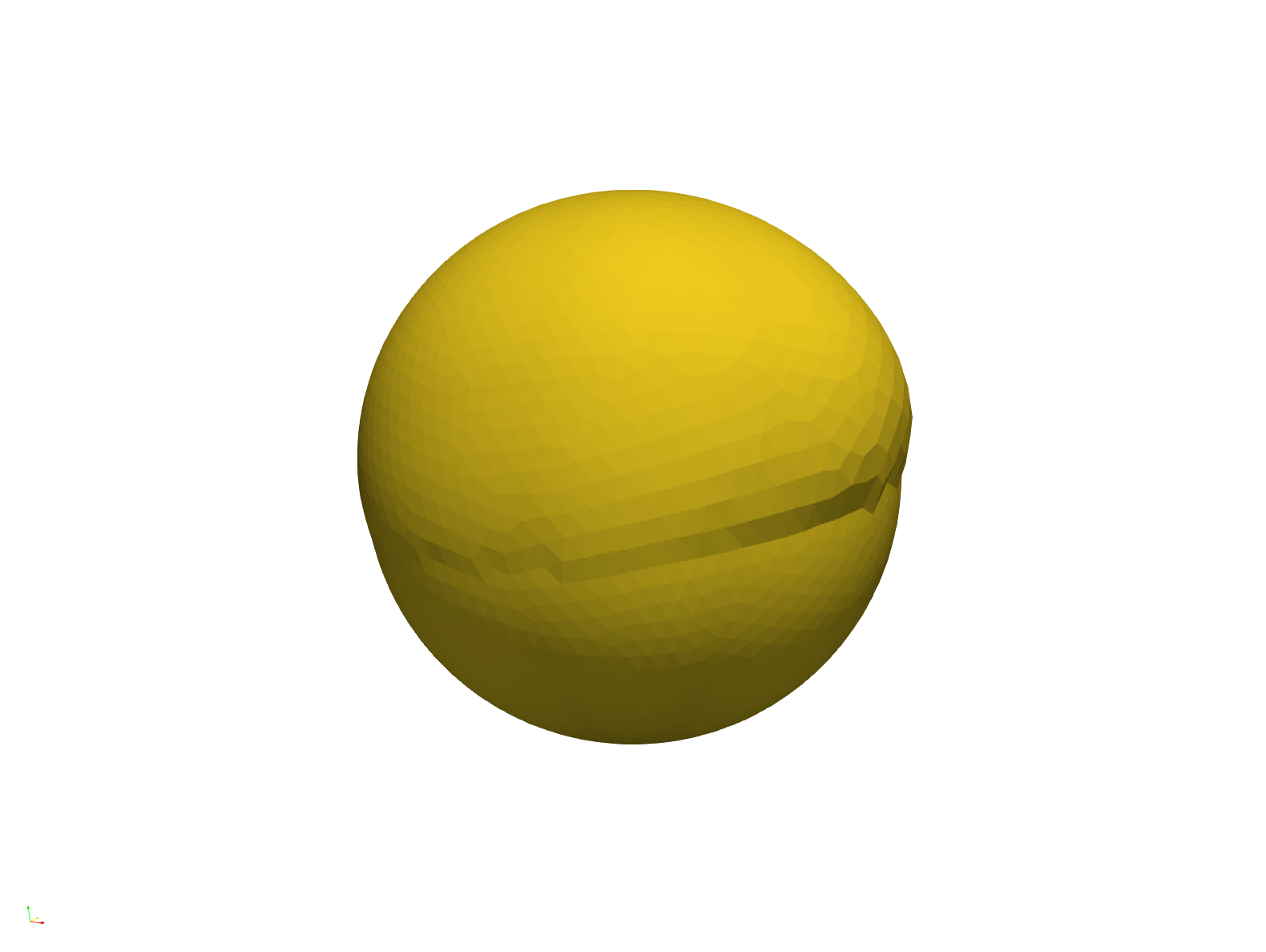} 
		\caption{No Reg: Iteration 2}
	\end{subfigure}
	\hfill 
	\begin{subfigure}[b]{0.45\textwidth}
		\centering
		\includegraphics[width=\textwidth]{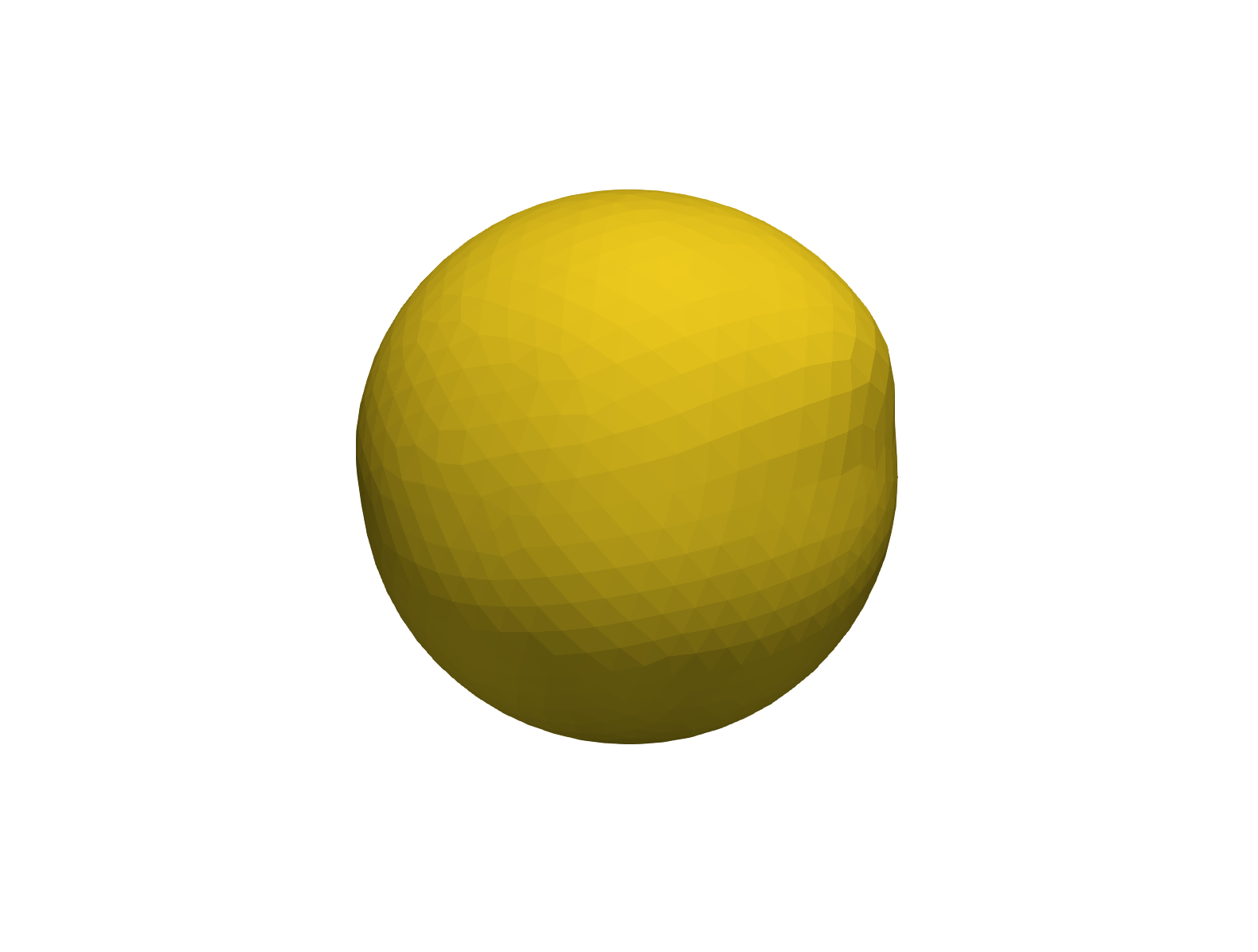}
		\caption{With Reg: Iteration 2}
	\end{subfigure}
	
	\vspace{0.3em}
	
	\begin{subfigure}[b]{0.45\textwidth}
		\centering
		\includegraphics[width=\textwidth]{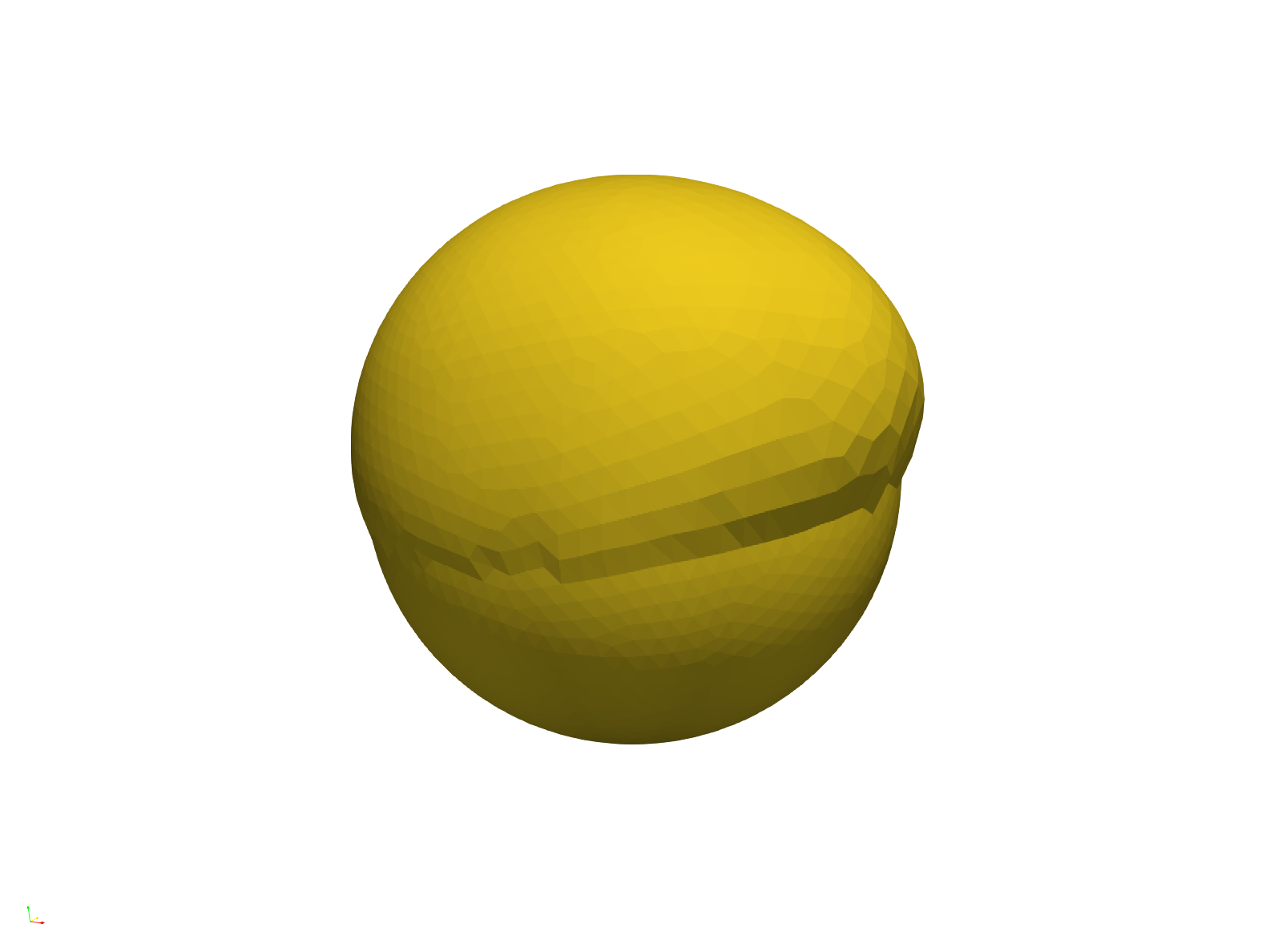} 
		\caption{No Reg: Iteration 3}
	\end{subfigure}
	\hfill 
	\begin{subfigure}[b]{0.45\textwidth}
		\centering
		\includegraphics[width=\textwidth]{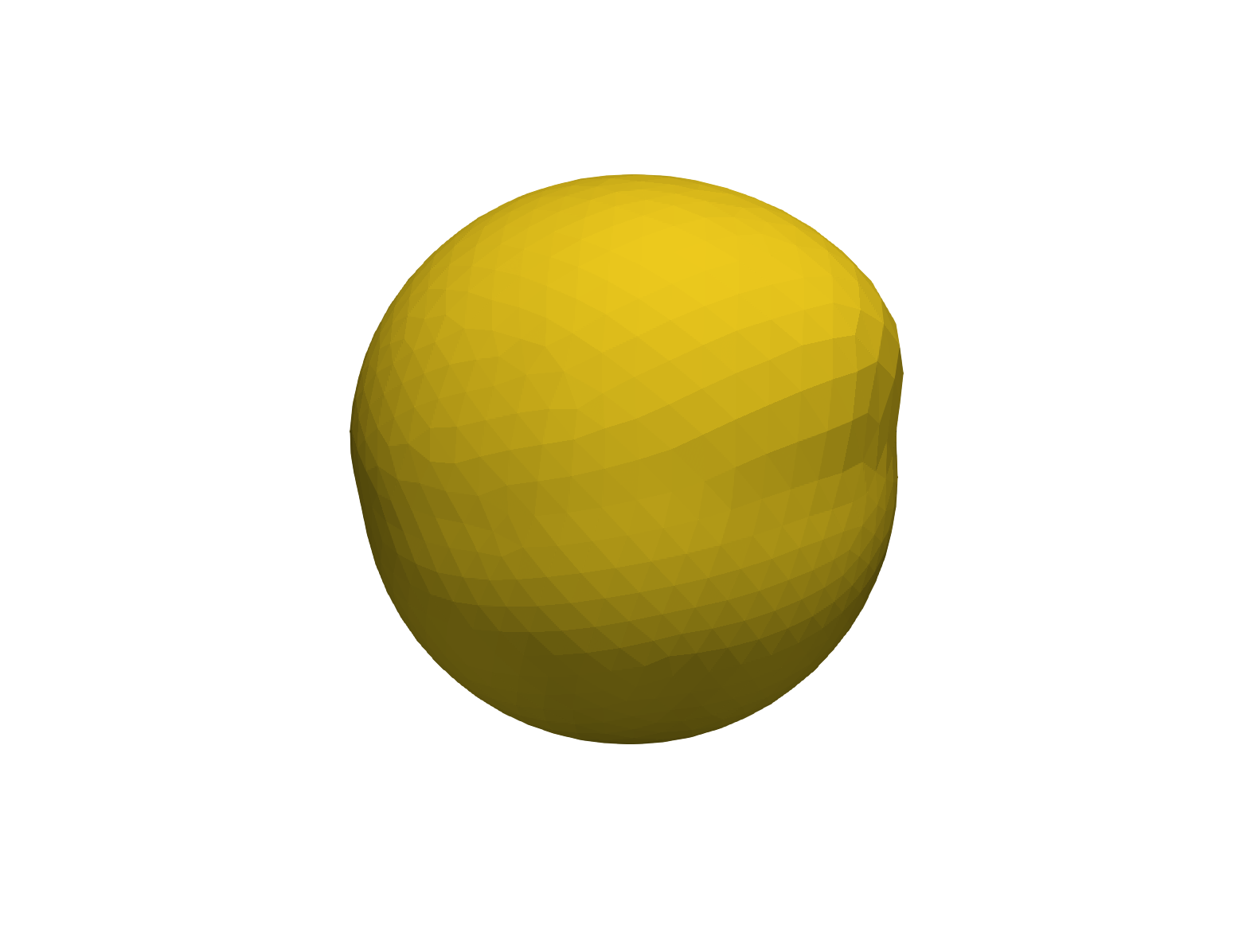}
		\caption{With Reg: Iteration 3}
	\end{subfigure}

	\caption[Step-by-step comparison for Example 5 with and without regularization]{Step-by-step comparison for Example 5. Left: Without regularization. Right: With an LDDMM operator as the regularization term. Set the parameters $\alpha$ and $\gamma$ equal to 0.01 and 0.001, $\kappa$ = 0.1  }
	\label{fig:compare_m8_full}
\end{figure}
\clearpage

\subsection{Extended Simulation Scenarios}

 While the beams and spheres serve as a fundamental benchmark, biological structures often exhibit more complicated shapes, some even with non-convex or non-simply connected features. In the following experiments, we maintain the mixed variational formulation derived in the previous section but vary the domain geometry $\Omega$ and boundary conditions. All the examples below use the LDDMM operator with $\beta = 2$ and set the parameters $\alpha$ and $\gamma$ equal to 0.01 and 0.001, $\kappa$ = 1.

\subsubsection{Sphere with Localized Boundary Constraint}
\begin{figure}[H] 
	\centering
	\setlength{\abovecaptionskip}{3pt} 
	\setlength{\belowcaptionskip}{0pt}
	
	\begin{subfigure}[b]{0.32\textwidth}
		\centering
		\includegraphics[width=\textwidth]{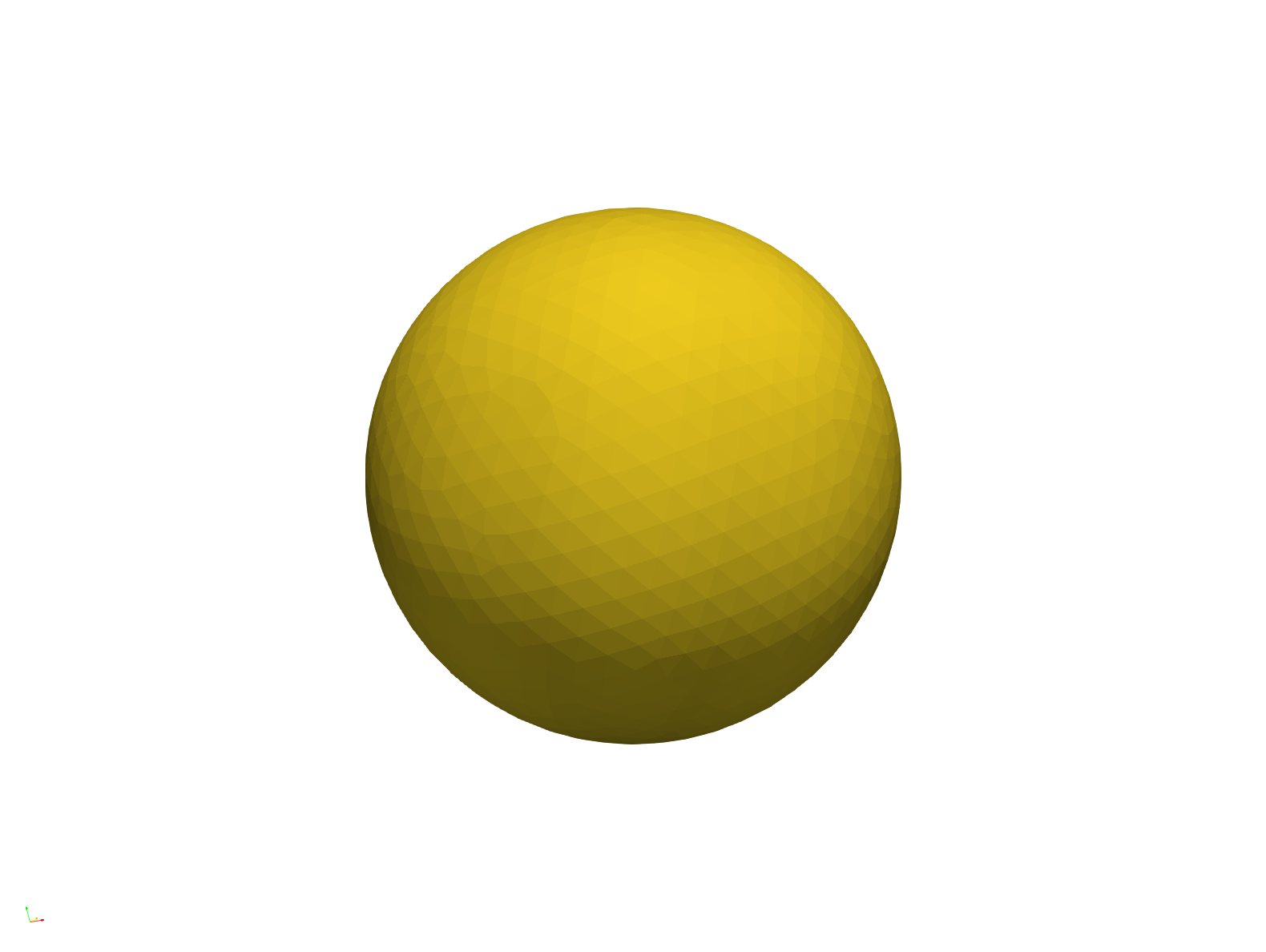} 
		\caption{Initial Shape}
	\end{subfigure}
	\hfill 
	\begin{subfigure}[b]{0.32\textwidth}
		\centering
		\includegraphics[width=\textwidth]{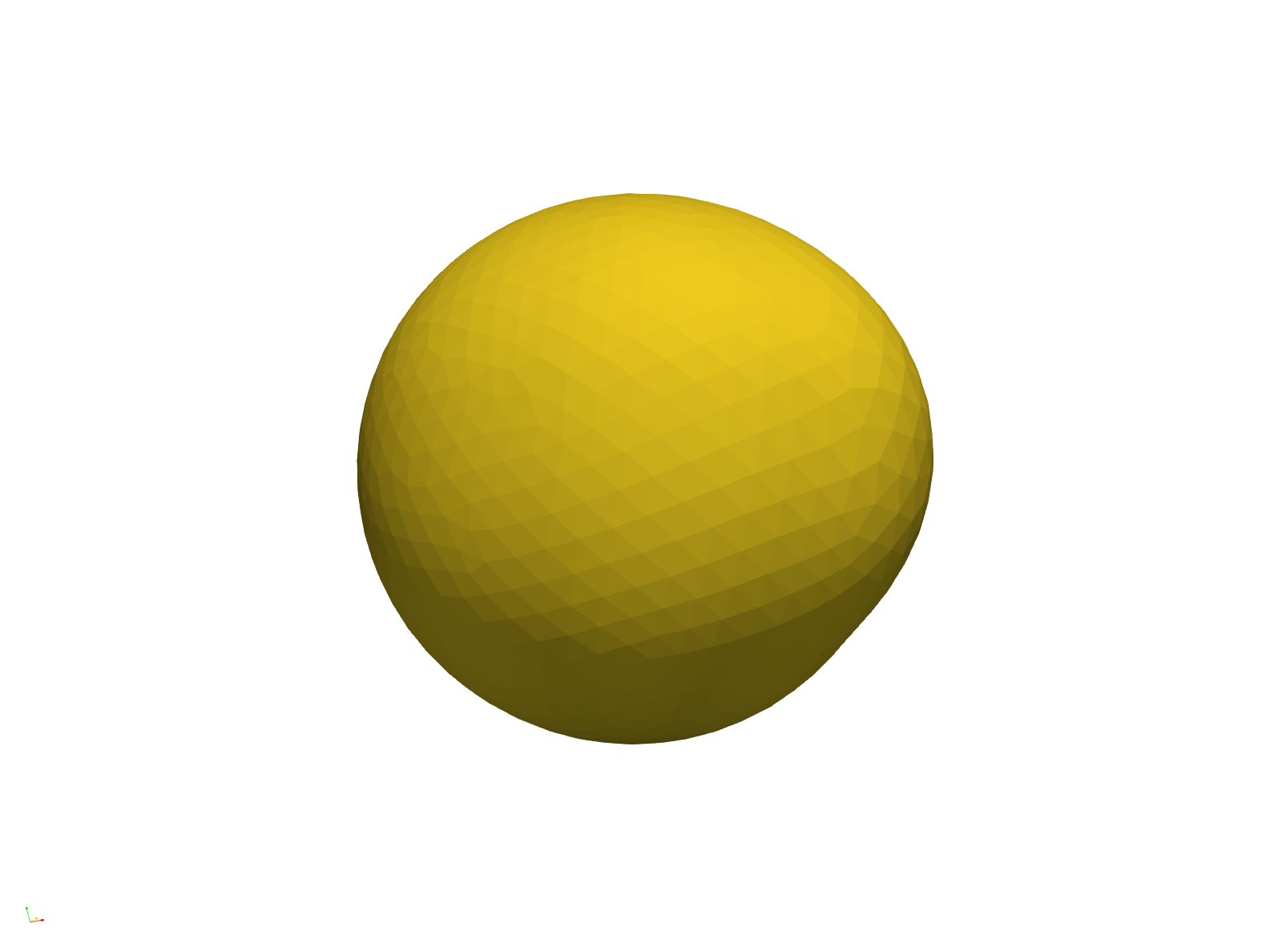} 
		\caption{Iteration 2}
	\end{subfigure}
	\hfill 
	\begin{subfigure}[b]{0.32\textwidth}
		\centering
		\includegraphics[width=\textwidth]{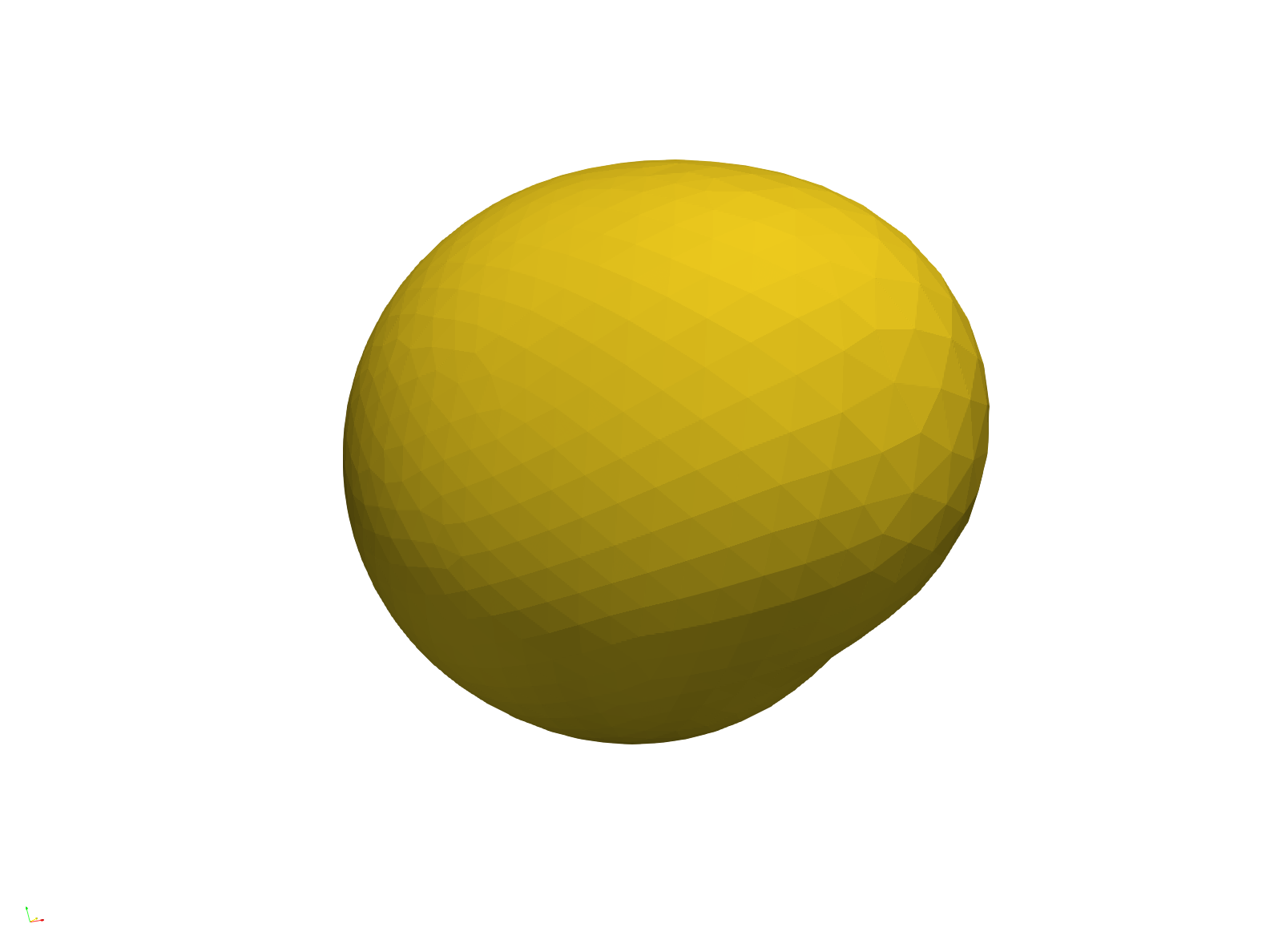}
		\caption{Final Result}
	\end{subfigure}
	
	\caption{Simulation of a sphere with a reduced fixed boundary area at the bottom.}
	\label{fig:sphere_small_bc}
\end{figure}

In the previous simulation (\cref{fig:my_five_plots}), the displacement was fixed on the entire lower hemisphere, which restricts the deformation modes significantly. Here, we relax this condition to simulate a structure supported by a small ``stem" or tissue connection.
The Dirichlet boundary condition $\mathbf{u} = \mathbf{0}$ is applied only to a small patch at the bottom of the sphere ($z < -0.4$), while the rest of the boundary is traction-free. The growth center is located at $(0.25, 0, 0)$. As shown in \cref{fig:sphere_small_bc}, the relaxed boundary allows the sphere to expand and tilt more freely, resulting in a potato-like deformation that was previously constrained.

\subsubsection{Ellipsoid}
\begin{figure}[H] 
	\centering
	\setlength{\abovecaptionskip}{3pt} 
	\setlength{\belowcaptionskip}{0pt}
	
	\begin{subfigure}[b]{0.32\textwidth}
		\centering
		\includegraphics[width=\textwidth]{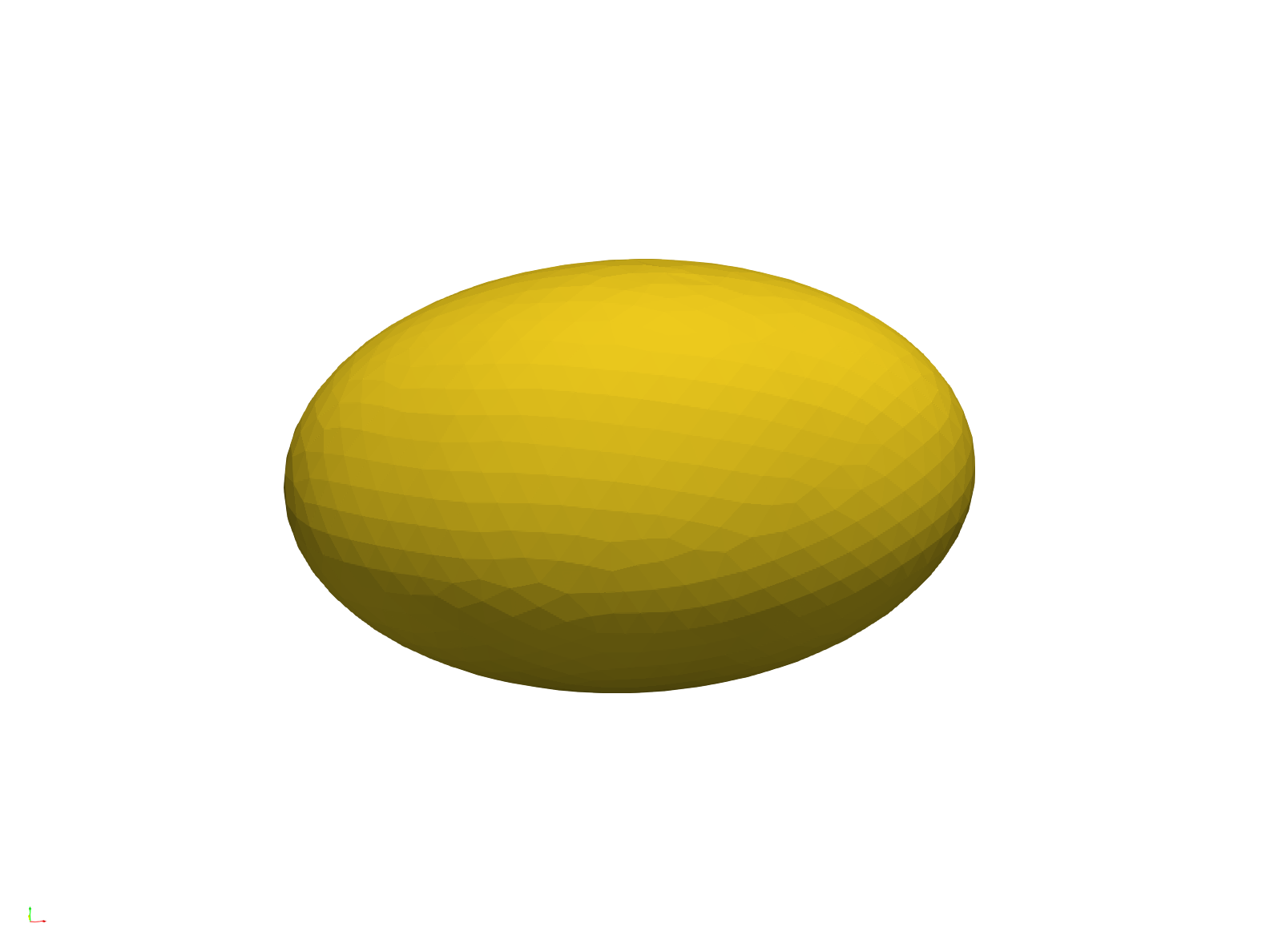}
		\caption{Initial Ellipsoid}
	\end{subfigure}
	\hfill 
	\begin{subfigure}[b]{0.32\textwidth}
		\centering
		\includegraphics[width=\textwidth]{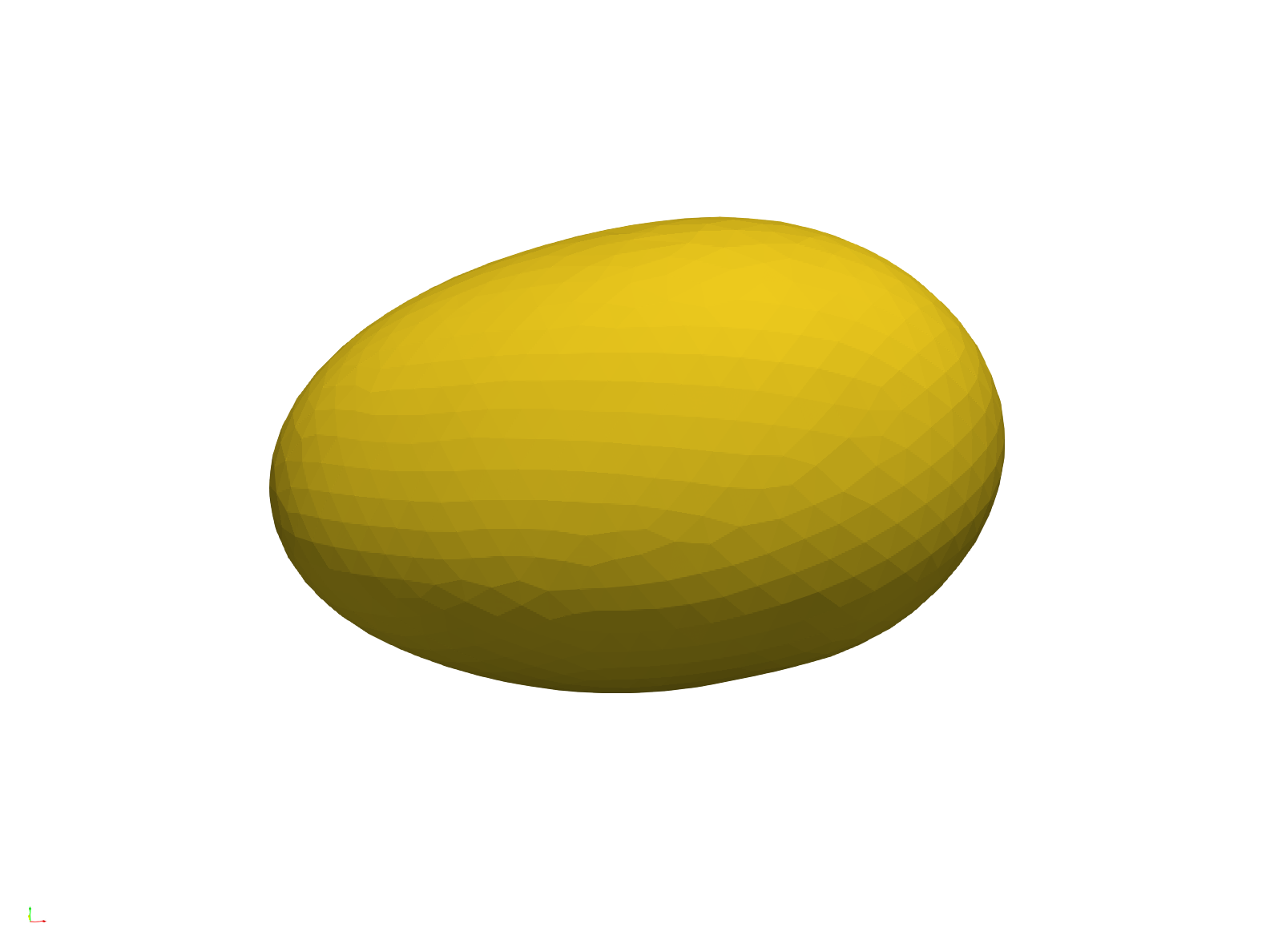}
		\caption{Iteration 2}
	\end{subfigure}
	\hfill 
	\begin{subfigure}[b]{0.32\textwidth}
		\centering
		\includegraphics[width=\textwidth]{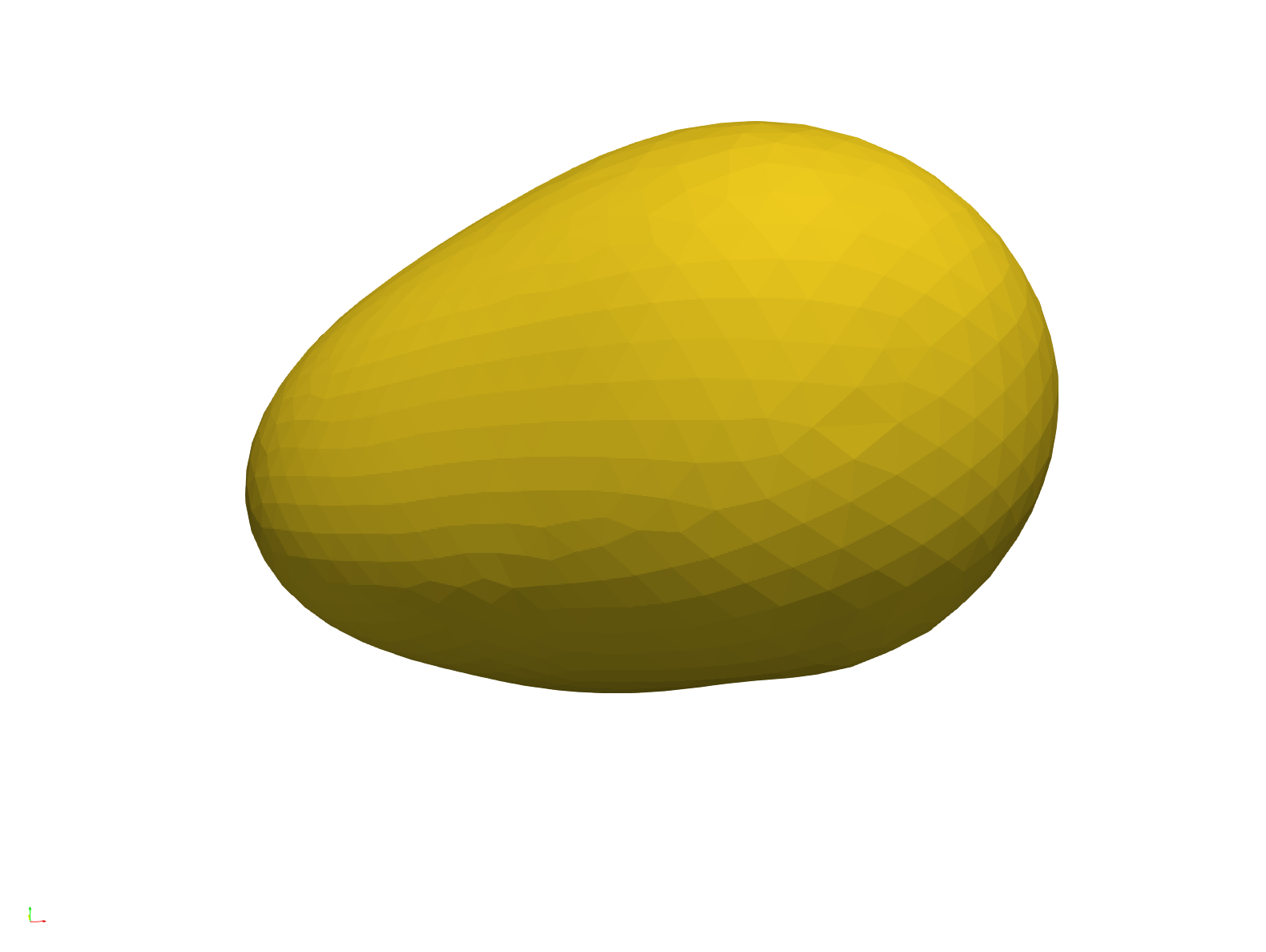}
		\caption{Final Result}
	\end{subfigure}
	
	\caption{Evolution of an ellipsoid.}
	\label{fig:ellipsoid}
\end{figure}

Next, we consider an ellipsoid generated by dilating a unit sphere with factors $(1.6, 1.0, 1.0)$. The growth center is placed offset along the long axis at $(0.3, 0, 0)$, and with a reduced fixed boundary area at the bottom. \cref{fig:ellipsoid} demonstrates that the growth induces an asymmetric bulging and bending of the ellipsoid. 

\subsubsection{Kidney Shape: Non-Convexity and Smoothness}
Also, we simulate a shape that looks like a "kidney" or "bean", representing a non-convex domain. A critical requirement for the LDDMM operator framework and high-order regularization is that the boundary must be sufficiently smooth. 

\begin{figure}[H] 
	\centering
	\setlength{\abovecaptionskip}{3pt} 
	\setlength{\belowcaptionskip}{0pt}
	
	\begin{subfigure}[b]{0.32\textwidth}
		\centering
		\includegraphics[width=\textwidth]{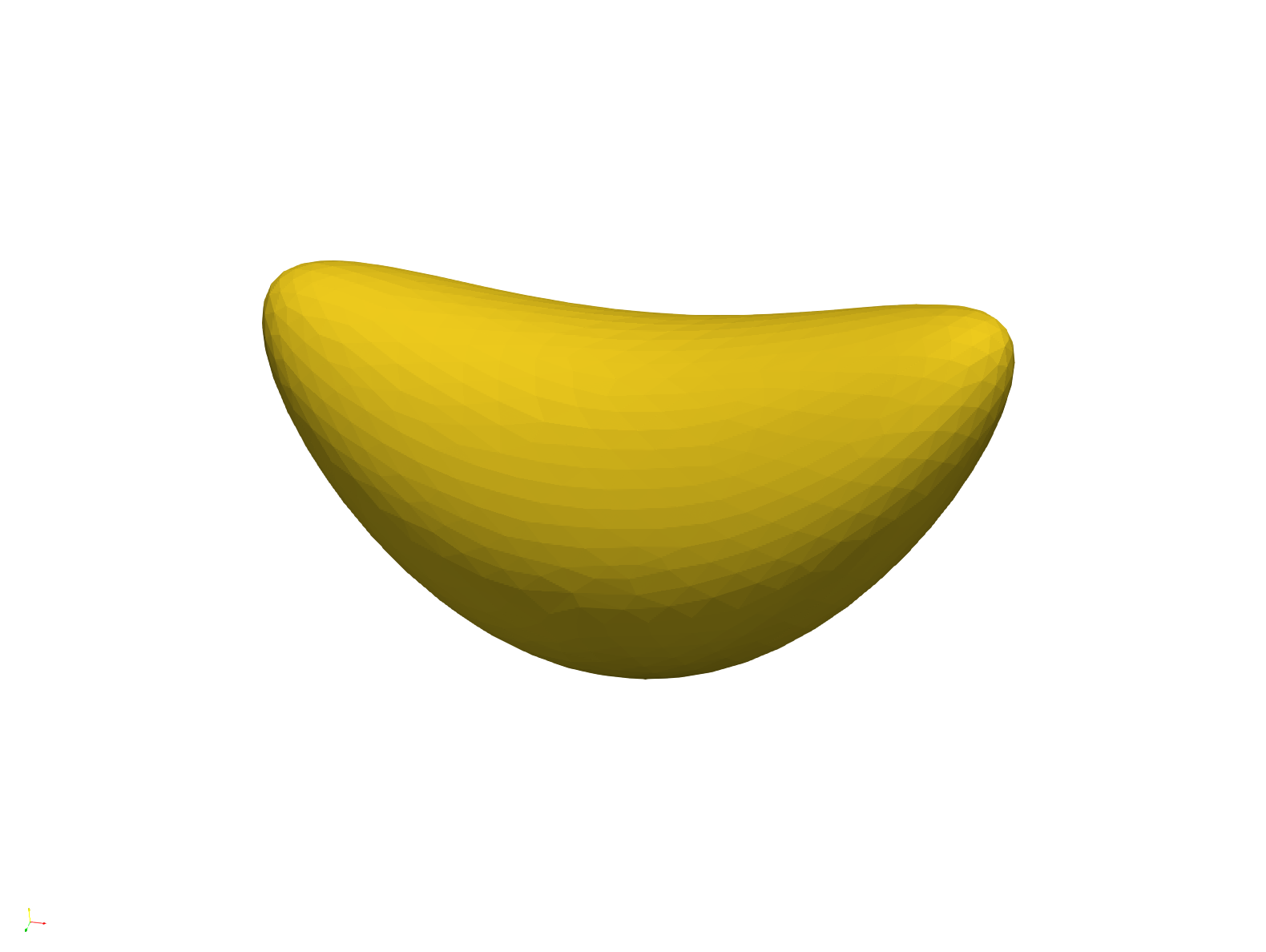}
		\caption{Initial Smooth Kidney}
	\end{subfigure}
	\hfill 
	\begin{subfigure}[b]{0.32\textwidth}
		\centering
		\includegraphics[width=\textwidth]{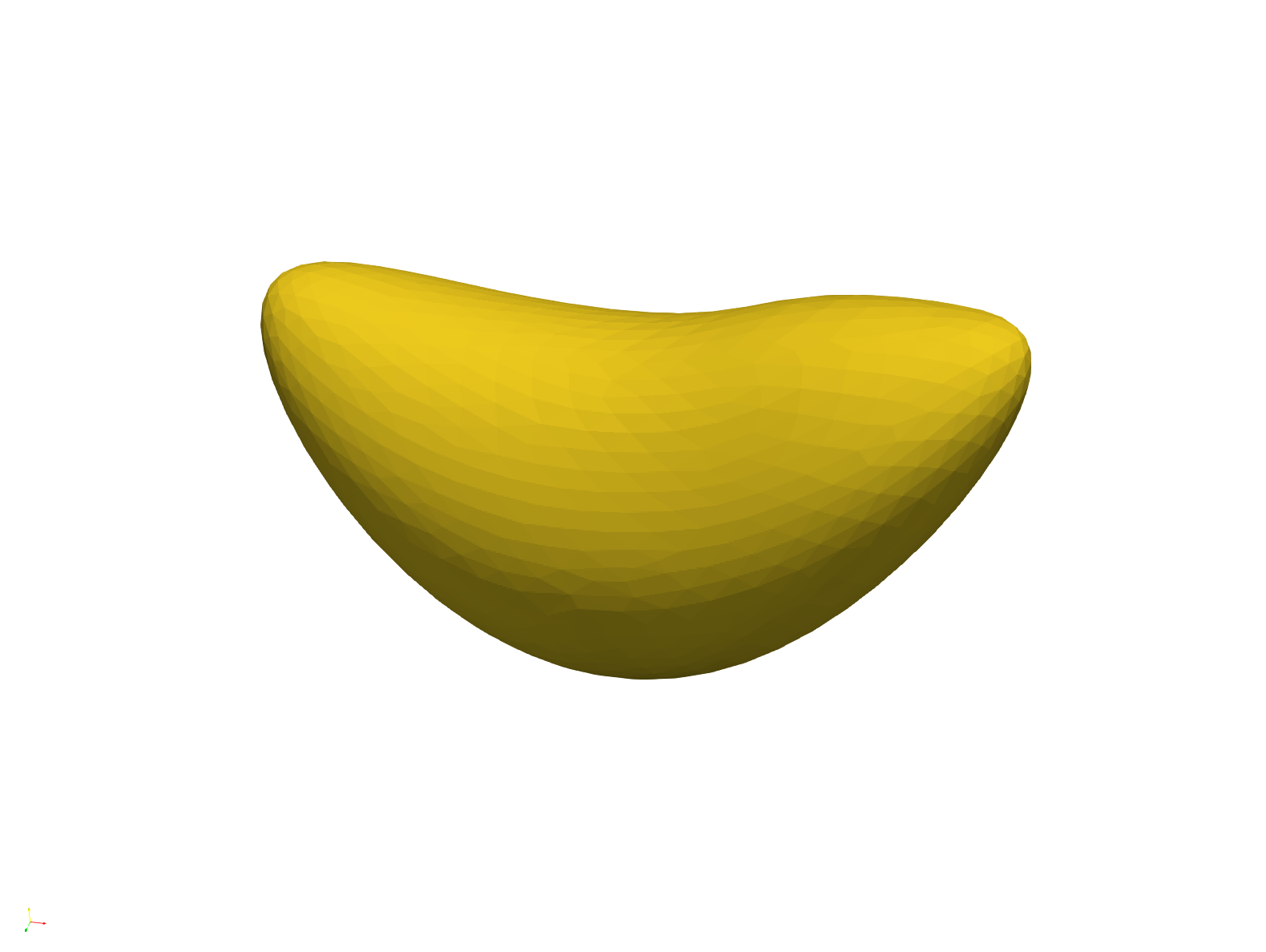} 
		\caption{Iteration 1}
	\end{subfigure}
	\hfill 
	\begin{subfigure}[b]{0.32\textwidth}
		\centering
		\includegraphics[width=\textwidth]{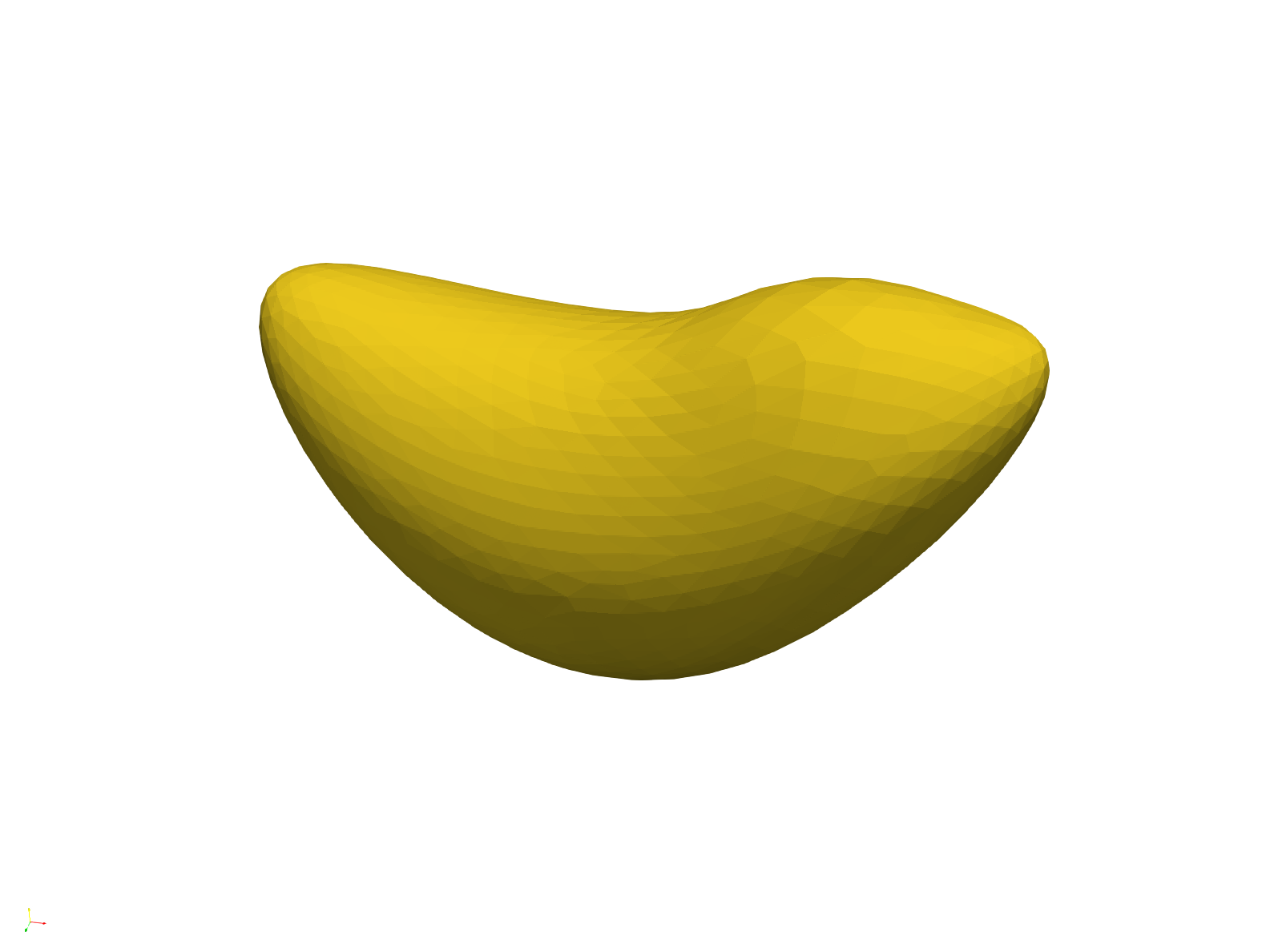}
		\caption{Iteration 2}
	\end{subfigure}
	\hfill 
	\begin{subfigure}[b]{0.32\textwidth}
		\centering
		\includegraphics[width=\textwidth]{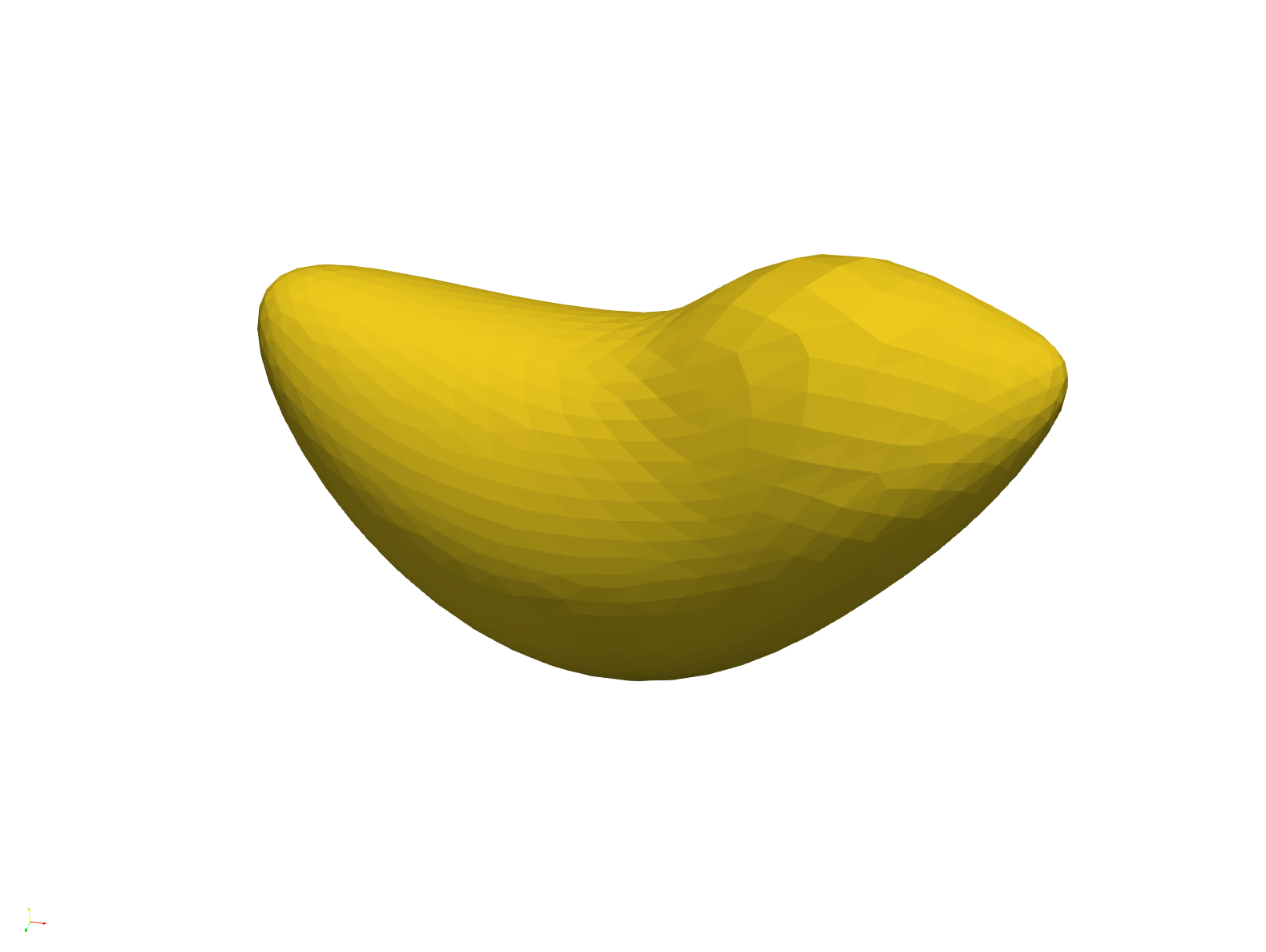} 
		\caption{Iteration 3}
	\end{subfigure}
	
	\caption{Simulation of a smooth, non-convex kidney shape generated by coordinate warping.}
	\label{fig:kidney}
\end{figure}

To ensure global $C^\infty$ smoothness, we generate this domain via a smooth coordinate transformation of the sphere defined by the mapping $\boldsymbol{\Phi}(X,Y,Z)$:
\begin{equation}
	\boldsymbol{\Phi}(X,Y,Z) = 
	\begin{pmatrix}
		x \\ y \\ z
	\end{pmatrix}
	= 
	\begin{pmatrix}
		1.2 X \\
		0.6 Y + \alpha_b (X^2 -\alpha_o) \\
		0.6 Z
	\end{pmatrix}
\end{equation}
where $\alpha_b=0.6$ controls the bending curvature, and parameter $\alpha_o = 0.5$. The growth center is placed inside the kidney-like shape. \cref{fig:kidney} shows the result. The solver successfully computes the deformation on this non-convex manifold. The regularization term prevents the mesh from collapsing in the concave region, demonstrating the method's applicability to complex biological organs.

\subsection{Negative Growth and Simulation}
Recall that the growth tensor $\tens{G}$ is defined as $\tens{G}(\vec{x}) = g(\vec{x})\mathbf{I} = A \exp\left(-\frac{|\vec{x} - \vec{x}_0|^2}{2\tau^2}\right) \mathbf{I}$ where $A>0$ typically represents volumetric expansion. However, the mathematical framework derived in Section 9 is general and naturally extends to the case of negative growth, or tissue atrophy, by allowing the amplitude coefficient to be negative:
\begin{equation}
	A < 0.
\end{equation}
Physically, a negative growth rate corresponds to a localized loss of mass or volume resorption, commonly observed in biological processes such as cell apoptosis, muscle atrophy, or the therapeutic shrinkage of tumors.
 
To demonstrate this, we apply the same higher-order LDDMM solver to the sphere and kidney shape but set the growth amplitude to $A = -0.5$ and $A = -0.3$ respectively. The regularization term is particularly crucial here to prevents mesh self-intersection or collapse at the center of the sink.

\begin{figure}[H] 
	\centering
	\setlength{\abovecaptionskip}{3pt} 
	\setlength{\belowcaptionskip}{0pt}
	
	\begin{subfigure}[b]{0.48\textwidth}
		\centering
		\includegraphics[width=\textwidth]{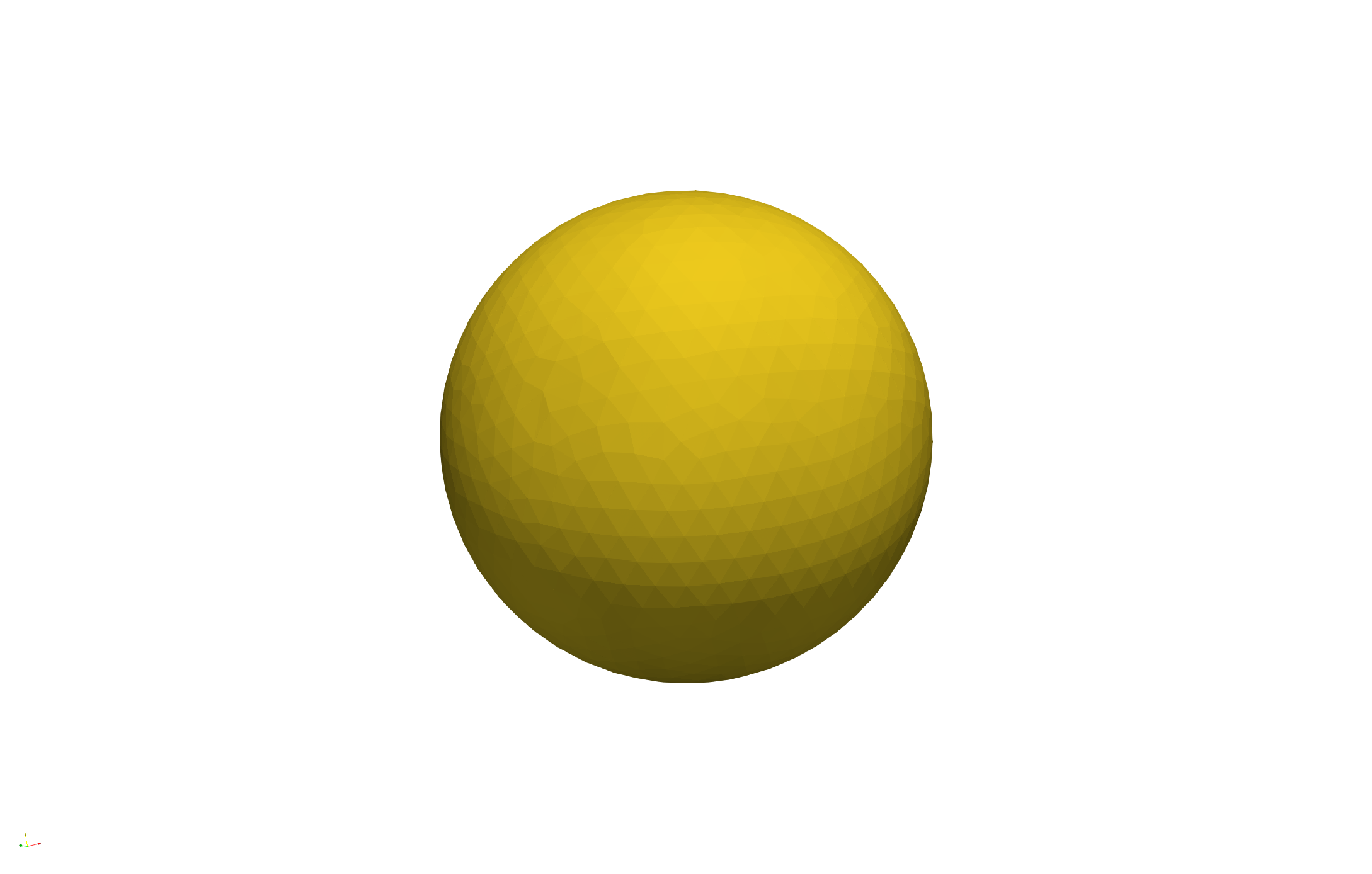}
		\caption{Initial Sphere}
	\end{subfigure}
	\hfill 
	\begin{subfigure}[b]{0.48\textwidth}
		\centering
		\includegraphics[width=\textwidth]{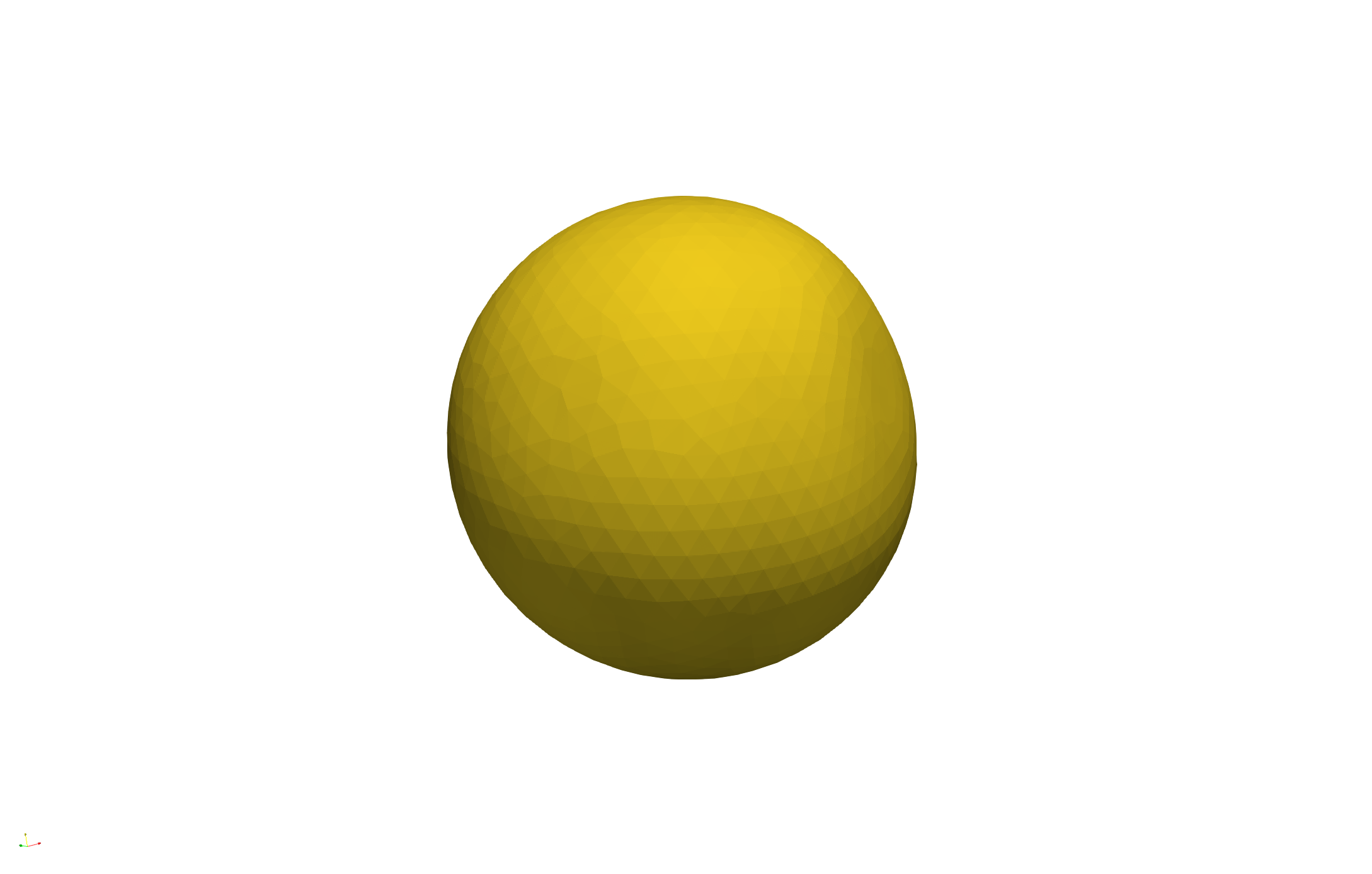}
		\caption{Iteration 2}
	\end{subfigure}

	\begin{subfigure}[b]{0.48\textwidth}
		\centering
		\includegraphics[width=\textwidth]{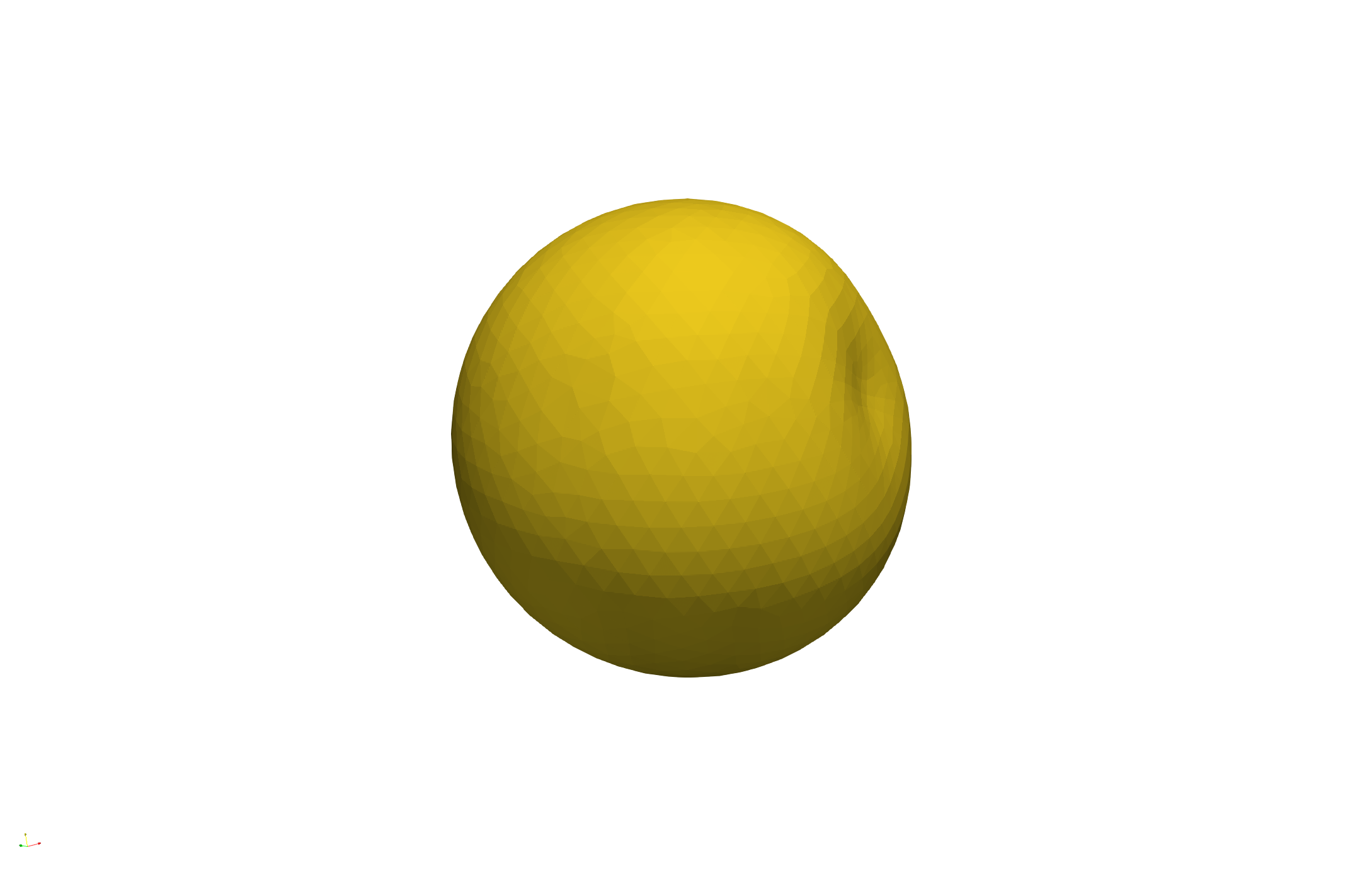}
		\caption{Final Result}
	\end{subfigure}
	
	\caption{Evolution of a sphere with negative growth.}
	\label{fig:ngsphere}
\end{figure}

\begin{figure}[H] 
	\centering
	\setlength{\abovecaptionskip}{3pt} 
	\setlength{\belowcaptionskip}{0pt}

	\begin{subfigure}[b]{0.48\textwidth}
		\centering
		\includegraphics[width=\textwidth]{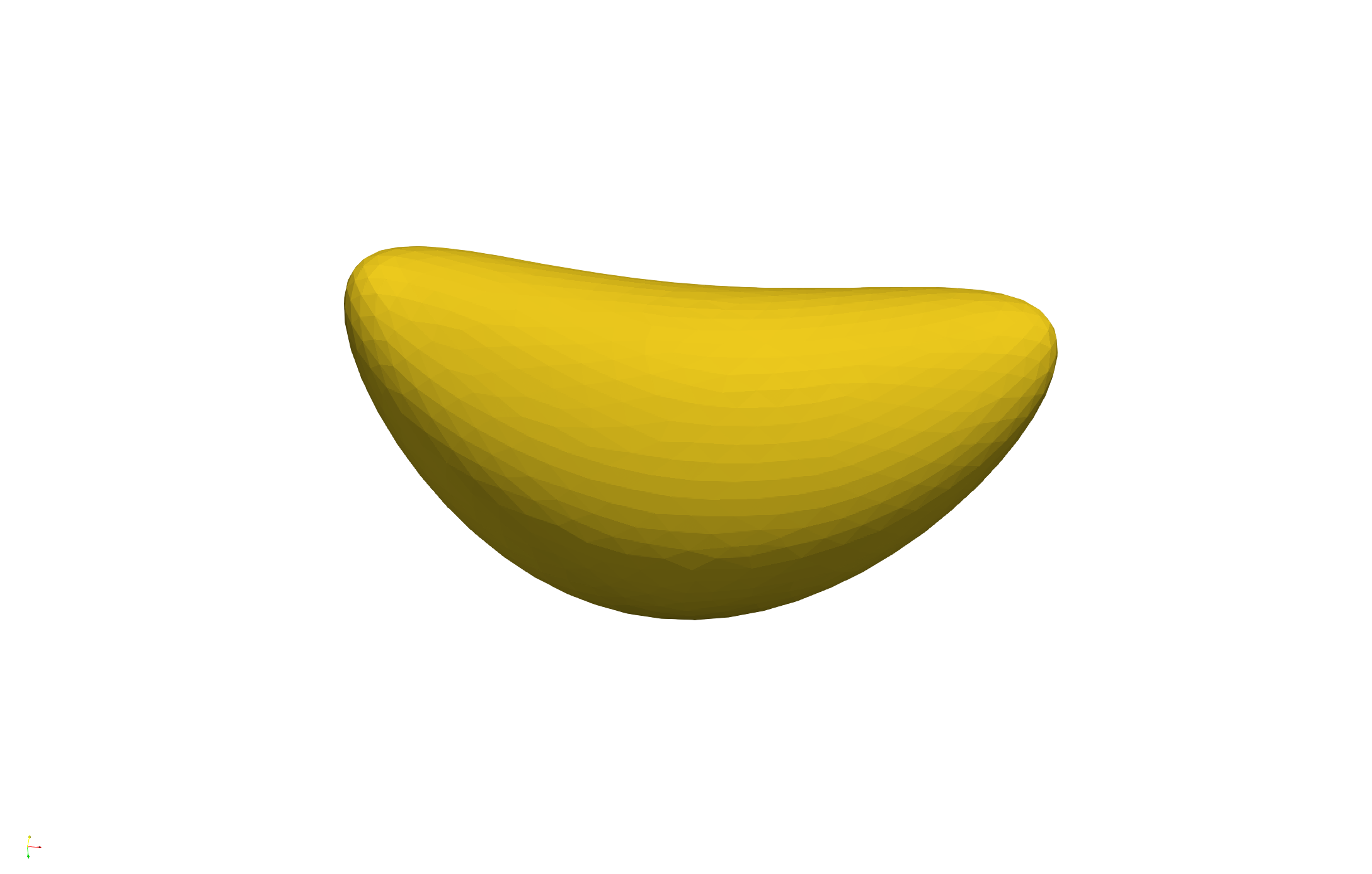}
		\caption{Initial Shape}
	\end{subfigure}
	\hfill 
	\begin{subfigure}[b]{0.48\textwidth}
		\centering
		\includegraphics[width=\textwidth]{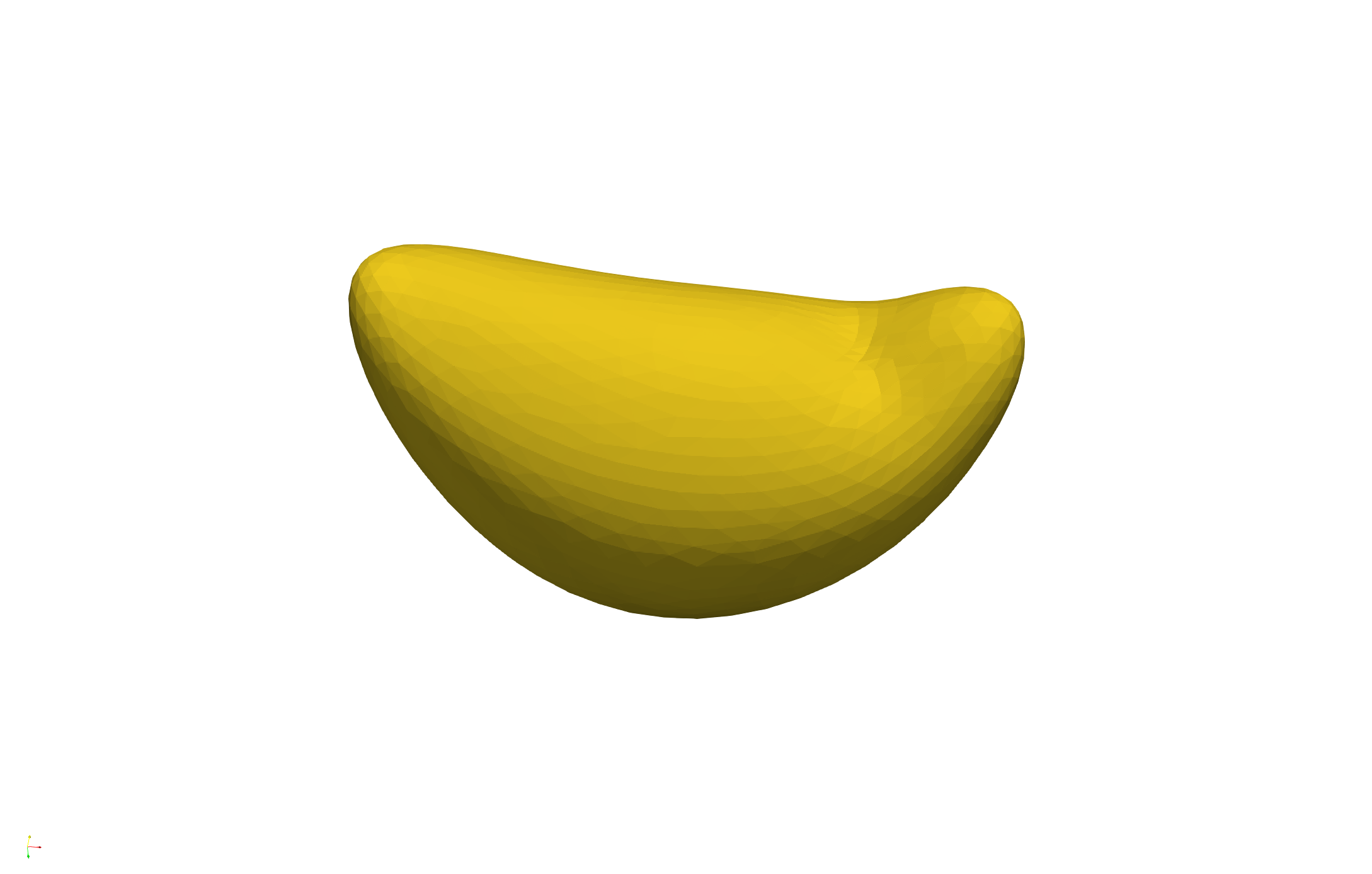}
		\caption{Iteration 2}
	\end{subfigure}

	\begin{subfigure}[b]{0.48\textwidth}
		\centering
		\includegraphics[width=\textwidth]{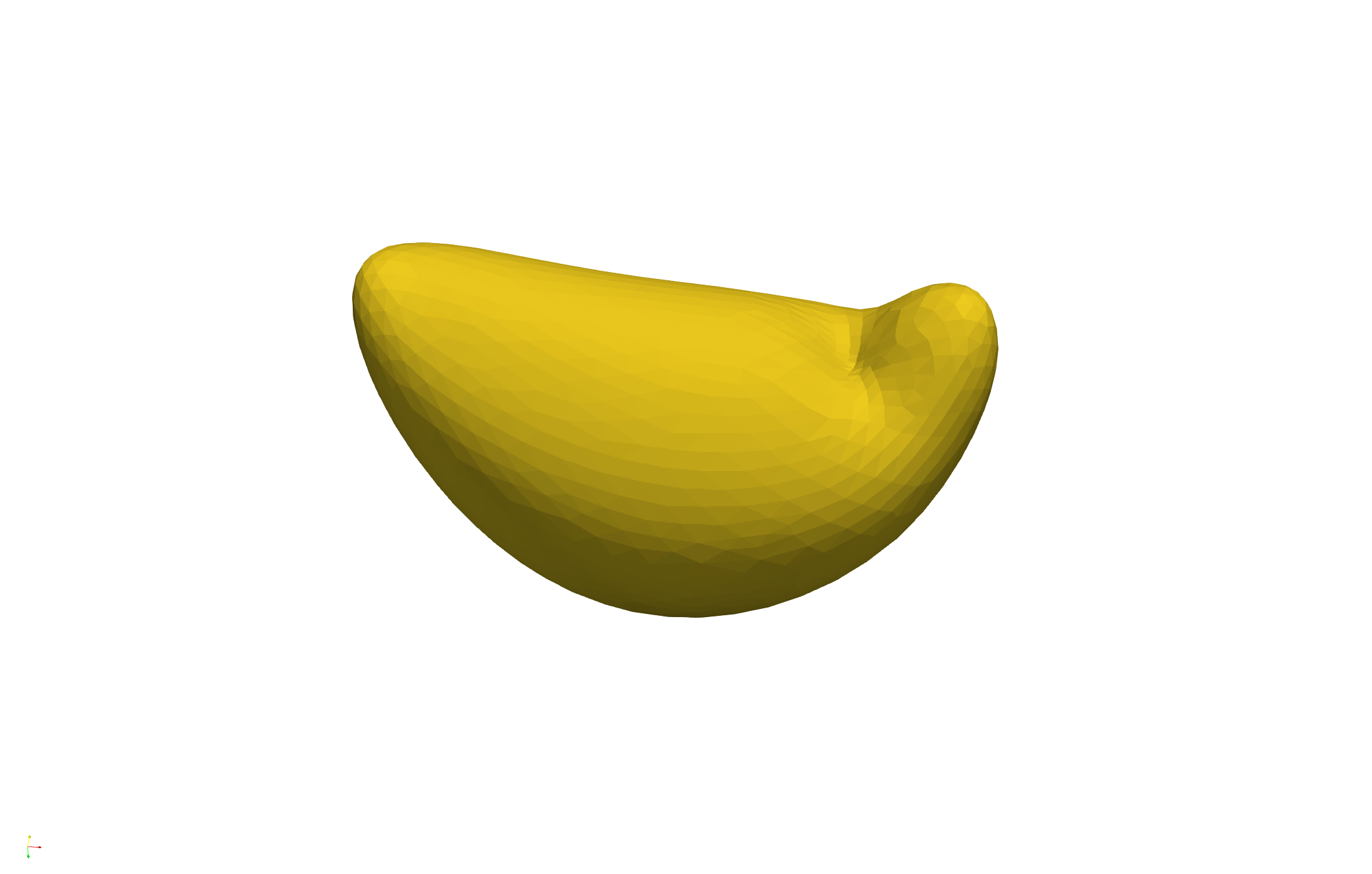}
		\caption{Final Result}
	\end{subfigure}
	
	\caption{Evolution of a kidney shape with negative growth.}
	\label{fig:ngkidney}
\end{figure}
\clearpage

\subsection{Summary of the Experiments}

\begin{table}[htbp]
    \centering
    \caption{Computational performance of the evolution with LDDMM operator across different geometries (Tested on MacBook M3 2024 CPU with the averaged reported execution times across multiple runs).}
    \label{tab:computational_performance}
    \begin{tabular}{>{\raggedright\arraybackslash}p{6cm} >{\centering\arraybackslash}p{2.5cm} >{\centering\arraybackslash}p{3cm}}
        \toprule
        \textbf{Configuration} & \textbf{Number of} & \textbf{Total Time} \\
                               & \textbf{Elements}  & \textbf{(4 steps) [s]} \\
        \midrule
        Sphere ($\kappa = 1$)   & 14,201 & 973 \\
        \addlinespace[0.5em]
        Sphere ($\kappa = 0.1$) & 14,201 & 775 \\
        \addlinespace[0.5em]
        Sphere (SFB, $\kappa = 1$) & 14,201 & 687 \\
        \addlinespace[0.5em]
        Negative $\tens{G}$ Sphere (SFB, $\kappa = 1$) & 14,201 & 864 \\
        \addlinespace[0.5em]
        Bean shape (SFB, $\kappa = 1$) & 14,201 & 729 \\
        \bottomrule
    \end{tabular}
    \vspace{1ex}
    \flushleft{\footnotesize \textit{Note: SFB denotes Small Fixed Boundary.}}
\end{table}
The results summarized in \cref{tab:computational_performance} demonstrate that the total execution time is related to several critical factors: the magnitude of the regularization parameter $\kappa$, the direction of growth, the geometric complexity of the domain, and the specific configuration of the boundary conditions.

In conclusion, despite the high-order nature of the mixed formulation, the computational cost remains manageable on a standard consumer-grade processor for a moderate number of iterations. For the static model, execution times are reduced to seconds, even when employing more elements than those utilized in the dynamic LDDMM experiments. 

\clearpage

\section{Discussion and Future Directions}

\subsection{The Challenge of Inverse Morphoelasticity}
Thus far, this work has focused on the \textit{forward problem}: given a prescribed growth tensor $\mathbf{G}$, we compute the resulting displacement field $\mathbf{u}$, then we get the estimated velocity field $\mathbf{v}(t)$. However, in many biological applications (e.g., tumor expansion, organ development), the scenario is reversed. We observe the morphological change which is the displacement $\mathbf{u}_{\text{obs}}$ and seek to infer the underlying biological growth mechanism - the tensor $\mathbf{G}$.

As discussed in the previous section, determining $\mathbf{G}$ from $\mathbf{u}$ is an ill-posed inverse problem due to the existence of non-trivial compatible growth modes that produce zero stress or displacement. Furthermore, in the context of our previously derived high-order regularization model, solving this inverse problem using traditional PDE-constrained optimization is computationally prohibitive. It requires iteratively deriving and solving complex adjoint equations for the eighth-order system. 

\subsection{Physics-Informed Neural Networks (PINNs) Approach}

To address these challenges, we propose a framework based on Physics-Informed Neural Networks \cite{Karniadakis2021}, which has demonstrated that it can be used to inverse problem of PDEs while learning from available data and offers a compelling alternative to traditional mesh-based methods \cite{RAISSI2019686}. Unlike the Finite Element Method that requires solving complex weak forms, PINNs approximate the solution using a deep neural network $\mathbf{u}_{\boldsymbol{\theta}}(\mathbf{x})$ parametrized by weights $\boldsymbol{\theta}$. 

A primary advantage of this approach in the context of the biharmonic equation is its ability to handle the eighth-order term $\Delta^4 \mathbf{u}$ directly via automatic differentiation. This eliminates the necessity for auxiliary variables common in mixed formulations, allowing the problem to remain a direct minimization over $\mathbf{u}$ alone. Furthermore, because PINNs are mesh-free, they bypass the requirement for specialized $C^1$-continuous Hermite elements or complex mixed function spaces. This meshless nature is particularly advantageous for modeling large-deformation morphoelasticity, as it naturally accommodates significant geometric changes without the mesh quality degradation that would otherwise necessitate expensive re-meshing or Arbitrary Lagrangian-Eulerian techniques.

Within this unified optimization framework, we can solve both forward and inverse problems simultaneously. We represent the unknown growth tensor $\mathbf{G}(\mathbf{x})$ either as a separate neural network or as a trainable parametric function (e.g., a Gaussian field with trainable amplitude and center). The training objective is to minimize a composite loss function:
\begin{equation}
	\mathcal{L}(\boldsymbol{\theta}, \mathbf{G}) = \mathcal{L}_{\text{data}} + \lambda \mathcal{L}_{\text{PDE}}
\end{equation}

By minimizing $\mathcal{L}(\boldsymbol{\theta}, \mathbf{G})$ in a certain design, the network is able to find an approximate solution that satisfies the biharmonic regularity and identifies the optimal growth tensor $\mathbf{G}$ that explains the observation.

\subsection{Operator Learning for Geodesic Shooting}
The current framework relies on solving the optimization problem for each new instance of growth or boundary condition. For the geodesic path problem in shape spaces, which requires solving the evolution equation iteratively over time, the computational cost can be high.

A promising future direction is the application of \textbf{Neural Operators}, such as Deep Operator Networks (DeepONet) \cite{Karniadakis2021}, which is proved to serve as universal approximators for nonlinear continuous operators \cite{Lu2021}. Instead of approximating a single solution, these models learn the solution operator mapping the functional space of growth tensors to the space of displacement fields:
\begin{equation}
	\mathcal{N}: \mathbf{G}(\mathbf{x}) \mapsto \mathbf{u}(\mathbf{x})
\end{equation}
Once trained offline, a Neural Operator can infer the morphoelastic deformation resulting from any arbitrary growth distribution efficiently. And it is also worth mentioning that in multiphysics applications like electro-convection, one can combine DeepONets with
physics encoded by PINNs, to accomplish
real-time accurate predictions with extrapolation \cite{CAI2021110296}. 

Therefore, utilizing DeepONets as an efficient surrogate provides a transformative solution to the computational bottlenecks of geodesic shooting. By significantly accelerating the iterative matching process, this approach may bridge the gap between rigorous shape analysis theory and the latency constraints of real-time biomedical applications.

\section*{Acknowledgments}

I would like to express my deepest gratitude to my advisor, Professor Laurent Younes, for his invaluable guidance, insightful mentorship, and continuous support throughout this independent research project. I am also immensely grateful to Professor Amitabh Basu and Professor Kobe Marshall-Stevens for their engaging discussions and constructive feedback, which greatly benefited the development of this research. 

Furthermore, I would like to extend my sincere thanks to all three professors for their generous support, unwavering encouragement, and invaluable advice during my Ph.D. application process.

\clearpage

\appendix
\crefalias{section}{appendix}

\section{Derivation and Settings for Lamé–Navier Equation} \label{appendix:Lamé–Navier}
\subsubsection*{1. Start from the Equilibrium Equation}
	The equation for static equilibrium (Cauchy's first law of motion with no acceleration) is:
	\begin{align*}
		\nabla \cdot \bm{\sigma} + \mathbf{f} &= \mathbf{0} \\
		% Index notation
		\text{or in index notation:} \quad \sigma_{ij,i} + f_j &= 0 
	\end{align*}
	where $\bm{\sigma}$ is the Cauchy stress tensor and $\mathbf{f}$ is the body force per unit volume.
	
	\subsubsection*{2. Substitute the Constitutive Relation}
	For a linear, isotropic, elastic material, Hooke's law relates stress to strain:
	\begin{equation*}
		\sigma_{ij} = \lambda \delta_{ij} \varepsilon_{kk} + 2\mu \varepsilon_{ij}
	\end{equation*}
	where $\lambda$ and $\mu$ are the Lamé constants, and $\delta_{ij}$ is the Kronecker delta. Substituting this into the equilibrium equation gives:
	\begin{align*}
		(\lambda \delta_{ij} \varepsilon_{kk} + 2\mu \varepsilon_{ij})_{,i} + f_j &= 0 \\
		% Using linearity of the derivative and assuming spatially constant Lamé parameters:
		\lambda (\delta_{ij} \epsilon_{kk})_{,i} + 2\mu (\varepsilon_{ij})_{,i} + f_j &= 0 \\
		% Since the Kronecker delta is constant:
		\lambda \delta_{ij} \varepsilon_{kk,i} + 2\mu \varepsilon_{ij,i} + f_j &= 0 
	\end{align*}
	By the sifting property of the Kronecker delta, the summation over $i$ in the first term retains only the $i=j$ component, simplifying $\delta_{ij} \varepsilon_{kk,i}$ to $\varepsilon_{kk,j}$. The equation becomes:
	\begin{equation}
		\lambda \varepsilon_{kk,j} + 2\mu \varepsilon_{ij,i} + f_j = 0 \label{eq:stress_in_strain_eng}
	\end{equation}
	
	\subsubsection*{3. Substitute the Kinematic Relation}
	Next, we express strain in terms of displacement $\mathbf{u}$ using the small-strain kinematic relation:
	\begin{equation*}
		\varepsilon_{ij} = \frac{1}{2} (u_{i,j} + u_{j,i})
	\end{equation*}
	The trace of the strain tensor (volumetric strain) is therefore:
	\begin{equation*}
		\varepsilon_{kk} = \frac{1}{2} (u_{k,k} + u_{k,k}) = u_{k,k}
	\end{equation*}
	Substituting these into Equation \cref{eq:stress_in_strain_eng}:
	\begin{align*}
		\lambda (u_{k,k})_{,j} + 2\mu \left( \frac{1}{2} (u_{i,j} + u_{j,i}) \right)_{,i} + f_j &= 0 \\
		\lambda u_{k,kj} + \mu (u_{i,j} + u_{j,i})_{,i} + f_j &= 0 \\
		\lambda u_{k,kj} + \mu (u_{i,ji} + u_{j,ii}) + f_j &= 0
	\end{align*}
	
	\subsubsection*{4. Rearrange to Obtain the Final Equation}
	We simplify the last expression by manipulating the dummy indices and assuming the displacement field is smooth enough to interchange the order of partial derivatives.
	In the term $\mu u_{i,ji}$, we rename the dummy index $i$ to $k$, yielding $\mu u_{k,jk}$. Then, we switch the order of differentiation: $u_{k,jk} = u_{k,kj}$. The equation becomes:
	\begin{align*}
		\lambda u_{k,kj} + \mu u_{k,kj} + \mu u_{j,ii} + f_j &= 0 \\
		(\lambda + \mu) u_{k,kj} + \mu u_{j,ii} + f_j &= 0
	\end{align*}
	This is the Navier-Cauchy equation in \textbf{index form}. To convert it to \textbf{vector form}, we identify the differential operators:
	\begin{itemize}
		\item $u_{k,k} = \nabla \cdot \mathbf{u}$ (Divergence of $\mathbf{u}$)
		\item $u_{k,kj} = (u_{k,k})_{,j}$ is the $j$-th component of $\nabla(\nabla \cdot \mathbf{u})$ (Gradient of the divergence)
		\item $u_{j,ii}$ is the $j$-th component of $\nabla^2 \mathbf{u}$ (Vector Laplacian of $\mathbf{u}$)
	\end{itemize}
	Combining these yields the final vector form of the Navier-Cauchy equation:
	\begin{equation}
		\boxed{
			(\lambda + \mu) \nabla(\nabla \cdot \mathbf{u}) + \mu \nabla^2 \mathbf{u} + \mathbf{f} = \mathbf{0}
		}
	\end{equation}
	This is a linear elliptic system of equations for the displacement $\mathbf{u}$. It can be shown that the system has a unique classical solution if the boundary conditions are well-posed. For instance, a portion of the boundary must be fixed to ensure the object cannot undergo rigid body motion, thus excluding the pure traction case. Also, this holds provided that the bulk and shear moduli are positive\cite{sadd2020elasticity},\cite{gould2018introduction}: \label{solution}
	\begin{align*}
		K &= \frac{E}{3(1 - 2\nu)} = \lambda + \frac{2}{3}\mu> 0 \\
		G &= \frac{E}{2(1 + \nu)} = \mu > 0
	\end{align*}
	which poses the following restrictions on the Poisson's ratio:
	\[
	-1 < \nu < 0.5
	\]
	Here we make everything satisfied as our settings and let $\mathbf{f} \in L^2(\Omega)\text{ and } \mathbf{T} \in L^2(\partial\Omega)$, where $\mathbf{T}$ is the traction already defined \cref{12}.

 \section{Derivation of the Variational Equations for the Linear Elastic Model by Stokes' Theorem}  

 \label{appendix: Variational_Linear}

 We start with the equilibrium equation. Pick an arbitrary vector test function $\vec{v} \in \hat{\mathcal{V}}$ (where $\hat{\mathcal{V}}$ is a suitable vector-valued function space where test functions are zero on the Dirichlet boundary), and integrate over the domain $\Omega$:
	\begin{equation}
		-\int_{\Omega} (\nabla \cdot \tens{\sigma(\varepsilon(\vec{u})}) \cdot \vec{v} \, dV = \int_{\Omega} \vec{f} \cdot \vec{v} \, dV
	\end{equation}
	We now apply integration by parts to the left-hand side. Using the chain rule, one can derive that $\nabla \cdot (\tens{\sigma} \vec{v}) = (\nabla \cdot \tens{\sigma}) \cdot \vec{v} + \tens{\sigma} : \nabla\vec{v}$, where the double dot is denoted as $\bm{A} : \bm{B}=tr(\bm{A^TB})$. We can rewrite the integral:
	\begin{equation}
		-\int_{\Omega} (\nabla \cdot \tens{\sigma}) \cdot \vec{v} \, dV = \int_{\Omega} \tens{\sigma} : \nabla\vec{v} \, dV - \int_{\Omega} \nabla \cdot (\tens{\sigma} \vec{v}) \, dV
	\end{equation}
	Applying the Stokes' theorem (Here, Divergence theorem) to the second term on the right-hand side gives:
	\begin{equation}
		\int_{\Omega} \nabla \cdot (\tens{\sigma} \vec{v}) \, dV = \oint_{\partial\Omega} (\tens{\sigma} \vec{v}) \cdot \vec{n} \, dS = \oint_{\partial\Omega} (\tens{\sigma}\vec{n}) \cdot \vec{v} \, dS
	\end{equation}
	where $\vec{n}$ is the outward unit normal to the boundary $\partial\Omega$. And $dS$ here is the area element. The second equal sign holds because here $\tens{\sigma}$ is symmetric. The term $\tens{\sigma}\vec{n} = \tens{\sigma}^T\vec{n} $ is defined as the traction vector $\vec{T}$, which represents the surface forces \cite{gould2018introduction}. Substituting this back gives:
	\begin{equation}
		\int_{\Omega} \tens{\sigma} : \nabla\vec{v} \, dV - \oint_{\partial\Omega} \vec{T} \cdot \vec{v} \, dS = \int_{\Omega} \vec{f} \cdot \vec{v} \, dV
		\label{12}
	\end{equation}
	Rearranging the terms, we arrive at the final variational formulation: find $\vec{u}$ such that for all test functions $\vec{v}$:
	\begin{equation}
		\int_{\Omega} \tens{\sigma(\varepsilon(\vec{u})}) : \nabla\vec{v} \, dV = \int_{\Omega} \vec{f} \cdot \vec{v} \, dV + \int_{\Gamma_N} \vec{T} \cdot \vec{v} \, dS
	\end{equation}
	The boundary integral is only evaluated over the Neumann boundary $\Gamma_N$, since $\vec{v} = \vec{0}$ on the Dirichlet boundary $\Gamma_D$. Where $\Gamma_D \cup \Gamma_N = \Gamma = \partial \Omega$.

\section{Derivation of the Variational Equations for the Morphoelastic Model}
\label{appendix:Variational_Morpho}

We evaluate the derivative term by term:

\begin{itemize}
	\item \textbf{Internal Strain Energy Term}\\
	We apply the chain rule. Let $\tens{\varepsilon}_g(\delta) = \tens{\varepsilon}(\vec{u} + \delta \vec{v}) - \tens{G} = \tens{\varepsilon}_g(\vec{u}) + \delta\tens{\varepsilon}(\vec{v})$.
	
	\begin{equation}
		\frac{d}{d\delta} \int_{\Omega} \Psi(\tens{\varepsilon}_g(\delta)) \,dV = \int_{\Omega} \frac{\partial \Psi}{\partial \tens{\varepsilon}_g} : \frac{d\tens{\varepsilon}_g}{d\delta} \,dV = \int_{\Omega} \tens{\sigma}(\tens{\varepsilon}_g) : \tens{\varepsilon}(\vec{v}) \,dV
		\label{27}
	\end{equation}
	
	Here we provide a detailed explanation for the strain energy term \cref{27}.

	The first step involves applying the chain rule to the expression $\frac{d}{d\delta} \Psi(\tens{\varepsilon}_g(\delta))$. This is a direct extension of the chain rule from multi-variable calculus.
	    For a scalar function of a vector $\vec{v}(t)$, the chain rule is $\frac{d}{dt}f(\vec{v}(t)) = \nabla f \cdot \frac{d\vec{v}}{dt}$. The dot product takes the place of multiplication. In \textbf{Tensor Calculus}, we extend this analogy. The scalar function is $\Psi$, the argument is the tensor $\tens{\varepsilon}_g(\delta)$, and the dot product is replaced by the double-dot product.
	Thus, the chain rule for our functional is:
	\begin{equation}
		\frac{d}{d\delta} \Psi(\tens{\varepsilon}_g(\delta)) = \frac{\partial \Psi}{\partial \tens{\varepsilon}_g} : \frac{d\tens{\varepsilon}_g}{d\delta}
	\end{equation}
	since $\tens{\varepsilon}_g(\delta) = \tens{\varepsilon}(\vec{u}) + \delta\tens{\varepsilon}(\vec{v}) - \tens{G}$, with respect to the scalar $\delta$. This derivative is straightforward:
	\begin{equation}
		\frac{d\tens{\varepsilon}_g}{d\delta} = \frac{d}{d\delta} \left( \tens{\varepsilon}(\vec{u}) + \delta\tens{\varepsilon}(\vec{v}) - \tens{G} \right) = \tens{\varepsilon}(\vec{v})
		\label{eq:inner_deriv}
	\end{equation}
	
	For tensor calculus,  $\frac{\partial \Psi}{\partial \tens{\varepsilon}_g}$ is a tensor whose components are $(\frac{\partial \Psi}{\partial \tens{\varepsilon}_g})_{ij} = \frac{\partial \Psi}{\partial \varepsilon_{g,ij}}$ \cite{goriely2017mathematics}.
	
	The strain energy density is:
	\begin{equation*}
		\Psi(\tens{\varepsilon}_g) = \frac{\lambda}{2}(\text{tr}(\tens{\varepsilon}_g))^2 + \mu\text{tr}(\tens{\varepsilon}_g^2)
	\end{equation*}
	
	One can verify that the resulting components are precisely the components of the stress tensor $\tens{\sigma}_g$. Therefore, we have proven that:
	\begin{equation*}
		\frac{\partial \Psi}{\partial \tens{\varepsilon}_g} = \lambda \text{tr}(\tens{\varepsilon}_g)\mathbf{I} + 2\mu\tens{\varepsilon}_g = \tens{\sigma}_g
	\end{equation*}
	
	\item \textbf{Body Force Potential Energy Term}
	\begin{equation*}
		\frac{d}{d\delta} \left(-\int_{\Omega} \vec{f} \cdot (\vec{u} + \delta \vec{v}) \,dV \right) = -\int_{\Omega} \vec{f} \cdot \vec{v} \,dV
	\end{equation*}
	
	\item \textbf{Surface Traction Potential Energy Term}
	\begin{equation*}
		\frac{d}{d\delta} \left(-\int_{\Gamma_N} \vec{T} \cdot (\vec{u} + \delta \vec{v}) \,dS \right) = -\int_{\Gamma_N} \vec{T} \cdot \vec{v} \,dS
	\end{equation*}
\end{itemize}

Combining these terms, the condition $\delta E = 0$ becomes:
\begin{equation*}
	\int_{\Omega} \tens{\sigma}(\tens{\varepsilon}_g) : \tens{\varepsilon}(\vec{v}) \,dV - \int_{\Omega} \vec{f} \cdot \vec{v} \,dV - \int_{\Gamma_N} \vec{T} \cdot \vec{v} \,dS = 0
\end{equation*}
Substituting $\tens{\varepsilon}_g = \tens{\varepsilon}(\vec{u}) - \tens{G}$ and rearranging, we get the final variational form:
\begin{equation}
	\int_{\Omega} \tens{\sigma}(\tens{\varepsilon}(\vec{u})) : \tens{\varepsilon}(\vec{v}) \,dV = \int_{\Omega} \vec{f} \cdot \vec{v} \,dV + \int_{\Gamma_N} \vec{T} \cdot \vec{v} \,dS + \int_{\Omega} \tens{\sigma}(\tens{G}) : \tens{\varepsilon}(\vec{v}) \,dV
\end{equation}

The displacement gradient $\nabla\mathbf{v}$ (a second-order tensor) describes the total local deformation of a material point. This deformation can be additively decomposed into two parts:
\begin{equation*}
	\nabla\mathbf{v} = \underbrace{\frac{1}{2} \left( \nabla\mathbf{v} + (\nabla\mathbf{v})^T \right)}_{\boldsymbol{\varepsilon}(\mathbf{v}) \text{ (Strain Tensor)}} +\underbrace{\frac{1}{2} \left( \nabla\mathbf{v} - (\nabla\mathbf{v})^T \right)}_{\boldsymbol{\omega}(\mathbf{v}) \text{ (Spin/Rotation Tensor)}}
\end{equation*}

Thus, we uses the term $\sigma(\mathbf{u}) : \nabla\mathbf{v}$ instead of $\sigma(\mathbf{u}) : \boldsymbol{\varepsilon}(\mathbf{v})$. The equivalence arises from the fundamental property that the stress tensor $\boldsymbol{\sigma}$ is symmetric ($\boldsymbol{\sigma} = \boldsymbol{\sigma}^T$). The double-dot product of a symmetric tensor with any skew-symmetric tensor is identically zero.
\begin{align*}
	\int_{\Omega} \boldsymbol{\sigma} : \nabla\mathbf{v} \, dV &= \int_{\Omega} \boldsymbol{\sigma} : (\boldsymbol{\varepsilon}(\mathbf{v}) + \boldsymbol{\omega}(\mathbf{v})) \, dV \\
	&= \int_{\Omega} \boldsymbol{\sigma} : \boldsymbol{\varepsilon}(\mathbf{v}) \, dV + \int_{\Omega} \underbrace{\boldsymbol{\sigma} : \boldsymbol{\omega}(\mathbf{v})}_{=0} \, dV \\
	&= \int_{\Omega} \boldsymbol{\sigma} : \boldsymbol{\varepsilon}(\mathbf{v}) \, dV
\end{align*}

\section{Proofs of Norm and Metric Properties}
\label{appendix:proofs}

\subsection*{Proof that $\|\mathbf{G}\|_{[\Omega]}$ Defines a Norm}

We demonstrate that $\|\mathbf{G}\|_{[\Omega]}$, derived from the minimization problem \cref{eq:growth_norm} satisfies the three defining properties of a norm on the vector space of growth rate tensors.

\paragraph{1. Positive Definiteness.} ($\|\mathbf{G}\|_{[\Omega]} \ge 0$, and $\|\mathbf{G}\|_{[\Omega]} = 0 \iff \mathbf{G} = \mathbf{0}$)

The term $\kappa \|\mathbf{v}\|^2_V$ is non-negative since $\kappa > 0$ and $\|\cdot\|_V$ is a norm. The elastic energy density $B(\cdot)$ is a positive semi-definite quadratic form, can be view as a metric tensor, so the integral is also non-negative. The minimum of non-negative values is non-negative, thus $\|\mathbf{G}\|_{[\Omega]} \ge 0$.

If $\mathbf{G} = \mathbf{0}$, the minimum is achieved at $\mathbf{v} = \mathbf{0}$, yielding $\|\mathbf{0}\|_{[\Omega]} = 0$.
Conversely, assume that $\|\mathbf{G}\|^2_{[\Omega]} = 0$. Because $v\mapsto \kappa \|\vec{v}\|_V^2 + \int_{\Omega} \Psi(\tens{\varepsilon}(\vec{v}) - \tens{G}) \, dV$ is a coercive strictly convex functional on $V$, it has a unique minimizer $\mathbf v_G$. Then for this minimizer both terms in the sum must be zero and $\kappa \|\mathbf{v}_{\mathbf{G}}\|^2_V = 0$ implies $\mathbf{v}_{\mathbf{G}} = \mathbf{0}$. The second term then becomes $\int_{\Omega} \Psi(-\mathbf{G}) dV = 0$. Since $\Psi$ is positive definite on the space of symmetric tensors, this requires $\mathbf{G}(\vec{x}) = \mathbf{0}$ for almost every $\vec{x} \in \Omega$.

\paragraph{2. Absolute Homogeneity.} ($\|\alpha\mathbf{G}\|_{[\Omega]} = |\alpha|\|\mathbf{G}\|_{[\Omega]}$ for any scalar $\alpha$)

Consider $\|\alpha\mathbf{G}\|^2_{[\Omega]} = \min_{\mathbf{v} \in V} \left( \kappa \|\mathbf{v}\|^2_V + \int_{\Omega} \Psi(\dot{\boldsymbol{\varepsilon}}(\mathbf{v}) - \alpha\mathbf{G}) \, dV \right)$. We perform a change of variables in the minimization, letting $\mathbf{v} = \alpha\mathbf{w}$.
\begin{align*}
	\|\alpha\mathbf{G}\|^2_{[\Omega]} &= \min_{\mathbf{w} \in V} \left( \kappa \|\alpha\mathbf{w}\|^2_V + \int_{\Omega}\Psi(\dot{\boldsymbol{\varepsilon}}(\alpha\mathbf{w}) - \alpha\mathbf{G}) \, dV \right) \\
	&= \min_{\mathbf{w} \in V} \left( \kappa \alpha^2 \|\mathbf{w}\|^2_V + \int_{\Omega} \alpha^2 \Psi(\dot{\boldsymbol{\varepsilon}}(\mathbf{w}) - \mathbf{G}) \, dV \right) \\
	&= \alpha^2 \min_{\mathbf{w} \in V} \left( \kappa \|\mathbf{w}\|^2_V + \int_{\Omega} \Psi(\dot{\boldsymbol{\varepsilon}}(\mathbf{w}) - \mathbf{G}) \, dV \right) = \alpha^2 \|\mathbf{G}\|^2_{[\Omega]}.
\end{align*}
Taking the square root yields the desired property.

\paragraph{3. Triangle Inequality.}
\[
    \|\mathbf{G}_1+\mathbf{G}_2\|_{[\Omega]}
    \le
    \|\mathbf{G}_1\|_{[\Omega]}
    +
    \|\mathbf{G}_2\|_{[\Omega]}.
\]

Let $\mathbf v_1$ and $\mathbf v_2$ be the unique minimizers associated with
$\mathbf G_1$ and $\mathbf G_2$, respectively. That is,
\[
    \|\mathbf G_i\|_{[\Omega]}^2
    =
    \kappa\|\mathbf v_i\|_V^2
    +
    \int_\Omega
    \Psi\bigl(\boldsymbol{\varepsilon}(\mathbf v_i)-\mathbf G_i\bigr)
    \,dV,
    \qquad i=1,2.
\]
Since $V$ is a linear velocity space, the function
$\mathbf v_1+\mathbf v_2$ is an admissible test velocity for the growth tensor
$\mathbf G_1+\mathbf G_2$. Therefore, by the definition of the minimum,
\begin{align*}
    \|\mathbf G_1+\mathbf G_2\|_{[\Omega]}^2
    &\le
    \kappa\|\mathbf v_1+\mathbf v_2\|_V^2
    +
    \int_\Omega
    \Psi\Bigl(
    \boldsymbol{\varepsilon}(\mathbf v_1+\mathbf v_2)
    -
    (\mathbf G_1+\mathbf G_2)
    \Bigr)
    \,dV \\
    &=
    \kappa\|\mathbf v_1+\mathbf v_2\|_V^2
    +
    \int_\Omega
    \Psi\Bigl(
    \bigl(\boldsymbol{\varepsilon}(\mathbf v_1)-\mathbf G_1\bigr)
    +
    \bigl(\boldsymbol{\varepsilon}(\mathbf v_2)-\mathbf G_2\bigr)
    \Bigr)
    \,dV,
\end{align*}
where we used the linearity of the strain operator
$\boldsymbol{\varepsilon}(\cdot)$.

We now introduce a product-space norm. Let
\[
    \mathrm{Sym}(3)
    :=
    \left\{
    A\in\mathbb R^{3\times 3}
    \mid A^T=A
    \right\}
\]
be the space of symmetric $3\times 3$ matrices, and define
\[
    L^2(\Omega;\mathrm{Sym}(3))
    =
    \left\{
    \mathbf S:\Omega\to\mathrm{Sym}(3)
    \ \middle|\
    \int_\Omega |\mathbf S(x)|^2\,dV<\infty
    \right\}.
\]
Here $|\mathbf S|^2=\mathbf S:\mathbf S$ is the Frobenius norm.

Define
\[
    \mathcal H
    =
    V\times L^2(\Omega;\mathrm{Sym}(3)).
\]
For $(\mathbf v,\mathbf S)\in\mathcal H$, set
\[
    \|(\mathbf v,\mathbf S)\|_{\mathcal H}
    =
    \left(
    \kappa\|\mathbf v\|_V^2
    +
    \int_\Omega \Psi(\mathbf S)\,dV
    \right)^{1/2}.
\]
Since $\kappa>0$ and $\Psi$ is a positive-definite quadratic form on
$\mathrm{Sym}(3)$, this expression defines a norm on $\mathcal H$. Equivalently,
it is the weighted product norm obtained from the $V$-norm and the elastic-energy
norm on $L^2(\Omega;\mathrm{Sym}(3))$.

Now define
\[
    h_1
    =
    \left(
    \mathbf v_1,
    \boldsymbol{\varepsilon}(\mathbf v_1)-\mathbf G_1
    \right),
    \qquad
    h_2
    =
    \left(
    \mathbf v_2,
    \boldsymbol{\varepsilon}(\mathbf v_2)-\mathbf G_2
    \right).
\]
Then the previous upper bound can be written as
\[
    \|\mathbf G_1+\mathbf G_2\|_{[\Omega]}
    \le
    \|h_1+h_2\|_{\mathcal H}.
\]
Since $\|\cdot\|_{\mathcal H}$ is a norm, it satisfies the triangle inequality. Therefore,
\begin{align*}
    \|\mathbf G_1+\mathbf G_2\|_{[\Omega]}
    &\le
    \|h_1+h_2\|_{\mathcal H} \\
    &\le
    \|h_1\|_{\mathcal H}
    +
    \|h_2\|_{\mathcal H} \\
    &=
    \left(
    \kappa\|\mathbf v_1\|_V^2
    +
    \int_\Omega
    \Psi\bigl(\boldsymbol{\varepsilon}(\mathbf v_1)-\mathbf G_1\bigr)
    \,dV
    \right)^{1/2} \\
    &\quad+
    \left(
    \kappa\|\mathbf v_2\|_V^2
    +
    \int_\Omega
    \Psi\bigl(\boldsymbol{\varepsilon}(\mathbf v_2)-\mathbf G_2\bigr)
    \,dV
    \right)^{1/2} \\
    &=
    \|\mathbf G_1\|_{[\Omega]}
    +
    \|\mathbf G_2\|_{[\Omega]}.
\end{align*}
Hence the triangle inequality holds.

All three properties are satisfied, confirming that $\|\mathbf{G}\|_{[\Omega]}$ is a norm.

\subsection*{Proof that the Functional Defines a Distance Metric}

We prove that the square root of the action function \cref{eq:action_functional} we defined satisfies the axioms of a metric. Thus, represents a geodesic from $\Omega_0$ to $\Omega_1$. 

First, we explain the existence of the optimal path to ensure the feasibility of the following proof. As shown in \cite{charon2023shape}, this problem is mathematically equivalent to simultaneously finding the optimal control $\tens{G}(t)$ and the optimal state trajectory (driven by the velocity field $\vec{v}(t)$) that minimize the combined action functional:
\begin{equation} \label{eq:morpho_action}
\min_{\vec{v}(\cdot), \, \tens{G}(\cdot) \in \mathcal{G}(\Omega(t))} \int_0^1 \left( \kappa \|\vec{v}(t)\|_V^2 + \int_{\Omega(t)} \Psi(\tens{\varepsilon}(\vec{v}(t)) - \tens{G}(t)) \, dV \right) dt
\end{equation}

This formulation corresponds to Eq.~(16) in \cite{charon2023shape}, representing a notational variant and a \textit{spatially homogeneous special case} of their original framework. Specifically, the explicit correspondence between our notations and theirs is as follows: our time-dependent velocity field $\vec{v}(t)$ and symmetric growth tensor $\tens{G}(t)$ directly correspond to their $v(t, \cdot)$ and $g(t, x)$, respectively. Furthermore, our strain tensor $\tens{\varepsilon}(\vec{v}(t))$ represents their linearized deformation tensor $(dv(t,x) + dv(t,x)^T)/2$, while our integral measure $dV$ and elastic energy density function $\Psi(\cdot)$ map to their spatial measure $dx$ and elastic tensor $B(x, \cdot)$.

To ensure that the minimization problem in Eq.~\eqref{eq:morpho_action} defines a valid geodesic distance $D$ and avoids trivial solutions, our model operates under several key mathematical assumptions. First, unlike the spatially varying tensor $B(x, \cdot)$ in the general model, our elastic energy density $\Psi(\cdot)$ does not explicitly depend on the spatial coordinate $x$, which inherently assumes that the material properties of the shape $\Omega(t)$ are spatially homogeneous. Second, we require $\Psi$ to be a positive semi-definite quadratic form satisfying the coercivity condition (e.g., $\Psi(S) \ge c|S|^2$). Finally, the optimal growth tensor $\tens{G}(t)$ is constrained within a specific feasible set $\mathcal{G}(\Omega(t))$. Without this restriction, one could simply choose $\tens{G}(t) = \tens{\varepsilon}(\vec{v}(t))$ to trivially reduce the elastic penalty to zero. In our work, we adopt a Gaussian growth model, which naturally satisfies this specific restriction. 

Under these conditions, the existence of the minimizer for Eq.~\eqref{eq:morpho_action} is mathematically guaranteed. The proof relies on the weak convergence of minimizing sequences and the weak lower semi-continuity of the action functional, which has been rigorously established in Appendix B of \cite{charon2023shape}.

\paragraph{1.  Positive Definiteness.}
Since $\|\mathbf{G}(t)\|_{[\Omega(t)]} \ge 0$, the integral is non-negative, and thus $d(\Omega_0, \Omega_1) \ge 0$.
If $\Omega_0 = \Omega_1$, the trivial path with $\mathbf{G}(t) = \mathbf{0}$ has zero length, so $d(\Omega_0, \Omega_0) = 0$. Conversely, if $d(\Omega_0, \Omega_1) = 0$, the optimal path must have zero length, which implies $\|\mathbf{G}(t)\|_{[\Omega(t)]} = 0$ for almost every $t$. This means $\mathbf{G}(t) = \mathbf{0}$ a.e., which in turn implies the velocity field is zero. A zero velocity field means the shape does not evolve, so $\Omega_1 = \Omega_0$.

\paragraph{2. Symmetry.} ($d(\Omega_0, \Omega_1) = d(\Omega_1, \Omega_0)$)

Given an optimal path from $\Omega_0$ to $\Omega_1$ parameterized by $t \in [0,1]$ with control $\mathbf{G}(t)$, we can construct a time-reversed path from $\Omega_1$ to $\Omega_0$ by setting the new parameter $\tau = 1-t$. The control for this reversed path is $\mathbf{G}'(\tau) = -\mathbf{G}(1-\tau)$. The length of this path is:
\begin{align*}
	L' &= \int_0^1 \|\mathbf{G}'(\tau)\|_{[\Omega(1-\tau)]} d\tau = \int_0^1 \|-\mathbf{G}(1-t)\|_{[\Omega(1-t)]} dt \\
	&= \int_0^1 \|\mathbf{G}(1-t)\|_{[\Omega(1-t)]} dt.
\end{align*}
With a change of variable $s=1-t$, we see $L'$ is equal to the original path length $L$. Since for every path from $\Omega_0$ to $\Omega_1$, there exists a path of equal length from $\Omega_1$ to $\Omega_0$, their infima must be equal.

\paragraph{3. Triangle Inequality.} ($d(\Omega_0, \Omega_2) \le d(\Omega_0, \Omega_1) + d(\Omega_1, \Omega_2)$)

This follows the standard construction for path-based metrics. For any $\epsilon > 0$, we can find a path $\mathcal{P}_{01}$ from $\Omega_0$ to $\Omega_1$ of length $L(\mathcal{P}_{01}) < d(\Omega_0, \Omega_1) + \epsilon/2$, and a path $\mathcal{P}_{12}$ from $\Omega_1$ to $\Omega_2$ of length $L(\mathcal{P}_{12}) < d(\Omega_1, \Omega_2) + \epsilon/2$.
Concatenating these paths creates a valid (though not necessarily optimal) path $\mathcal{P}_{02}$ from $\Omega_0$ to $\Omega_2$ with length $L(\mathcal{P}_{02}) = L(\mathcal{P}_{01}) + L(\mathcal{P}_{12})$.
The distance $d(\Omega_0, \Omega_2)$ is the infimum over all such paths, so it must be less than or equal to the length of this particular path:
\begin{align*}
	d(\Omega_0, \Omega_2) &\le L(\mathcal{P}_{02}) = L(\mathcal{P}_{01}) + L(\mathcal{P}_{12}) \\
	&< d(\Omega_0, \Omega_1) + d(\Omega_1, \Omega_2) + \epsilon.
\end{align*}
Since this holds for any $\epsilon > 0$, the inequality must hold without $\epsilon$. Thus, the functional defines a valid distance metric.

\section{Derivation of the Generalized High-Order PDE in Strong Form}
\label{appendix:general_strong_form}

This section explains why the inclusion of the generalized LDDMM-type
regularization transforms the equation into a system of order \(4\beta\).
The derivation below should be understood as a formal derivation of the
interior strong form. In the actual numerical formulation, we use the weak
or mixed variational form, where the corresponding high-order natural boundary
conditions are handled weakly.

The regularized weak form of the problem is to find a displacement field
\(\vec u\) such that for all admissible test functions \(\vec h\),
\begin{equation}
    \int_{\Omega}
    \tens{\sigma}(\tens{\varepsilon}(\vec u)):\nabla\vec h\,dV
    +
    \kappa
    \int_{\Omega}
    L\vec u\cdot L\vec h\,dV
    =
    \int_{\Omega}
    \tens{\sigma}(\tens G):\nabla\vec h\,dV.
    \label{eq:appendix_general_weak}
\end{equation}
Equivalently, since the stress tensors are symmetric, one may replace
\(\nabla\vec h\) by \(\tens{\varepsilon}(\vec h)\) in the stress terms. We keep
the present notation to remain consistent with the preceding derivation.

Let
\[
    A=-\alpha\Delta+\gamma I,
    \qquad
    L=A^\beta,
\]
where \(I\) denotes the identity operator acting componentwise on vector fields.
Thus,
\[
    L^\dagger L
    =
    (A^\beta)^\dagger A^\beta.
\]

We first combine the elastic and growth stress terms:
\begin{equation}
    \int_{\Omega}
    \left(
    \tens{\sigma}(\tens{\varepsilon}(\vec u))
    -
    \tens{\sigma}(\tens G)
    \right):\nabla\vec h\,dV.
\end{equation}
Applying integration by parts once gives
\begin{equation}
    -\int_{\Omega}
    \nabla\cdot
    \left(
    \tens{\sigma}(\tens{\varepsilon}(\vec u))
    -
    \tens{\sigma}(\tens G)
    \right)
    \cdot \vec h\,dV
    +
    \int_{\partial\Omega}
    \left(
    \tens{\sigma}(\tens{\varepsilon}(\vec u))
    -
    \tens{\sigma}(\tens G)
    \right)\vec n\cdot \vec h\,dS.
    \label{eq:appendix_general_elastic_boundary}
\end{equation}
For the purpose of identifying the interior strong form, we focus on the volume
term. The boundary term corresponds to the usual traction-type natural boundary
condition.

The regularization term is
\begin{equation}
    \kappa
    \int_{\Omega}
    L\vec u\cdot L\vec h\,dV
    =
    \kappa
    \int_{\Omega}
    (A^\beta\vec u)\cdot(A^\beta\vec h)\,dV.
\end{equation}
Formally, if the boundary conditions are chosen so that all boundary terms
generated by repeated integration by parts vanish, then \(A\) is self-adjoint
with respect to the \(L^2\)-inner product. In that case,
\[
    (A^\beta)^\dagger=A^\beta,
    \qquad
    L^\dagger L=A^\beta A^\beta=A^{2\beta}.
\]
Therefore,
\begin{align}
    \kappa
    \int_{\Omega}
    (A^\beta\vec u)\cdot(A^\beta\vec h)\,dV
    &=
    \kappa
    \int_{\Omega}
    (A^{2\beta}\vec u)\cdot\vec h\,dV
    \nonumber\\
    &=
    \kappa
    \int_{\Omega}
    \left[
    (-\alpha\Delta+\gamma I)^{2\beta}\vec u
    \right]\cdot\vec h\,dV.
    \label{eq:appendix_general_reg}
\end{align}
This is the formal interior contribution of the regularization term.

\paragraph{Conclusion: Interior Strong Form.}

Combining the volume contributions from the elastic-growth part and the
regularization part, we obtain
\begin{equation}
    \int_{\Omega}
    \left(
    \kappa(-\alpha\Delta+\gamma I)^{2\beta}\vec u
    -
    \nabla\cdot\tens{\sigma}(\tens{\varepsilon}(\vec u))
    +
    \nabla\cdot\tens{\sigma}(\tens G)
    \right)\cdot\vec h\,dV
    =
    0.
\end{equation}
Since this identity holds for arbitrary test functions supported in the interior
of \(\Omega\), the interior strong form is
\begin{equation}
    \kappa(-\alpha\Delta+\gamma I)^{2\beta}\vec u
    -
    \nabla\cdot\tens{\sigma}(\tens{\varepsilon}(\vec u))
    +
    \nabla\cdot\tens{\sigma}(\tens G)
    =
    \vec 0
    \qquad
    \text{in } \Omega.
    \label{eq:appendix_general_strong_form}
\end{equation}

This equation is a differential system of order \(4\beta\), since
\((-\alpha\Delta+\gamma I)^{2\beta}\) contains derivatives up to order \(4\beta\).
For the specific case \(\beta=2\) used in the simulations, the strong-form
regularization operator is
\[
    (-\alpha\Delta+\gamma I)^4,
\]
which is an eighth-order differential operator.

\paragraph{Remark on Boundary Conditions.}

The above derivation identifies the interior strong form. The complete strong
boundary value problem would also contain boundary conditions arising from both
the physical elasticity model and the high-order regularization term. The
physical boundary conditions are imposed as
\[
    \vec u=\vec 0
    \qquad
    \text{on } \Gamma_D,
\]
on the fixed boundary, while the remaining part \(\Gamma_N\) carries the
Neumann-type traction condition.

The high-order regularization term may induce additional natural boundary
conditions if the problem is written entirely in strong form. These conditions
involve traces of auxiliary high-order quantities such as \(A^k\vec u\) and their
normal derivatives. In this work, however, the numerical formulation is based on
the weak and mixed variational forms. Therefore, these additional boundary
conditions are not imposed separately; they are encoded weakly by the choice of
admissible spaces, test spaces, and the regularized bilinear form.

\bibliographystyle{plain}
\bibliography{references}

\end{document}